\documentclass[11pt]{amsart}

\usepackage{amsmath}
\usepackage{amsfonts}
\usepackage{amssymb}
\usepackage{latexsym}

\usepackage[cp866]{inputenc}

 \advance\hoffset by -0.0 truecm

\newtheorem{theorem}{Theorem}[section]
\newtheorem{corollary}[subsection]{Corollary}
\newtheorem{lemma}[subsection]{Lemma}
\newtheorem{proposition}[subsection]{Proposition}
\theoremstyle{definition}
\newtheorem{definition}[subsection]{Definition}
\theoremstyle{remark}
\newtheorem{remark}[theorem]{Remark}
\numberwithin{equation}{section}
 
 \DeclareMathOperator{\capac}{cap}

\DeclareMathOperator{\dist}{dist}
\DeclareMathOperator{\supp}{supp}
\DeclareMathOperator{\loc}{loc}
\DeclareMathOperator{\card}{card}

\DeclareMathOperator{\spn}{span}
\DeclareMathOperator{\diam}{diam}
\DeclareMathOperator{\diverg}{div}
\DeclareMathOperator{\mes}{mes}
\DeclareMathOperator{\ACL}{ACL}
\DeclareMathOperator{\dv}{div}
\DeclareMathOperator{\End}{End}
\DeclareMathOperator{\Id}{Id}
\DeclareMathOperator{\sign}{sign}
\newcommand{\esssup}{\operatornamewithlimits{ess\,sup}}
\newcommand{\Om}{\Omega}

\newcommand{\Real}{\mathbb{R}}

\newcommand{\G}{\mathbb{G}}
\newcommand{\Heis}{\mathbb{H}}

\newcommand{\CG}{\mathcal{G}}

\begin{document}

\title[Value Distribution]
{On Value Distributions for Quasimeromorphic Mappings on $\mathbb
H$-type Carnot Groups}

\author{Irina Markina,\ Sergey Vodopyanov}

\address{Departamento de Matem\'atica,\newline
Universidad T\'ecnica Federico Santa Mar\'{\i}a,\newline av.
Espa\~{n}a 1680, Valpara\'{\i}so, Chile}

\email{irina.markina@usm.cl}

\address{Sobolev Institute of Mathematics,\newline pr. Koptyuga 4,
Novosibirsk, Russia}

\email{vodopis@math.nsc.ru}

\thanks{This research  was carried out with the partial support
of the Russian Foundation for Basic Research, \#03--01--00899, and
the Projects FONDECYT (Chile), \#1030373, \#1040333}




\begin{abstract}
In the present paper we define quasimeromorphic mappings on
homogeneous groups and study their properties. We prove an
analogue of results of L.~Ahlfors, R.~Nevanlinna and S.~Rickman,
concerning the value distribution for quasimeromorphic mappings on
$\mathbb H$-type Carnot groups for parabolic and hyperbolic
 type domains.
\end{abstract}

\maketitle

\section*{Introduction}

The classical  value distribution theory for analytic functions
$w(z)$ studies the system of sets $z_a$ of a domain $G_z$ where
the function $w(z)$ takes the value $w=a$ for an arbitrary $a$. A
central result in the distribution theory is the Picard theorem,
stating that a meromorphic function in the punctured plane assumes
all except for at most two values $a_1,a_2$, $a_1\neq a_2$
infinitely often. In an equivalent way, we can say that an
analytic function $w(z):\mathbb R^2\to\mathbb
R^2\setminus\{a_1,a_2\}$ must be constant if $a_1$ and $a_2$ are
distinct points in $\mathbb R^2$. We mention results of
J.~Hadamard, E.~Borel, G.~Julia, A.~Beurling, L.~V.~Ahlfors and
others~\cite{Ahl,AB,Bor,Jul} in general value distribution theory,
going far beyond Picard-type theorems. Nevertheless, in those
extensions and deep generalizations the nature of conformal
mappings was not actively involved. New ideas of the function
theory and potential theory point of view were incoming by
R.~Nevanlinna~\cite{Nev1,Nev2}. The most important achievements of
the Nevanlinna theory were not only analytic deep results, but its
geometric aspects and relations with Riemannian surfaces of
analytic functions. Such principal notions, as a characteristic
function, a defect function, a branching index connect the
asymptotic behavior of an analytic function $w(z)$ with properties
of the Riemannian surface which is the conformal image of the
domain of $w(z)$.

A natural generalization of an analytic function of one complex
variable to the Euclidean space of the dimension $n>2$ was firstly
introduced and studied by Yu.~G.~Reshetnyak in
1966---1968~\cite{Resh2,Resh3,Resh1}. In some sense this is a
quasiconformal mapping admitting branch points. Such mappings were
called in Russian school the {\it mappings with bounded
distortion}. The main contribution of Yu.~G.~Reshetnyak to the
foundation of this theory is a discovery of a connection between
mappings with bounded distortion and non-linear partial
differential equations.  Yu.~G.~Reshetnyak has proved also that an
analytic definition of mappings with bounded distortion implies
the topological properties: the continuity, the openness, and the
discreteness. Later these mappings, under the name {\it
quasiregular mappings}, were investigated intensively by
O.~Martio, S.~Rickman, J.~V\"ais\"al\"a, F.~W.~Gehring,
M.~Vuorinen, B.~Bojarski, T.~Iwaniec and
others~\cite{BI,BI1,Geh,MRV2,MRV1,Rick1,Sr,Vuor}. Briefly a
quasiregular mapping can be defined as an appropriate Sobolev
mapping with nonnegative Jacobian and such that an infinitesimal
ball is transformed into infinitesimal ellipsoid with bounded
ratio of the largest and the smallest semi-exes.

The Picard theorem is true for quasiregular mappings in $\mathbb
R^2$. In fact, an arbitrary quasiregular mapping $f$ in $\mathbb
R^2$ has a representation $f=g\circ h$ where $h:\mathbb
R^2\to\mathbb R^2$ is a quasiconformal mapping and $g$ is an
analytic function of $\mathbb R^2$ omitting two points~\cite{LV2}.
In 1967 V.~A.~Zorich~\cite{Zo} asked whether a Picard-type theorem
exists for quasiregular mappings in higher dimensions. S.~Rickman
has given a complete answer to the question and developed the
value distribution theory for quasimeromorphic mappings in
$\mathbb R^n$, $n>2$, based on the potential theory, metric and
topological properties of quasiregular
mappings~\cite{Rick3,Rick4,Rick5,Rick6,Rick1}. Quasimeromorphic
mappings $f:\mathbb R^n\to\overline{\mathbb R}^n$ generalize
quasiregular mappings in the same way as meromorphic mappings do
the analytic functions.

A stratified nilpotent group (of which $\Real^n$ and the
Heisenberg group are the simplest examples) is a Lie group
equipped with an appropriate family of dilations. Thus, this group
forms a natural habitat for extensions of many of the objects
studied in the Euclidean space. The fundamental role of such
groups in analysis was noted by E.~M.~Stein~\cite{St1,St2}. There
has been since a wide development in the analysis of the so-called
stratified nilpotent Lie groups, nowadays, also known as Carnot
groups. The theory of quasiconformal and quasiregular mappings on
Carnot groups is presented in the
works~\cite{Dair,HH,HK,KR,Mar1,Rick2,Vod2,Vod10,Vod7}.

In the present paper we define quasimeromorphic mappings on Carnot
groups and study their properties. The main difference with the
Euclidean definition of a quasimeromorphic mapping is an absence
of inversions on general Carnot groups. Therefore, we are not able
to use neither a stereographic projection nor conformal metric on
Carnot groups. Nevertheless, the definition of quasimeromorphic
mappings on Carnot groups we give, allows us not only to obtain
the analogues of their Euclidean properties but to adopt also the
ideas and methods of the value distribution theory for
quasimeromorphic mappings in the Euclidean spaces, developed by
S.~Rickman. The main difference with respect to Rickman's approach
is that we do not use the inversion as a conformal mapping defined
on the one-point compactification of a Carnot group. We present
some results concerning the value distribution of
$K$-quasimeromorphic mappings on $\mathbb H$(eisenberg)-type
Carnot groups in a domain with one boundary point and for
$K$-quasimeromorphic mappings defined on the unit ball.

The paper is organized as follows. In Section~1 we give the
necessary definitions. Section~2 is devoted to the properties of
quasimeromorphic mappings and capacity estimates on an arbitrary
Carnot group. In Section~3 we consider module inequalities playing
fundamental role in the proofs of the main theorems. Section~4 is
dedicated to the relationships between the module of a family of
curves, a counting function, and averages of the counting function
over spheres. In Section~5 we state the first main theorem and
prove auxiliary lemmas. Sections~6 and~7 are devoted to proofs of
the first and the second principal theorems respectively. The
theorems are stated and showed for the $\mathbb H$-type Carnot
groups.

 It is well-known, that S.~Rickman employed a special family
of curves in order to find a method for estimating their modules.
A suitable counterpart of such families on Carnot groups exists in
frame of "polar coordinates" in the $\mathbb H$-type Carnot
groups. The key property is that the radial curves have the finite
length in Carnot-Carath\'eodory metric. It allows us to involve
the classical module methods~\cite{KR1,KR,Vas}. In~\cite{BT},
"polarizable" Carnot groups were introduced. These groups admit
the analogue of polar coordinates. Unfortunately, nowadays, it is
unknown an example of polarizable Carnot group, which is not of
$\mathbb H$-type. By the way, there are no non-trivial examples of
quasiregular mappings on an arbitrary Carnot groups. Nevertheless,
if the theory of quasiregular mappings is not degenerate on the
polarizable Carnot groups, then our results are true also for this
setting.

Now we state the principal result of our work.

\begin{theorem}\label{45}
Let $\mathbb G$ be a $\mathbb H$-type Carnot group, $f:\ \mathbb
G\to\overline{\mathbb G}$ be a nonconstant $K$-quasimeromorphic
mapping. Then there exists a set $E\subset[1,\infty[$ and a
constant $C(Q,K)<\infty$ such that
\begin{equation}\label{222}
\lim\limits_{r\to \infty}\sup\limits_{r\notin E}\sum\limits_{j=0
}^{q}\Big(1-\frac{n(r,a_j)}{\nu(r,1)}\Big)_+\leq
C(Q,K)\quad\text{with}\quad \int\limits_{E}\frac{dr}{r}<\infty,
\end{equation}
whenever $a_0,a_1,\ldots,a_q$ are distinct points in
$\overline{\mathbb G}$.
\end{theorem}

The definitions of the counting function $n(r,a_j)$ and the
average $\nu(r,1)$ see in Section~4.

We would like to call the attention on the difference between our
assertion and  Rickman's one~(see for
instance~\cite[p.~80]{Rick1}). S.~Rickman employed a version
of~\eqref{222} which is conformally invariant and used essentially
this property in his proof. R.~Nevanlinna pointed out that
averages of the counting function with respect to distinct
measures can find different applications and physic-geometrical
meaning (see also O.~Frostman~\cite{Fr1,Fr2}). For this reason,
possessing only limited geometrical and analytical tools, we deal
with expression~\eqref{222} which is not conformally invariant but
still carries an information sufficient to effectively control the
distribution of values of a quasimeromorphic mapping. As a
corollary of our main result we get the Picard theorem.

\begin{theorem}
Let $\mathbb G$ be a $\mathbb H$-type Carnot group. For each
$K\geq 1$, there exists a constant $q(\mathbb G,K)$ such that
every  $K$-quasiregular mapping $f:\mathbb G\to\mathbb
G\setminus\{a_1,\ldots,a_q\}$, where $q\geq q(\mathbb G,K)$ and
$a_1,\ldots,a_q$ are distinct, is constant.
\end{theorem}

Another way of proving this assertion can be found in~\cite{Vod7}.

The next theorem is stated for $K$-quasimeromorphic mappings in
the unit ball $B(0,1)$. The proof of the statement essentially
uses the method developed for Theorem~\ref{45}.

\begin{theorem}\label{46}
Let $\mathbb G$ be a $\mathbb H$-type Carnot group,
$f:B(0,1)\to\overline{\mathbb G}$ be a nonconstant
$K$-quasimeromorphic mapping such that
\begin{equation*}
\limsup\limits_{r\to 1}(1-r)A(r)^{\frac{1}{Q-1}}=\infty.
\end{equation*}
Then there exists a set $E\subset(0,1)$ satisfying
\begin{equation*}
\liminf\limits_{r\to 1}\frac{\mes_1(E\cap [r,1))}{(1-r)}=0,
\end{equation*}
and a constant $C(Q,K)<\infty$ such that
\begin{equation*}
\lim\limits_{r\to 1}\sup\limits_{r\notin E}\sum\limits_{j=0
}^{q}\Big(1-\frac{n(r,a_j)}{\nu(r,1)}\Big)_+\leq C(Q,K),
\end{equation*}
whenever $a_0,a_1,\ldots,a_q$ are distinct points in
$\overline{\mathbb G}$.
\end{theorem}

\section{Notations and definitions}

The Carnot group is a connected and simply connected nilpotent Lie
group $\G$ whose Lie algebra $\CG$ decomposes into the direct sum
of vector subspaces $V_1\oplus V_2\oplus\ldots\oplus V_m$
satisfying the following relations:
$$[V_1,V_k] = V_{k+1},\qquad 1\leq k<m,\quad\quad [V_1,V_m]
=\{0\}.$$

We identify the Lie algebra $\CG$ with a space of left-invariant
vector fields. Let $X_{11},\ldots,X_{1n_1}$ be a basis of $V_1$,
$n_1=\dim V_1$, and $\langle\cdot,\cdot\rangle_0$ be a
left-invariant Riemannian metric on $V_1$ such that $$ \langle
X_{1i},X_{1j}\rangle_0=\left\lbrace \begin{array}{c ll} 1 &
\text{if}\  & i=j, \\ 0 & \text{if} & i\neq j.
\end{array}\right.
$$ Then, $V_1$ determines a subbundle $HT$ of the tangent bundle
$T\G$ with fibers $$
HT_q=\spn\,\{X_{11}(q),\ldots,X_{1n_1}(q)\},\quad q\in\G. $$ We
call $HT$ the {\it horizontal tangent bundle} of $\G$ with $HT_q$
as the {\it horizontal tangent space} at $q\in\G$. Respectively,
the vector fields $X_{1j}$, $j=1,\ldots,n_1$, are said to be  {\it
horizontal vector fields}.

Next, we extend $X_{11},\ldots,X_{1n_1}$ to a basis
$$
X_{11},\ldots,X_{1n_1},X_{21},\ldots,X_{2n_2},\ldots,X_{m1},\ldots,X_{mn_m}
$$ of $\CG$. Here, each vector field $X_{ij}$, $2\leq i\leq m$,
$1\leq j\leq n_i=\dim V_i$, is a commutator $$
X_{ij}=[\ldots[[X_{1k_1},X_{1k_2}],X_{1k_3}],\ldots ,X_{1k_i}] $$
of the length $i-1$ of  basic vector fields of the chosen basis of
$V_1$.

It is known  (see, for instance,~\cite{FS}) that the exponential
map $\exp :\mathcal G\to \G$ from the Lie algebra $\mathcal G$
into the Lie group $\G$ is a global diffeomorphism. We can
identify the points $q\in \G$ with the points $x\in \Real^N$,
$N=\sum\limits_{i=1}^{m}\dim V_i$, by means of the mapping
$q=\exp\bigr(\sum\limits_{i,j}x_{ij}X_{ij}\bigl)$. The collection
$\{x_{ij}\}$ is called the {\it normal coordinates} of $q\in\G$. The
number $N=\sum\limits_{i=1}^{m}\dim V_i$ is the topological
dimension of the Carnot group. The bi-invariant Haar measure on
$\G$ is denoted by $dx$; this is the push-forward of the Lebesgue
measure in $\Real^N$ under the exponential map. {\it The family of
dilations} $\{\delta_{\lambda}(x):\lambda>0\}$ on the Carnot group
is defined as $$
\delta_{\lambda}x=\delta_{\lambda}(x_{ij})=(\lambda
x_{1},\lambda^2x_{2},\ldots, \lambda^mx_{m}),$$ where
$x_i=(x_{i1},\ldots,x_{in_i})\in V_i$. Moreover,
$d(\delta_{\lambda}x)=\lambda^Qdx$ and the quantity
$Q=\sum\limits_{i=1}^{m}i\dim V_i$ is called  {\it the homogeneous
dimension} of $\G$.

{\sc Example 1}. The Euclidean space $\Real^n$ with the standard
structure  exemplifies an~Abelian group: the exponential map is
the identical mapping and the vector fields
$X_i=\frac{\partial}{\partial x_i}$, $i=1,\ldots, n$, have trivial
commutators only and constitute a~basis for the corresponding Lie
algebra.

{\sc Example 2}. The simplest example of a non-abelian Carnot
group is the Heisenberg group~$\Heis^n$. The non-commutative
multiplication is defined as $$
pq=(x,y,t)(x^{\prime},y^{\prime},t^{\prime})=
(x+x^{\prime},y+y^{\prime},t+t^{\prime}-2xy^{\prime}+2yx^{\prime}),
$$ where $x,x^{\prime},y,y^{\prime}\in \Real^n$,
$t,t^{\prime}\in\Real$. Left translation $L_p(\cdot)$ is defined
as $L_p(q)=pq$. The left-invariant vector fields
$$ X_i =\frac{\partial}{\partial x_i}+2y_i\frac{\partial}{\partial
t},\quad   Y_i =\frac{\partial}{\partial
y_i}-2x_i\frac{\partial}{\partial t}, \quad i=1,\ldots, n, \quad
\quad T =\frac{\partial}{\partial t},
$$
constitute the basis of the Lie algebra of the Heisenberg group.
All non-trivial relations are only of the form $[X_i,Y_i]=-4T$,
$i=1,\ldots, n$, and all other commutators vanish. Thus, the
Heisenberg algebra has the dimension $2n+1$ and splits into the
direct sum $\CG=V_1\oplus V_2$. The vector space $V_1$ is
generated by the vector fields $X_i$, $Y_i$, $i=1,\ldots, n$, and
the space $V_2$ is the one-dimensional center which is spanned by
the vector field $T$. More information see~\cite{K1,K2}.

{\sc Example 3}. A Carnot group is said to be of $\mathbb H$-type
if the Lie algebra $\mathcal G=V_1\oplus V_2$ is two-step and if
the inner product $\langle\cdot,\cdot\rangle_0$ in $V_1$ can be
extended to an inner product $\langle\cdot,\cdot\rangle$ in all of
$\mathcal G$ so that the linear map $J:V_2\to \End(V_1)$ defined
by $\langle J_ZU,V\rangle=\langle Z,[U,V]\rangle$ satisfies
$J^2_Z=-\langle Z,Z\rangle \Id$ for all $Z\in V_2$. For the moment
we introduce the notation $\|Z\|^2=\langle Z,Z\rangle$. Then
$\|J_ZV\|=\|Z\|\cdot\|V\|$ and $\langle V,J_ZV\rangle=0$ for all
$V\in V_1$ and $Z\in V_2$. More details and information see
in~\cite{Cow1,Cow2,Kap1,Kap2,Rei}.

A~homogeneous norm on $\mathbb G$ is, by definition, a continuous
function $|\cdot|$ on~$\G$ which is smooth on $\G\setminus\{0\}$
and such that $|x|=|x^{-1}|$, $|\delta_{\lambda}(x)|=\lambda|x|$,
and $|x|=0$ if and only if $x=0$. All homogeneous norms are
equivalent. We choose one of them that admits an analogue of polar
coordinates on $\mathbb H$-type Carnot groups (see~\cite{BHT}).
This norm $|\cdot|=u_2^{1/(1-Q)}$ is associated to Folland's
singular solution $u_2$ for the sub-Laplacian
$\Delta_0=\sum_{j=1}^{n_1}X_{1j}^2$ at $0\in \mathbb G$. Another
advantage of this norm is that it gives the exact value for the
$Q$-capacity of spherical ring domains. The norm $|\cdot|$ defines
a pseudo-distance: $d(x,y)=|x^{-1}y|$ satisfying the generalized
triangle inequality $d(x,y)\leq\varpi(d(x,z)+d(z,y))$ with a
positive constant $\varpi$. By $B(x,r)$ we denote an open ball of
radius $r>0$ centered at $x$ in the metric~$d$. Note that
$B(x,r)=\{y\in\G:d(x,y)<r\}$ is the left translation of the ball
$B(0,r)$ by~$x$, which is the image of the "unit ball" $B(0,1)$
under $\delta_r$. By $\mes(E)$ we denote the measure of the set
$E$. Our normalizing condition is such that the balls of radius
one have measure one: $\mes(B(0,1))=\int_{B(0,1)}dx=1$. We have
$\mes(B(0,r))=r^Q$ because the Jacobian of the dilation $\delta_r$
is~$r^Q$.

A continuous map $\gamma: I\to \mathbb G$ is called a curve. Here
$I$ is a (possibly unbounded) interval in $\mathbb R$. If
$I=[a,b]$ then we say that $\gamma: [a,b]\to \mathbb G$ is a
closed curve. A closed curve $\gamma: [a,b]\to \mathbb G$ is
rectifiable if
$$\sup\Big\{\sum\limits_{k=1}^{p-1}d\big(\gamma(t_k),\gamma(t_{k+1})\big)\Big\}<\infty,$$
where the supremum ranges over all partitions $a=t_1< t_2<\ldots<
t_p=b$ of the segment $[a,b]$. P.~Pansu proved in~\cite{Pansu1}
that any rectifiable curve is differentiable almost everywhere in
$[a,b]$ in the   Riemannian sense and there exist measurable
functions $a_{j}(s)$, $s\in [a,b]$, such that
$$\dot\gamma(s)=\sum\limits_{j=1}^{n_1}a_j(s)X_{1j}(\gamma(s))\quad
\text{and}\quad
d\big(\gamma(s+\tau),\gamma(s)\exp(\dot\gamma(s)\tau)\big)=o(\tau)\
\text{as}\ \tau\to 0$$ for almost all $s\in (a,b)$.

An absolutely continuous curve is a continuous map $\gamma: I\to
\mathbb G$ satisfying the following property: for an arbitrary
number $\varepsilon>0$, there exists $\delta>0$ such that for
arbitrary disjoint collections of segments
$(\alpha_i,\beta_i)\subset I$ with
$\sum_i(\beta_i-\alpha_i)\leq\delta$ we have
$\sum_id(\gamma(\beta_i),\gamma(\alpha_i))\leq\varepsilon$.  A
closed absolutely continuous curve  $\gamma: [a,b]\to \mathbb G$
is always rectifiable. Its length $l(\gamma)$ can be calculated by
the formula
$$
l(\gamma)= \int\limits_{a}^b
\langle\dot\gamma(s),\dot\gamma(s)\rangle_0^{1/2}ds =
\int\limits_{a}^b \biggl(\sum\limits_{j=1}^{n_1} |a_j(s)|^2
\biggr)^{1/2}\,ds
$$
where $\langle\cdot,\cdot\rangle_0$ is the left invariant
Riemannian metric on $V_1$.

A result of~\cite{Ch} implies that one can connect two arbitrary
points $x$, $y\in\G$ by a rectifiable curve. The
Carnot-Carath\'eodory distance $d_c(x,y)$ is the infimum of the
lengths over all rectifiable curves with endpoints $x$ and $y\in
\mathbb G$. Since $\langle\cdot,\cdot\rangle_0$ is left-invariant,
the Carnot-Carath\'eodory metric is also left-invariant. The
metric $d_c(x,y)$ is finite since the points $x,y\in\G$ can be
joined by a rectifiable curve with endpoints $x,y$. The Hausdorff
dimension of the metric space $(\G,d_c)$ coincides with the
homogeneous dimension $Q$ of the group~$\mathbb G$. More
information see in~\cite{Mitch,Pansu1,Str}.

The Sobolev space $W^1_p(\Om)$ ($L^1_p(\Om)$), $1\leq p<\infty$,
consists of locally summable functions $u:\Om\to\mathbb R$,
$\Omega\subset\mathbb G$, having distributional derivatives
$X_{1j}u$ along the vector fields $X_{1j}$:
$$\int\limits_{\Omega}X_{1j}u\varphi\,dx=-\int\limits_{\Omega}uX_{1j}\varphi\,dx,\qquad
j=1,\ldots,n_1,$$ for any test function $\varphi\in C_0^{\infty}$,
and the finite norm
$$
\|u\mid W^1_p(\Om)\|=\Bigl(\int\limits_{\Om}|u|^p\,dx\Bigr)^{1/p}+
\Bigl(\int\limits_{\Om}|\nabla_{0}u|_0^p\,dx\Bigr)^{1/p}
$$
(semi-norm
$$
\|u\mid L^1_p(\Om)\|=
\Bigl(\int\limits_{\Om}|\nabla_{0}u|_0^p\,dx\Bigr)^{1/p}\Bigr).
$$
Here $\nabla_{0}u=(X_{11}u,\ldots,X_{1n_1}u)$ is the {\it
subgradient} of $u$ and
$|\nabla_{0}u|_0=\langle\nabla_{0}u,\nabla_{0}u\rangle_0$. We say,
that $u$ belongs to $W^1_{p,\loc}(\Omega)$ if for an arbitrary
bounded domain $U$, $\overline U\subset\Om$, the function $u$
belongs to $W^1_p(U)$. Henceforth, for a bounded domain
$U\subset\Omega$ whose closure $\overline U$ belongs to $\Omega$,
we write $U\Subset\Omega$ and say that $U$ is a compact domain in
$\Omega$.

\begin{definition}[\cite{Resh4,Vod6,Vod8}] \label{0}
Suppose that $(\mathbb X, r)$ is a complete metric space, $r$ is a
metric on~$\mathbb X$, and $\Omega $ is a domain on a Carnot
group~$\mathbb G$. We say that a mapping $f:\Omega\to\mathbb X$
{\it belongs to Sobolev class\/}\ $W_{p, \loc}^1(\Omega; \mathbb
X)$ if the following conditions~hold.
\begin{itemize}
\item[(A)] {For each $z\in\mathbb X$, the function
 $[f]_z:x\in\Omega\mapsto r(f(x), z)$
 belongs to the class $W_{p,\loc}^1(\Omega)$.}
\item[(B)] {The  family of functions
$\big(\nabla_{0}[f]_z\big)_{z\in\mathbb X}$ has a dominant
belonging to~$L_{p,\loc}(\Omega)$, i.e., there is a function $g\in
L_{p,\loc}(\Omega)$ independent of~$z$ and such that
$\big\vert\nabla_{0}[f]_z(x)\big\vert_0\le g(x)$ for almost all
$x\in\Omega $.}
\end{itemize}
\end{definition}

\begin{definition}\label{1} A function $u:\Omega\to\mathbb R$, $\Omega\subset\mathbb G$,
is said to be {\it absolutely continuous on lines} $(u\in
\ACL(\Omega))$ if for any domain $U\Subset\Omega$, and any
fibration $\mathcal X_j$ defined by the left-invariant vector
fields $X_{1j}$, $j=1,\ldots,n_1$, the function $u$ is absolutely
continuous on $\gamma\cap U$ with respect to the $\mathcal
H^1$-Hausdorff measure for $d\gamma$-almost all curves
$\gamma\in\mathcal X_j$. (Recall that the measure $d\gamma$ on
$\mathcal X_j$ equals the inner product $i(X_j)$ of the vector
field $X_j$ by the bi-invariant volume form $dx$.)
\end{definition}

For a function $u\in \ACL(\Omega)$, the derivatives $X_{1j}u$
along the vector fields $X_{1j}$, $j=1,\ldots,n_1$, exist almost
everywhere in $\Omega$. It is known that a function
$u:\Om\to\mathbb R$ belongs to $W^1_p(\Om)$ ($L^1_p(\Om)$), $1\leq
p<\infty$, if and only if it can be modified on a set of measure
zero by such a way that $u\in L_p(\Om)$ ($u$ is locally
$p$-summable), $u\in \ACL(\Omega)$, and $X_{1j}u\in L_p(\Om)$
hold, $j=1,\ldots,n_1$. The reader can find More information on
$\ACL$-functions in~\cite{KR,UV2,Vod6}.

\begin{proposition}[\cite{Vod6,Vod8}]
A mapping $f:\Om\to\mathbb G$, $\Om\subset\mathbb G$, belongs to
the Sobolev class $W^1_{p,\loc}(\Om)$, $1\leq p<\infty$, if and
only if it can be modified on a set of measure zero by such a way
that
\begin{itemize} \item[1)]{$|f(x)|\in L_{p,\loc}(\Omega)$;}
\item[2)] {the coordinate functions $f_{ij}$ belong to
$\ACL(\Omega)$ for all $i$ and $j$;} \item[3)] {$f_{1j}\in
W^1_{p,\loc}(\Om)$ for $1\leq j \leq n_1$;} \item[4)] {the vector
$$ X_{1k}(f(x))=\sum\limits_{1\leq l\leq m, 1\leq\omega\leq
n_l}X_{1k}(f_{l\omega}(x)) \frac{\partial}{\partial x_{l\omega}}
$$ belongs to $HT_{f(x)}$ for almost all $x\in\Om$ and all
$k=1,\ldots,n_1$.}
\end{itemize}
\end{proposition}

In~\cite{HajK,Vod6,Vod7}, one can find various definitions of the
Sobolev space on Carnot groups and their correlations. The matrix
$X_{1k}f=(X_{1k}f_{1j})_{k,j=1,\ldots,n_1}$ defines a linear
operator $D_Hf:\ V_1\to V_1$~\cite{Pansu1,Pansu2} which is called
a {\it formal horizontal differential}. A norm of the operator
$D_Hf$ is defined by $$ |D_Hf(x)|=\sup\limits_{\xi\in
V_1,|\xi|_0=1} |D_Hf(x)(\xi)|_0.
$$ The norm $|D_Hf|$ is equivalent to $|\nabla_{0}f|_0=
\Bigl(\sum\limits_{i=1}^{n_1}|X_{1i}f|_0^2\Bigr)^{\frac{1}{2}}$.
It has been proved in~\cite{UV2,Vod6} that the formal horizontal
differential $D_Hf$ generates a homomorphism $Df:\mathcal G\to
\mathcal G$ of Lie algebras which is called a {\it formal
differential}. The determinant of the matrix $Df(x)$ is denoted by
$J(x,f)$ and called a {\it $($formal$)$ Jacobian}.

A continuous mapping $f:\Omega \to\G$, $\Omega\subset\G$, is {\it
open} if the image of an open set is open and {\it discrete} if
the pre-image $f^{-1}(y)$ of each point $y\in f(\Omega)$ consists
of isolated points. We say that $f$ is sense-preserving if a
topological degree $\mu(y,f,U)$ is strictly positive for all
domains $U\Subset\Omega$ and $y\in f(U)\setminus f(\partial U)$.
The precise definition of the topological degree see in
Subsection~\ref{degree}.

\begin{definition}\label{3}
Let $\Omega$ be a domain on the group $\mathbb G$. A mapping
$f:\Omega\to{\mathbb G}$ is said to be a {\it
quasiregular mapping} if
\begin{itemize}
\item[1)] {$f$ is continuous open discrete and  sense-preserving
;} \item[2)] {$f$ belongs to $W^1_{Q,\loc}(\Omega)$;} \item[3)]
{the formal horizontal differential $D_Hf$ satisfies the condition
\begin{equation}\label{21} \max\limits_{|\xi|_0=1,\xi\in
V_1}|D_Hf(x)(\xi)|_0\leq K\min\limits_{|\xi|_0=1,\xi\in
V_1}|D_Hf(x)(\xi)|_0
\end{equation} for almost all $x\in \Omega$.}
\end{itemize}
\end{definition}

It is known \cite{Vod6} that the pointwise inequality~(\ref{21})
is equivalent to the following one: {\it the formal horizontal
differential $D_Hf$ satisfies the condition
\begin{equation}\label{22} |D_Hf(x)|^{Q}\leq
K^{\prime}J(x,f)\end{equation} for almost all $x\in \Omega$ where
$K^{\prime}$ depends on $K$.} The smallest constant $K^{\prime}$
in inequality (\ref{22}) is called the {\it outer distortion} and
denoted by $K_O(f)$. It is not hard to see that for a quasiregular
mapping the inequality
\begin{equation}\label{id}
0\leq J(x,f)\leq  K^{''}\min\limits_{|\xi|_0=1,\xi\in
V_1}|D_Hf(x)(\xi)|_0^Q
\end{equation} also holds
for almost all $x\in \Omega$ where $K^{''}$ depends on $K$. The
smallest constant $K^{''}$ in inequality~(\ref{id}) is called the
{\it inner distortion} and denoted by $K_I(f)$.

It is established in~\cite{Vod12,Vod11} that the conditions 2 and
3 of Definition~\ref{3} provide for a non-constant mapping on a
two-step Carnot group to be continuous open discrete and
sense-preserving if there exists a singular solution $w \in
W^1_{\infty,\loc}(\mathbb G\setminus 0)$ to the equation
$\dv(|\nabla_0w|^{Q-2}\nabla_0w)=0$. In~\cite{Dair}, the same
result is proved under stronger assumption that the solution $w$
belongs to $C^1$. Such a singular solution exists on the $\mathbb
H$-type Carnot groups~\cite{HH}. By another words, on the $\mathbb
H$-type Carnot groups a mapping with bounded distortion (that is a
mapping satisfying conditions 2 and 3 of Definition~\ref{3}) is
also a quasiregular one. As soon as on Carnot groups, there is no
a complete counterpart of the Euclidean theory of mappings with
bounded distortion we will distinguish mappings with bounded
distortion and quasiregular mappings.

\begin{definition} A continuous mapping $f:\Omega\to\mathbb G$ is $\mathcal
P$-differentiable at $x\in\Omega$ if the family of maps
$f_t=\delta_{1/t}(f(x)^{-1}f(x\delta_ty))$ converges locally
uniformly to an automorphism of $\mathbb G$ as $t\to 0$.
\end{definition}

In the following theorem we formulate analytic properties of
quasiregular mappings \cite{Vod10,Vod6,Vod7,Vod9}. In the
statement of the theorem we use notions of a topological degree
$\mu(y,f,D)$ of the mapping $f$ and a multiplicity function
$N(y,f,A)=\card\{x\in f^{-1}(y)\cap A\}$ (see the precise
definitions in Subsection~\ref{degree}).

\begin{theorem}\label{226} Let $f: \Omega\to\mathbb G$,
$\Omega\subset\mathbb G$, be a quasiregular mapping. Then it
possesses the following properties:
\begin{itemize}
\item[1)]{$f$ is $\mathcal P$-differentiable almost everywhere in
$\Omega$;} \item[2)]{$\mathcal N$-property: if $\mes(A)=0$ then
$\mes(f(A))=0$;} \item[3)] {$\mathcal N^{-1}$-property: if
$\mes(A)=0$ then $\mes(f^{-1}(A))=0$;}
\item[4)]{$\mes(B_f)=\mes(f(B_f))=0$;}
\item[5)] {$J(x,f)>0$ almost everywhere in $\Om$;}
\item[6)]{for every compact
domain $D\Subset\Omega$ such that $\mes(f(\partial D))=0$ (every
measurable set $A\subset \Omega$) and every measurable function
$u$, the function $y\mapsto u(y)\mu(y,f,D)$ $(y\mapsto
u(y)N(y,f,D))$ is integrable in $\mathbb G$ if and only if the
function $(u\circ f)(x)J(x,f)$ is integrable on $D$ $(A)$;
moreover, the following change of variable formulas hold:
\begin{equation}\label{in1}
\int\limits_{D}(u\circ f)(x)J(x,f)\,dx  = \int\limits_{\mathbb G
}u(y)\mu(y,f,D)\,dy,
\end{equation}
\begin{equation}\label{in2}
\int\limits_{A}u(x)J(x,f)\,dx  =  \int\limits_{\mathbb G
}\sum\limits_{x\in f^{-1}(y)\cap A}u(x)\,dy,
\end{equation}
\begin{equation}\label{in3}
\int\limits_{A}(u\circ f)(x)J(x,f)\,dx  =  \int\limits_{\mathbb G
}u(y)N(y,f,A)\,dy.
\end{equation}
}
\end{itemize}
\end{theorem}

We use the notation $\overline{\mathbb G}=\mathbb G\cup\{\infty\}$
for the one-point compactification of the Carnot group $\mathbb
G$.  The system of neighborhoods for $\{\infty\}$ are generated by
the complement to homogeneous closed balls. It is evident that
$\overline{\mathbb G}$ is topologically equivalent to the unit
Euclidean sphere $S^N$ in the Euclidean space $\mathbb R^{N+1}$.
Later on, we use the symbol $\Omega$ to denote a domain (open
connected set) on the Carnot group $\mathbb G$. It is not excluded
that $\Omega$ coincides with $\mathbb G$.

\begin{definition}\label{4}
A continuous mapping $f:\Omega\to\overline{\mathbb G}$ is said to
be a {\it quasimeromorphic mapping} if
\begin{itemize}
\item[1)] { $f:\Omega\setminus f^{-1}(\infty)\to \mathbb G$ is a
quasiregular mapping;} \item[2)] {for any domain
$\omega\Subset\Omega$, the multiplicity function $N(y,f,\omega)$
is essentially bounded:
$$N(f,\omega)=\esssup\limits_{y\in{\mathbb G}}N(y,f,\omega)=
\esssup\limits_{y\in{\mathbb
G}}\card\{f^{-1}(y)\cap\omega\}<\infty.$$}
\end{itemize}
\end{definition}

An ordered triplet $(F_0,F_1;\Omega)$ of nonempty sets, where
$\Omega$ is open in $\overline{\mathbb G}$, $F_0$ and $F_1$ are
compact subsets of $\overline{\Omega}$, is said to be a {\it
condenser} on  $\overline{\mathbb G}$. We define the $p$-{\it
capacity}, $1\leq p<\infty$, of the condenser $E=(F_0,F_1;\Omega)$
as
\begin{equation}\label{defcap}
\capac_p(E)=\capac_p(F_0,F_1;\Omega)=
\inf\int\limits_{\Omega\setminus\{\infty\}}|\nabla_{0} v|_0^p\,dx,
\end{equation}
where the infimum is taken over all nonnegative functions $v\in
C(\Omega\cup F_0\cup F_1)\cap L^1_{p}(\Omega\setminus\{\infty\})$
such that $v=0$ in a neighborhood of $F_0\cap\overline{\Omega}$
and $v\geq 1$ in a neighborhood of $F_1\cap\overline{\Omega}$.
Functions taking part in the definition of the $p$-capacity of a
condenser are said to be {\it admissible} for this condenser. If
the set of admissible functions is empty then the $p$-capacity of
a condenser equals infinity, by definition.

If $\Omega\subset\overline{\mathbb G}$ is an open set and $C$ is a
compact set in $\Omega$ then, for brevity, we denote the condenser
$E=(C,\partial \Omega;\Omega)$ by $E=(C,\Omega)$, and we shall
write $\capac_p(E)= \capac_p(C,\Omega)$ instead of
$\capac_p(C,\partial \Omega;\Omega)$. The notion of the
$p$-capacity $\capac_p(C,\Omega)$ is extended to an arbitrary set
$E\subset\Omega$ by the usual way (see, for
instance~\cite{Choq,HKM} in the case $\mathbb G=\mathbb R^n$
and~\cite{ChV1,ChV2} in the geometry of vector fields satisfying
the H\"{o}rmander hypoellipticity condition).

We say that a compact $C\subset\overline{\mathbb G}$ has the
$p$-capacity zero and write $\capac_p C=0$, if $\capac_p(C,U)=0$
for some open set $U\subset \overline{\mathbb G}$ such that
$\capac_p(\overline{B},U)>0$ for some ball $B\Subset U$. One can
prove
\begin{itemize}
\item[1)]{if $\capac_p(\overline{B}_1,U)>0$ for some ball
$B_1\Subset U$ then $\capac_p(\overline{B}_2,U)>0$ for an
arbitrary other ball $B_2\Subset U$;} \item[2)]{if
$\capac_p(C,U)=0$ then $\capac_p(C,V)=0$ whenever $V$ is any
bounded open set containing $C$.}
\end{itemize}
An arbitrary Borel set $E$ has the $p$-capacity zero, if the same
holds for any compact subset of $E$, otherwise $\capac_p E>0$.

The chosen homogeneous norm gives the following exact value for
the $p$-capacity of spherical rings $(\overline B(x,r),B(x,R))$,
$0<r<R<\infty$, ~\cite{BHT}:

\begin{equation}\label{272}\capac_{p}(\overline
B(x,r),B(x,R))=\left\{\array{ll} \kappa(\mathbb
G,p)\Big(\frac{|p-Q|}{p-1}\Big)^{p-1}\Big|R^{\frac{p-Q}{p-1}}-r^{\frac{p-Q}{p-1}}\Big|^{1-p},\quad
& p\neq Q,
\\
\kappa(\mathbb G,Q)\Big(\ln\frac{R}{r}\Big)^{1-Q}, & p=Q,
\endarray\right.\end{equation} where $\kappa(\mathbb G,p)$ is a positive
constant whose an exact value was obtained in~\cite{BT} (we give
it in Section~\ref{relations}).

\section{Properties of quasimeromorphic mappings}

\begin{lemma}[{\cite{Vod7}}]\label{indop} Let $f:\Omega\to\overline{\mathbb
G}$ be a quasimeromorphic mapping. For any open set
$U\subset\Omega$ such that $N(f,U)<\infty$ the operator
$f^*:L_Q^1(f(U)\setminus\{\infty\})\to L_Q^1(U\setminus
f^{-1}(\infty))$ where $f^*(u)=u\circ f$, is bounded$:$
\begin{equation}\label{est}
\|f^*(u)\mid L_Q^1(U\setminus
f^{-1}(\infty))\|\leq(K_O(f)N(f,U))^{1/Q}\|u\mid
L_Q^1(f(U)\setminus\{\infty\})\|,
\end{equation}
and the chain rule works:
$\nabla_{0}f^*(u)(x)=D_Hf(x)^T\nabla_0u(f(x))$ almost everywhere
in~$U$.
\end{lemma}

\begin{proof}
The set $U\setminus f^{-1}(\infty)$ is an open set in $\Omega$.
Consider an arbitrary function
$$
u\in C^1(f(U)\setminus\{\infty\})\cap
L^1_Q(f(U)\setminus\{\infty\}).
$$
Then $v=u\circ f\in \ACL(U\setminus f^{-1}(\infty))$ (since the
function~$u$ is locally Lipschitz) and $\nabla_{0} v(x)=D_Hf (x)
^T\nabla_{0}u\big(f(x)\big)$ almost everywhere in $U\setminus
f^{-1}(\infty)$ (since the mapping $f$ is $\mathcal
P$-differentiated a.~e.). Using~(\ref{in3}) and the property
$N(f,U\setminus f^{-1}(\infty))\leq N(f,U)<\infty$, we obtain
\begin{eqnarray*}
\int\limits_{U\setminus f^{-1}(\infty)}\big\vert\nabla_{0}(u\circ
f)\big\vert_0^Q(x)\, dx & \leq & \int\limits_{U\setminus
f^{-1}(\infty)}
|\nabla_{0}u|_0^Q\big(f (x)\big)|D_Hf|^Q(x)\,dx\\
& \leq & K_O(f)\int\limits_{U\setminus
f^{-1}(\infty)}|\nabla_{0}u|_0^Q\big
(f(x)\big)J(x,f)\,dx\\
& \leq & K_O
(f)N(f,U)\int\limits_{f(U)\setminus\{\infty\}}|\nabla_{0}u|_0^Q(y)\,dy.
\end{eqnarray*}
Since the composition operator $f^*:C^1\big
(f(U)\setminus\{\infty\}\big)\cap
L^1_Q\big(f(U)\setminus\{\infty\}\big)\mapsto L^1_Q(U\setminus
f^{-1}(\infty))$ is bounded, this operator can be continuously
extended to~$L^1_Q\big(f(U)\setminus\{\infty\})$ making use of
arguments of~\cite{UV1},  and the extended operator will be also
the composition operator.
\end{proof}

\begin{lemma}\label{12}
If $f: \Omega\to\overline{\mathbb G}$, $\Omega\subseteq\G$, is a
quasimeromorphic mapping and $S=f^{-1}(\infty)$, then
$$\capac_{Q}(S)=0.$$
\end{lemma}

\begin{proof}
According to the definition of the quasimeromorphic mapping, we
have~$S\not\equiv\Omega$. Therefore, there exists a ball $B_0$
such that $\overline{B}_0\subset\Omega\setminus S$. Let $\omega$
be a domain satisfying $\omega\Subset\Omega$, $B_0\Subset\omega$,
and $\omega\cap S\neq\emptyset$. We shall prove that
$\capac_{Q}(\overline B_0,S\cap\overline{B};\omega)=0$ for an
arbitrary ball $B\Subset\omega$ such that $\overline B_0\cap
\overline{B}=\emptyset$ and $S\cap\overline B\neq\emptyset$.

Fix some domain $\omega\Subset\Omega$, $\omega\cap
S\neq\emptyset$, a ball $B\Subset\omega$ with $S\cap\overline
B\neq\emptyset$, and a point $y\in f(\omega\setminus\overline
B)\setminus\{\infty\}$. For any $R>1$ we consider a condenser
$E_R=(\mathbb C B(y,R), \mathbb C \overline B(y,r))$ where $r<1$
is small enough to provide $B(y,r)\subset
f(\omega\setminus\overline B)\setminus\{\infty\}$ and
$f^{-1}(B(y,r))\subset\Omega$ contains some open ball
$B_0\Subset\omega$ satisfying $\overline B_0\cap
S\cap\overline{B}=\emptyset$. Notice that since $\overline
B(y,r)\subset B(y,R)$ we may choose as an admissible function for
the $Q$-capacity of the condenser $E_R$ a function $\varphi_{R}$
such that $\varphi_{R}\vert_{\overline B(y,r)}=0$ and
$\varphi_{R}\vert_{\partial B(y,R)}=1$. By this we define
$$
\varphi_R(z)=\left\{\array{llll} 0 & \qquad & \text{if}\qquad & z
\in \overline B(y,r),
\\
\frac{\ln\frac{|y^{-1}z|}r}{\ln \frac Rr} & \qquad &\text{if} &
z\in B(y,R)\setminus\overline B(y,r),
\\
1 & \qquad & \text{if}\qquad & z\notin B(y,R).
\endarray\right.
$$
Then $\capac_Q(E_R)\leq \int\limits_{B(y,R)\setminus\overline
B(y,r)}|\nabla_{0}\varphi_R(z)|_0^Q\,dz\leq C\bigl(\ln \frac
Rr\bigr)^{1-Q}$, where $C$ is the Lipschitz constant of the
function $z\mapsto |y^{-1}z|$. We denote by $F_R$ the set $\{x\in
\Omega: \varphi_R(f(x))=1\}$. Then $S\subset F_R$ for any real
$R>1$ and therefore $S\subset\bigcap\limits_{R\geq 2}F_R$. It is
clear also that $\omega\supset F_R\cap\overline{B}$ and
$f^*(\varphi_R)=\varphi_R\circ f$ is an admissible function for
the $Q$-capacity of the condenser $(\overline
B_0,S\cap\overline{B};\omega)$ for all $R\geq k_0$ where $k_0$ is
some number greater than one. Now, we use Lemma~\ref{indop} to
derive
\begin{eqnarray*}
\capac_Q(\overline B_0,S\cap\overline{B};\omega) & \leq &
\bigl\|f^*(\varphi_R)\mid
L_Q^1(\omega)\bigr\|^Q=\bigl\|f^*(\varphi_R)\mid
L_Q^1(\omega\setminus S)\bigr\|^Q
\\
& \leq & K_O(f)N(f,\omega) \bigl\|\varphi_R\mid
L_Q^1(f(\omega)\setminus\{\infty\})\bigr\|^Q
\\ & \leq & K_O(f)N(f,\omega)
C\Bigl(\ln\frac Rr\Bigr)^{1-Q}.
\end{eqnarray*}
The right-hand side of this inequality goes to 0 as $R\to \infty$.
Therefore,
$$\capac_Q(\overline B_0,S\cap\overline{B};\omega)=0$$ and the lemma is
proved.
\end{proof}

As a consequence of Lemma \ref{12} we have the following
property~\cite{Vod2}: if $f: \Omega\to\overline{\mathbb G}$,
$\Omega\subseteq\G$, is a quasimeromorphic mapping, then
$S(x,t)\cap f^{-1}(\infty)=\emptyset$ for an arbitrary point
$x\in\Om$ and for almost all $t$ such that the sphere $S(x,t)$
belongs to $\Om$.

We say that a mapping $f$ is {\it light} if $f^{-1}(y)$ is totally
disconnected for all $y$. Thus, from the previous considerations
we have the following statement.

\begin{corollary}
A quasimeromorphic mapping is light.
\end{corollary}

\subsection{Topological degree}\label{degree}

Recall that we identify the Carnot group $G$ with its Lie algebra
$\mathcal G$ and thus with $\mathbb R^N$,
$N=\sum\limits_{i=1}^{m}\dim V_i$. Moreover,  the one-point
compactification of $\overline{\mathbb G}$ is topologically
equivalent to the unit  sphere $S^{N}$ centered at $0$ in $\mathbb
R^{N+1}$. Therefore the topological degree $\mu(y,f,D)$ of a
continuous mapping $f:\Omega\to\overline{\mathbb G}$ where
$D\Subset\Omega$ is a compact domain, can be treated as the
topological degree of the continuous mapping $f:\Omega\to S^{N}$
with the standard orientation in $\Omega\subset \mathbb R^{N}$ and
$S^{N}$. The topological degree $\mu(y,f,D)$ of the continuous
mapping $f:\Omega\to \overline{\mathbb G}$ at $y$ is well-defined
whenever $D$ is a compact domain in $\Omega$ and
$y\in\overline{\mathbb G}\setminus f(\partial D)$. The degree is
integer-valued function and has the following properties:
\begin{itemize}
\item[1)]{the function $y\mapsto \mu(y,f,D)$ is a constant in
every connected component of $\overline{\mathbb G}\setminus
f(\partial D)$ and $\mu(y,f,D)=0$ if $y\notin f(\overline D)$;}
\item[2)]{if $U$ is a connected component of $\overline{\mathbb
G}\setminus f(\partial D)$ such that $\mu(y,f,D)\ne0$ for some
point  $y\in U$ then for any $z\in U$ there exists $x$ such that
$f(x)=z$;}
\item[3)]{if $y\in f(D)\setminus
f(\partial D)$ and the restriction of $f$ to
$\overline D$ is one-to-one then $|\mu(y,f,D)|=1$.}
\item[4)]{if $D_1,\ldots,D_k\Subset \Om$ are disjoint open sets and if
$D\cap f^{-1}(y)\subset\bigcup\limits_{i=1}^k D_i\subset D\Subset\Om$, then
$$
\mu(y,f,D)=\sum\limits_{i=1}^k\mu(y,f,D_i),\quad y\notin
f(\partial D),\ \text{and}\ y\notin f(\partial D_i),\
i=1,\ldots,k.
$$}
\end{itemize}

Other properties  of the mapping degree can be found
in~\cite{BI,RR,Resh1}.

\begin{lemma}\label{11}
Let $f:\Omega\to\overline{\mathbb G}$, $\Omega\subseteq\mathbb G$,
be a quasimeromorphic mapping. If $f(x_0)=\infty$ then the image
of any neighborhood of $x_0$ is a neighborhood of $\{\infty\}$.
\end{lemma}

\begin{proof}
Let $x_0$ be a point such that $f(x_0)=\infty$. Since $f$ is light
we can find a sphere $S(x_0,r)\in\Omega$ such that
$\{\infty\}\notin f(S(x_0,r))$. We choose an open connected
component $U_{\infty}\in\overline{\mathbb G}\setminus f(S(x_0,r))$
containing  $\{\infty\}$. There exists a point $z\in U_{\infty}$
such that $z=f(x)$ for some point $x\in B(x_0,r)$. According to
the properties of quasiregular mappings, the image
$W=f(B(x_0,r)\setminus f^{-1}(\infty))$ is an open neighborhood of
$z$. Then the following  properties hold
(see~\cite{Vod6,Vod7,Vod9}):
\begin{itemize}
\item[a)]{for  all  $y\in W$, the pre-image $f^{-1}(y)\cap
B(x_0,r)$ contains  finitely many points;} \item[b)]{for almost
all points $y\in W$, the $\mathcal P$-differential exists in all
points $x\in f^{-1}(y)$ and $J(x,f)$ does not vanish;}
\item[c)]{for almost all points $y\in W$,
$$\mu(y,f,B(x_0,r))
=\sum\limits_{x\in f^{-1}(y)}\sign J(x,f)> 0.$$}
\end{itemize}

By properties of the topological degree, the last expression
implies that the degree $\mu(y,f,B(x_0,r))$ does not vanish at all
points $y\in U_{\infty}$. Thus, for any $y\in U_{\infty}$ there
exists $x\in B(x_0,r)$ such that $f(x)=y$.
\end{proof}

\begin{corollary} A quasimeromorphic mapping is open and discrete.
\end{corollary}

\begin{proof} The openness follows from Lemma~\ref{11}
and the definition of quasimeromorphic mappings. If a map is open
and light, then it is discrete. The complete proof can be found
in~\cite{Resh1,Rick1,TY}.
\end{proof}

\begin{lemma}\label{estcap}
Let $f:\Omega\to\overline{\mathbb G}$ be a quasimeromorphic
mapping and $U\subset\Omega$ be a domain such that
$N(f,U)<\infty$. Then the condenser $E=(F_0,F_1;U)$ meets the
inequality
$$
\capac_{Q}(F_0,F_1;U)\leq
K_O(f)N(f,U)\capac_{Q}(f(F_0),f(F_1);f(U)).
$$
\end{lemma}

\begin{proof} We have to consider two cases: $F_1\cap f^{-1}(\infty)=\emptyset$
and $F_1\cap f^{-1}(\infty)\ne\emptyset$. The first case is well
known (see, for instance,~\cite{Vod7}). The second one is more
interesting for us. We note that since a quasimeromorphic mapping
is open, the triplet $(f(F_0),f(F_1);f(U))$ is a condenser. Let
$u$ be an admissible function for $(f(F_0),f(F_1);f(U))$. Then, in
view of Lemma~\ref{indop}, the function $u\circ f$ is admissible
for the condenser $(F_0,F_1;U)$ and $u\circ f=1$ in some
neighborhood of $f^{-1}(\infty)$. Therefore, at the same
neighborhood we have $\nabla_{0} (u\circ f)=0$. Applying estimate
(\ref{est}), we obtain
\begin{eqnarray}\label{vopros}
\capac_Q(F_0,F_1;U) & \leq & \int\limits_{U}|\nabla_{0}(u\circ
f)|_0^Q(x)\,dx= \int\limits_{U\setminus
f^{-1}(\infty)}|\nabla_{0}(u\circ f)|_0^Q(x)\,dx\nonumber
\\
& \leq &
K_O(f)N(f,U)\int\limits_{f(U)\setminus\{\infty\}}|\nabla_{0}u|_0^Q(z)\,dz.
\end{eqnarray}
Since $u$ is an arbitrary admissible function, the lemma is
proved.
\end{proof}

We need the following $Q$-capacity estimate.

\begin{theorem}[\cite{MV3,Vod7}]\label{thce}
Let $f:\Omega\to \mathbb G$ be a non-constant quasiregular mapping
and $E=(C,U)$ be a condenser such that $C$ is a compact in $U$ and
$U\Subset\Omega$. Then $f(E)=(f(C),f(U))$ is also a condenser and
\begin{equation}\label{capest}
\capac_Q(f(C),f(U))\leq K_I(f)\capac_Q(C,U).
\end{equation}
\end{theorem}

\begin{proof}
In the case $\mathbb G =\mathbb  R^n$ the estimate~\eqref{capest}
is proved in \cite{MRV2}. The proof in our case is based on the
following construction. Since $\overline U$ is compact then
$N(f,U)<\infty$. We define the \emph{pushforward}  of a non-negative function
$u\in C_0(U)$  to be the function
 $v = f_{\sharp}u :f(\Omega)\rightarrow \mathbb R$, given  by
$$
 v(y):= \begin{cases}
 \sup\{ u(x)\, : \,f(x)=y \} & \text{if } y\in f(U),
\\ 0 & \text{otherwise.}
 \end{cases}
$$

By the same way as in~\cite[Lemma 7.6]{MRV2}, one can prove that
if $f$ is continuous discrete and open, and the non-negative
function $u : U \to \mathbb R$ is continuous with compact support,
then the function $v = f_{\sharp}u  :f(\Omega)\rightarrow \mathbb
R$ is also continuous and $\supp {v} \subset f(\supp u)$.
Moreover, if additionally $u\in C_{0}^{1}(U)$ and the mapping $f$
is quasiregular then $v = f_{\sharp}u$ belongs to
$W^1_Q(f(\Omega))$. Below the precise statement
follows~\cite{MV3,Vod7}.

{\it Let $f:\Omega\to \mathbb G$ be a non-constant quasiregular
mapping. Then the operator $f_{\sharp}$ possesses the following
properties:
\begin{enumerate}
\renewcommand{\theenumi}{\arabic{enumi}}
\def\labelenumi{\arabic{enumi})}
\item  $f_{\sharp}: C^1_0(U)^+\to W^1_Q(f(\Omega))\cap
C_0(f(\Omega))$ where the symbol $C^1_0(U)^+$ denotes all
non-negative functions of $C^1_0(U)$, \item
$\int\limits_{f(\Omega)}|\nabla_{0} f_{\sharp}(u)|_0^Q\, dx \leq
K_I\int\limits_{U}|\nabla_{0} u|_0^Q\,dx$ for any $u\in C^1_0(U)$,
\item  if the function $u$ is admissible for the condenser
$E=(U,C)$ then $f_{\sharp}u$  is admissible for the condenser
$f(E)=(f(U),f(C))$.
\end{enumerate}}

To prove the proposition one needs to check that $f_{\sharp}u\in
\ACL (f(\Omega))$ (see details in \cite{Vod7} where
$\ACL$-property is verified  for a function of similar nature).

>From the last two properties everyone can deduce the
inequality~\eqref{capest}.
\end{proof}

We use the estimate~(\ref{capest}) to prove the removability
property of quasimeromorphic mappings. Before to formulate it we
prove some auxiliary assertions.

Let $E\subset\overline{\mathbb G}$ be a closed set of positive $Q$-capacity.
We say that the set $E$ has the essentially positive $Q$-capacity at a point
$x\in E$, $x\ne\infty$, if
\begin{equation}\label{ce}
\capac_Q(E\cap\overline{B(x,r)},B(x,2r))>0
\end{equation}
for any positive $r$. One is able to check that
\begin{itemize}
\item[1)]{it is sufficient to verify~(\ref{ce}) for $r\in(0,r_0)$,
where $r_0$ is a positive number;} \item[2)] {the set
$$
\widetilde E=\{x\in E:\ \text{the set $E$ has the essentially positive $Q$-capacity
at $x$}\}
$$
is not empty and closed.}
\end{itemize}
Then there exists a point $x_0$ such that
$|x_0|=\inf\{|x|:x\in\widetilde E\}$. Let us denote the
intersection $E\cap B(x_0,1)$ by the symbol $E_0$. By definition, we have
$\capac_Q(E_0,B(x_0,2))$ is positive.

\begin{lemma}[\cite{MV3}]\label{17}
Let $E$ be a closed subset of $\overline{\mathbb G}$ with
$\capac_Q(E)>0$. Then for every $a>0$ and $d>0$ there exists
$\delta>0$ such that $\capac_Q(C,\mathbb C E)\geq \delta$ whenever
$C\subset\mathbb C E$ is a continuum such that $\diam(C)\geq a>0$
and $\dist(C, E_0)\leq d$.
\end{lemma}

\begin{proof} It is enough to prove the assertion under assumption
that $E$ is a non-empty bounded set.

We use the rule of contraries. Then there exist $a>0$ and $d>0$
such that for any $\delta_n=\frac{1}{n}$, $n\in\mathbb N$, we can
find a continuum $C_n$ with the diameter $\diam(C_n)\geq a>0$ and
$\dist(C_n, E_0)\leq d$ but $\capac_Q(C_n,\mathbb C E)\leq
\delta_n$. By these assumptions we derive existence of a real
number $R$, $R\geq d>0$, such that some connected part of the
intersection $\gamma_n=C_n\cap\overline{B(x_0,R)}$ has the
diameter $\diam(\gamma_n)\geq a/2>0$ and
\begin{equation}\label{ce1}
\capac_Q(E_0\cap\overline{B(x_0,R)},B(x_0,2R))>0.
\end{equation}

Since \begin{eqnarray*}\capac_Q(C_n,\mathbb C E) & \geq &
\capac_Q(\gamma_n,\mathbb C E)\geq\capac_Q(\gamma_n,\mathbb
C(E_0\cap\overline{B(x_0,R)})) \\ & \geq &
\capac_Q(\gamma_n,E_0\cap\overline{B(x_0,R)};B(x_0,2R))\end{eqnarray*}
we can choose admissible functions $\varphi_n(x)\in
C(\overline{B(x_0,2R)})\cap L^1_Q(B(x_0,2R)))$ for condensers
$(\gamma_n,E_0\cap\overline{B(x_0,R)};B(x_0,2R))$ such that
$\varphi_n(x)\in(0,1)$ when $x\in B(x_0,2R)$,
$$
\varphi_n(x)=\left\{\array{lll}
0\quad & \text{if}\quad & x\in \gamma_n,
\\
1\quad & \text{if}\quad & x\in E_0\cap\overline{B(x_0,R)},
\endarray\right.
$$
and
$$
\int\limits_{B(x_0,2R)}|\nabla_{0}\varphi_n|_0^Q\,dx\to 0\quad
\text{as $n\to\infty$}.
$$
Using Poincar'e inequality we can extract a subsequence (that we
denote by the same symbol) such that
$\varphi_n(x)\to\alpha\in[0,1]$ almost everywhere in $B(x_0,2R)$
as $n\to\infty$. Additionally, we can also assume that
$|\nabla_{0}\varphi_n(x)|_0\to 0 $ almost everywhere in
$B(x_0,2R)$ as $n\to\infty$.

Let $\psi\in C_0(B(x_0,2R))\cap L^1_Q(B(x_0,2R))$ be a function
such that $ \psi(x)=1$ if $x\in
\big(E_0\cap\overline{B(x_0,R)}\big)\bigcup\gamma_n$ and
$\psi(x)\in[0,1]$.

{\sc 1st case}: $\alpha<1$. The product
$g_n=(1-\alpha)^{-1}(\varphi_n-\alpha)\psi$ is an admissible
function for the condenser
$(E_0\cap\overline{B(x_0,R)},B(x_0,2R))$. Since
$|\nabla_{0}\varphi_n|_0\to 0$ and $|\varphi_n-\alpha|\to 0$
almost everywhere as $n\to 0$ we derive
\begin{eqnarray*}
\int\limits_{B(0,2R)}|\nabla_{0}g_n|_0^{Q}\,dx & \leq &
\frac{2^{Q-1}}{(1-\alpha)^Q}\Big(
\int\limits_{B(0,2R)}|\psi\nabla_{0}\varphi_n|_0^{Q}\,dx \\ & + &
\int\limits_{B(0,2R)}|(\varphi_n-\alpha)\nabla_{0}\psi|_0^{Q}\,dx\Big)\to
0
\end{eqnarray*}
as $n\to\infty$ by the Lebesgue dominated theorem. This
contradicts to (\ref{ce1}).

{\sc 2nd case}: $\alpha=1$. In this case the product
$g_n=\psi(1-\varphi_n)$ is an admissible function for a condenser
$(\gamma_n,B(x_0,2R))$, $n\in\mathbb N$. According to the above
estimates $\int\limits_{B(0,2R)}|\nabla_{0}g_n|_0^{Q}\,dx\to0$ as
$n\to\infty$. Results of \cite{H} (see also \cite{ChV2,KV}) imply
that $\diam(\gamma_n)\to0$ as $n\to\infty$ that contradicts to the
choice of $\gamma_n$.
\end{proof}

\begin{theorem}[\cite{MV3}]\label{114}
Let $\Omega$ be a domain in $\mathbb G$, $E\subset\Omega$ be a
closed set with $\capac_Q(E)=0$. If $f: \Omega\setminus
E\to\overline{\mathbb G}$ is a quasimeromorphic mapping and
$\capac_Q(\mathbb C f(\Omega\setminus E))$ is positive, then $f$
can be extended to a continuous mapping $f^*:
\Omega\to\overline{\mathbb G}$. Moreover, if the domain $\Omega$
is unbounded then there exists also a limit
\begin{equation}\label{infin}
\lim\limits_{x\to\infty,\,x\in\Omega}f^*(x)\in\overline{\mathbb
G}.
\end{equation}
\end{theorem}

\begin{proof} We may assume that
$E$ contains the set $f^{-1}(\infty)$. Since the $Q$-capacity of
$E$ is zero, the set $\Omega\setminus E$ is connected. To show
that $f$ has a limit at a point $b\in E$ we choose a sphere
$S(b,R)\subset\Omega\setminus E$ and two different sequences
$\{x_j\}\in\Omega\cap B(b,R/4)$, $\{x^{\prime}_j\}\in\Omega\cap
B(b,R/4)$ going to $b$ as $j\to\infty$. Let
$r_j=2\max\{|b^{-1}x_j|,|b^{-1}x^{\prime}_j|\}$. By $C_j\subset
B(b,r_j)\setminus E$, we denote a rectifiable curve with endpoints
$x_j$ and $x^{\prime}_j$. In view of $\capac_Q(E)=0$ the set $E$
is removable~\cite{VodMar} and we have
\begin{eqnarray}\label{zero}
\capac_Q(C_j,B(b,R)\setminus E)) & = &
\capac_Q(C_j,B(b,R))\nonumber
\\
&  \leq & \capac_Q(\overline B(b,r_j),B(b,R))=\kappa(\mathbb
G,Q)\Big(\ln \frac{R}{r_j}\Big)^{1-Q}.
\end{eqnarray}
Since the right-hand side of this relation tends to $0$ as $j\to
\infty$ then the left-hand side of it does the same.

Suppose that $f: \Omega\setminus E\to\overline{\mathbb G}$ has no
limit in $b\in E$. It follows that for some subsequences (that we
denote by the same symbols $\{x_j\}$ and $\{x_j^{\prime}\}$) we
have simultaneously
\begin{itemize}
\item[a)] {at least one of the sequences $f(x_j)$ and
$f(x^{\prime}_j)$ is bounded in $f(\Omega\setminus E)$,}
\item[b)]{$\diam (f(C_j))\geq\alpha>0$ for some constant
$\alpha>0$ and for all $j\in\mathbb N$.}
\end{itemize}
Applying Lemma \ref{17}, we obtain the inequality
\begin{equation}\label{positive}
\capac_Q(f(C_j),f(\Omega\setminus E))\geq \delta>0
\end{equation}
for some $\delta$ and all $j\in\mathbb N$.

On the other hand, for fixed $j$, we can exhaust a domain
$B(b,R)\setminus  E$ by compact domains $\omega_k$ such that
$C_j\subset\omega_1\Subset\ldots\Subset\omega_k\Subset\ldots
\Subset B(b,R)\setminus E$, $\bigcup\limits_{k}\omega_k=
B(b,R)\setminus E$. By Theorem \ref{thce} and properties of
capacity, we deduce
\begin{equation}\label{capest0}
\capac_Q(f(C_j),f(\Omega\setminus E))\leq
\capac_Q(f(C_j),f(\omega_k))\leq K_I(f)\capac_Q(C_j,\omega_k).
\end{equation}
Letting $k\to\infty$ in the right-hand side of~(\ref{capest0}), we
obtain
\begin{equation}\label{capest1}
\capac_Q(f(C_j),f(\Omega\setminus E))\leq K_I(f)\capac_Q(C_j,
B(b,R)\setminus E)
\end{equation} by
properties of the capacity. The inequalities (\ref{positive}) and
(\ref{capest1}) imply
\begin{equation}\label{capest2}
0<\delta\leq\capac_Q(f(C_j),f(\Omega\setminus E)) \leq
K_I(f)\capac_Q(C_j,B(b,R)\setminus E).
\end{equation}
We have a contradiction, since, by~(\ref{zero}), the right-hand side of
(\ref{capest2}) goes to 0 as $j\to 0$.

It remains to show~(\ref{infin}). Since the set $\Omega\setminus
E$ is open we can find a ball $B(x_0,R)\Subset\Omega\setminus E$
and two different sequences $\{x_j\}\in(\Omega\setminus E)\cap
\mathbb C\overline B(x_0,4R)$,
$\{x^{\prime}_j\}\in(\Omega\setminus E)\cap \mathbb C\overline
B(x_0,4R)$ going to $\infty$ as $j\to\infty$. Put
$r_j=\frac{1}{2}\min\{|x_0^{-1}x_j|,|x_0^{-1}x_j^{\prime}|\}$. We
denote by $C_j$ a rectifiable curve connecting points $x_j$ and
$x_j^{\prime}$ in $(\Omega\setminus E)\cap \mathbb C\overline
B(x_0,r_j)$ if the points $x_j$ and $x_j^{\prime}$ belong to the
same connected component of $(\Omega\setminus E)\cap \mathbb
C\overline B(x_0,r_j)$. In the case when $x_j$ and $x_j^{\prime}$
are in the different components of $(\Omega\setminus E)\cap
\mathbb C\overline B(x_0,r_j)$, then $C_j$ will denote the union
of rectifiable curves joining $x_j$ and $x_j^{\prime}$ with the
sphere $S(x_0,r_j)$ in $(\Omega\setminus E)\cap \mathbb C\overline
B(x_0,r_j)$. Then
\begin{eqnarray}\label{268}\capac_Q(C_j,(\Omega\setminus E)\cap \mathbb
C\overline B(x_0,R)) & = & \capac_Q(C_j,\Omega\cap \mathbb
C\overline B(x_0,R))\nonumber
\\ & \leq &
\capac_Q(S(x_0,r_j)\cap\Omega,S(x_0,R);\Omega)\nonumber
\\ & \leq &
\capac_Q(\overline B(x_0,R), B(x_0,r_j)) \\ & = & \kappa(\mathbb
G,Q)\Big(\ln\frac{r_j}{R}\Big)^{1-Q}.\nonumber\end{eqnarray}

If~(\ref{infin}) does not exist then the following limit
$$\lim\limits_{x\to\infty,\,x\in\Omega\setminus E}f(x)$$ does the
same and the properties~a), ~b) hold. Then, by Lemma \ref{17}, we
obtain the inequality (\ref{positive}) for some positive $\delta$
and all $j\in\mathbb N$. Arguing like above, we get an analogue of
(\ref{capest2}):
\begin{eqnarray*}
0 & < & \delta\leq\capac_Q(f(C_j),f(\Omega\setminus E))\leq
\capac_Q\big(f(C_j),f((\Omega\setminus E)\cap\mathbb C\overline B(x_0,R))\big)\\
& \leq & K_I(f) \capac_Q(C_j,(\Omega\setminus E)\cap\mathbb C
\overline B(x_0,R)).
\end{eqnarray*}
We come to a contradiction, since the right-hand side goes to $0$
as $j\to\infty$ by~(\ref{268}).
\end{proof}

\begin{definition}\label{6} If $f: \Omega\to \overline{\mathbb G}$ is a quasimeromorphic
mapping, and if $b$ is an isolated point of $\partial\Omega$ such
that $f$ has no limits at $b$, then we call $b$ the {\it
$($isolated$)$ essential singularity} of $f$. \end{definition}

\begin{corollary}\label{32}
Let $b$ be an isolated essential singularity of a quasimeromorphic
mapping $f:\Omega \to\overline{\mathbb G}$. Then $\capac_Q(\mathbb
C f(U\setminus b))=0$ for an arbitrary neighborhood
$U\subset\Omega\cup\{b\}$ of the point $b$.
\end{corollary}

\begin{proof}
If we suppose that there exists a neighborhood
$U\subset\Omega\cup\{b\}$ such that $\capac_Q(\mathbb C
f(U\setminus b))>0$ then $f\vert_{U\setminus\{b\}}$ is extended to
the point $b$ by Theorem~\ref{114}. This contradicts to the
assumption that $b$ is the essential singularity of $f$.
\end{proof}

\begin{lemma}[\cite{MV3}]\label{19}
Let $b$ be an isolated essential singularity of a quasimeromorphic
mapping $f:\Omega \to\overline{\mathbb G}$. Then there exists an
$\mathcal F_{\sigma}$-set $E\subset\overline{\mathbb G}$ with
$\capac_Q(E)=0$, such that $N(y,f,U\setminus \{b\})=\infty$ for
every $y\in\overline{\mathbb G}\setminus E$ and all neighborhoods
$U\subset\Omega\cup\{b\}$ of the point $b$.
\end{lemma}

\begin{proof}
We consider two possibilities. If $b\neq\infty$ then we can assume
that $B(b,1)\subset\Omega\cup\{b\}$ and denote the set
$B(b,\frac{1}{k})$  by $V_k$, $k=1,2,\ldots$. In the case
$b=\infty$, we will use the notation $V_k=\mathbb C \overline
B(0,k)\cap\Omega$, $k=1,2,\ldots$. Set
$E=\bigcup\limits_{k=1}^{\infty}\mathbb C f(V_k)$. Then $E$
contains $f(U\setminus\{b\})$ for any neighborhood
$U\subset\Omega\cup\{b\}$ of the point $b$. The properties of the
capacity imply that $\capac_Q(E)=0$ because of $\capac_Q(\mathbb C
f(U\setminus \{b\}))=0$ by Corollary~\ref{32}. For an arbitrary
neighborhood $U\subset\Omega\cup\{b\}$ and given
$y\in\overline{\mathbb G}\setminus E=\bigcap\limits_{k=1}^{\infty}
f(V_k)$, we can find a sequence $\{x_j\}$ with pairwise disjoint
elements such that $x_j\in V_{k_j}$, $f(x_j)=y$. The lemma is
proved. \end{proof}

\section{Main inequalities for modulus}
Here and subsequently $\langle a,b\rangle$ stands for an interval
of one of the following type $(a,b)$, $[a,b)$, $(a,b]$, and
$[a,b]$. We say that a curve $\gamma:\langle
a,b\rangle\to\overline\G$ is locally rectifiable, if $\gamma$ is
locally rectifiable on $\langle a,b\rangle\setminus
\gamma^{-1}(\infty)$. A restriction
$\gamma^{\prime}=\gamma\vert_{[\alpha,\beta]}$,
$[\alpha,\beta]\subset\langle a,b\rangle\setminus
\gamma^{-1}(\infty)$, is said to be a closed part of $\gamma$. The
closed part $\gamma^{\prime}$ is rectifiable and we may use the
length arc parameter $s$ on $\gamma^{\prime}$. The linear integral
is defined by $$
\int\limits_{\gamma}\rho\,ds=\sup\int\limits_{\gamma^{\prime}}\rho\,ds=
\sup\int\limits_{0}^{l(\gamma^{\prime})}\rho(\gamma^{\prime}(s))\,ds,$$
where the supremum is taken over all closed parts
$\gamma^{\prime}$ of $\gamma$ and $l(\gamma^{\prime})$ is the
length of~$\gamma^{\prime}$. Let $\Gamma$ be a family of curves in
$\overline\G$. Denote by $\mathcal F(\Gamma)$ the set of Borel
functions $\rho :\overline\G\to[0;\infty]$ such that the
inequality
$$
\int\limits_{\gamma}\rho\,ds\geq 1$$ holds for a locally
rectifiable curve $\gamma\in \Gamma$. Otherwise we put
$\int\limits_{\gamma}\rho\,ds=\infty$. An element of the family
$\mathcal F(\Gamma)$ is called an {\it admissible density}
for~$\Gamma$.

\begin{definition}\label{5} Let $\Gamma$ be a family of curves in
$\overline\G$ and $p\in (1,\infty)$. The quantity
$$M_p(\Gamma)=\inf\int\limits_{\mathbb G}\rho^p\,dx$$ is called
the $p$-{\it module of the family of curves} $\Gamma$. The infimum
is taken over all admissible densities $\rho \in\mathcal
F(\Gamma)$.
\end{definition}

It is known that the $p$-module of a family of non-rectifiable
curves vanishes~\cite{Fug}. If some property fails to hold for a
family of curves whose $p$-module vanishes, then we say that the
property holds {\it $p$-almost everywhere}.

Let $F_0$, $F_1$ be disjoint compacts in $\overline\Omega$. We say
that a curve $\gamma:\langle a,b\rangle\to\Omega$ connects $F_0$
and $F_1$ in $\Omega$ (starts on $F_0$ in $\Omega$) if
\begin{itemize}\item[1.]{$\overline{\gamma(\langle a,b\rangle)}\cap F_i\neq\emptyset$, $i=0,1$,
($\overline{\gamma(\langle a,b\rangle)}\cap F_0\neq\emptyset$),}
\item[2.]{$\gamma(t)\in\Omega$ for all $t\in(a,b)$.}
\end{itemize}
A family of curves connecting $F_0$ and $F_1$ (starting at $F_0$)
in $\Omega$ is denoted by $\Gamma(F_0,F_1;\Omega)$
($\Gamma(F_0;\Omega)$). In the next theorem the relation between
the $p$-capacity of the condenser $(F_0,F_1;\Omega)$ and the
$p$-module of the family $\Gamma(F_0,F_1;\Omega)$ is given.

\begin{theorem}\cite{Mar1}
Let $\Omega$ be a bounded domain in the Carnot group $G$. Suppose
that $K_0$ and $K_1$ are disjoint non-empty compact sets in the
closure of $\Omega$. Then
$$M_p(\Gamma(F_0,F_1;\Omega))=\capac_p(F_0,F_1;\Omega),\qquad p\in(1,\infty).$$
\end{theorem}

\begin{remark}\label{r3} Let $f:\Omega\to\overline{\mathbb G}$ be
a quasimeromorphic mapping and $\Gamma$ be a family of curves in
$\Omega$. We correlate the parametrization of the curves in
$\Gamma\subset\Omega$ and in $\Gamma^{\star}=f(\Gamma)\subset
f(\Omega)$. Let $\gamma^{\star}\in\Gamma^{\star}$ be a rectifiable
curve. We introduce the length arc parameter $s^{\star}$ in the
curve $\gamma^{\star}\in\Gamma^{\star}$. Thus $s^{\star}\in
I^{\star}=[0,l(\gamma^{\star})]$ where $l(\gamma^{\star})$ is the
length of the curve $\gamma^{\star}$. If $t$ is any other
parameter on $\gamma^{\star}$: $\gamma^{\star}(t)=f(\gamma(t))$,
then the function $s^{\star}(t)$ is strictly monotone and
continuous, so the same holds for its inverse function
$t(s^{\star})$. For the curve $\gamma(t)\in\Gamma$ such that
$f(\gamma(t))=\gamma^{\star}$ the parameter $s^{\star}$ can be
introduced by the following way
$$
f(\gamma(t(s^{\star})))=f(\gamma(s^{\star}))=\gamma^{\star}(s^{\star}),\quad
s^{\star}\in I^{\star}.$$ We note that if we take the length arc
parameter $s$ on $\gamma$, $s\in I=[0,l(\gamma)]$, then
\begin{equation}\label{234}
1=\Big|\frac{d\gamma(s)}{ds}\Big|_0=
\Big|\frac{d\gamma(s^{\star})}{ds^{\star}}\Big|_0\cdot\Big|\frac{ds^{\star}}{ds}\Big|.
\end{equation}
\end{remark}
>From now on, we use the letters $s$ and $s^{\star}$ to denote the
length arc parameters on curves $\gamma$ and
$\gamma^{\star}=f(\gamma)$. The corresponding domains of $s$ and
$s^{\star}$ are denoted by $I=[0,l(\gamma)]$ and
$I^{\star}=[0,l(\gamma^{\star})]$, respectively.

\begin{theorem}\label{K0ineq}
Let $f:\Omega \to\overline{\mathbb G}$ be a nonconstant
quasimeromorphic mapping. Then, for a Borel set $A$,
$A\subset\Omega$, with $N(f,A)<\infty$ and a family of curves
$\Gamma$ in $A$, we have
$$M_Q(\Gamma)\leq K_O(f)N(f,A)M_Q(f(\Gamma)).$$
\end{theorem}
\begin{proof}
Since the set $S=f^{-1}(\infty)$ has a $Q$-capacity zero we have
$\mes(S)=0$. The set of points, where $f$ is not $\mathcal
P$-differentiable, is included into a Borel set $F$ with
$\mes(F)=0$. We write $E$ for the union $F\cup B_f\cup S$.

Take an admissible function $\rho^{\star}$ for a family of curves
$f(\Gamma)$ and set
$$\rho(x)=\left\{\array{cll}
\rho^{\star}(f(x))|D_Hf(x)|\quad &\text{for}\quad & x\in
A\setminus E,
\\
\infty &\text{for}& x\in E,
\\
0 &\text{for}& x\notin A.
\endarray\right.$$
The function $\rho$ is Borel and nonnegative. Denote by
$\Gamma_0\subset\Gamma$ the subfamily of locally rectifiable
curves such that there exists an image $f(\gamma^{\prime})$ of a
closed part $\gamma^{\prime}$ of $\gamma$ that is not absolutely
continuous. A result of B.~Fuglede~\cite{Fug} implies that
$M_Q(\Gamma_0)=0$. Notice also, that since $\mes(E)=0$ the family
of curves $\Gamma_1\subset\Gamma$ where
$\int\limits_{\gamma}\chi_E\,ds>0$, $\gamma\in\Gamma_1$, has the
$Q$-module zero. In fact, define $\tilde\rho:\ A\to\mathbb R^1$ as
$$\tilde\rho(x)=\left\{\array{cll}
\infty &\text{for}& x\in E,
\\
0 &\text{for}& x\notin E.
\endarray\right.$$ Then $\int_{\gamma}\tilde\rho\,ds=\infty$ for
$\gamma\in\Gamma_1$. Hence, $\tilde\rho$ is admissible for
$\Gamma_1$ and $$M_Q(\Gamma_1)\leq\int_{E}\tilde\rho^Q\,dx=0.$$

Suppose that the closed parts $\gamma^{\prime}$ and
$f(\gamma^{\prime})$ of curves
$\gamma\in\Gamma\setminus(\Gamma_0\cup\Gamma_1)$ and $f(\gamma)\in
f(\Gamma\setminus(\Gamma_0\cup\Gamma_1))$ are parameterized as in
Remark~\ref{r3}. We have
\begin{equation}\label{237}
1=\Big|\frac{df(\gamma^{\prime}(s^{\star}))}{ds^{\star}}\Big|_0\leq
\Big|D_Hf(\gamma^{\prime})\Big|\cdot\Big|\frac{d\gamma^{\prime}(s^{\star})}{ds^{\star}}\Big|_0
=\Big|D_Hf(\gamma^{\prime})\Big|\cdot\Big|\frac{ds}{ds^{\star}}\Big|,
\end{equation} where $s$ is the length arc parameter on
$\gamma^{\prime}$.
By~(\ref{237}), we deduce
\begin{eqnarray}\label{235}
\int\limits_{\gamma}\rho\,ds & \geq
&\int\limits_{\gamma^{\prime}}\rho\,ds  \geq
\int\limits_{I}\rho^{\star}(f(\gamma^{\prime}(s)))|D_Hf(\gamma^{\prime}(s))|\,ds\nonumber
\\ & = &
\int\limits_{I^{\star}}\rho^{\star}(f(\gamma^{\prime}(s^{\star})))|D_Hf(\gamma^{\prime})|
\Big|\frac{ds}{ds^{\star}}\Big|\,ds^{\star}\geq
\int\limits_{f(\gamma^{\prime})}\rho^{\star}\,ds^{\star}
\end{eqnarray} for any closed part $\gamma^{\prime}$ of $\gamma\in\Gamma\setminus(\Gamma_0\cup\Gamma_1)$.
Taking supremum over all closed parts of $f(\gamma)$ we see
$\int_{\gamma}\rho\,ds
\geq\int_{f(\gamma)}\rho^{\star}\,ds^{\star}\geq 1$. We conclude
that $\rho$ is admissible for the family
$\Gamma\setminus(\Gamma_0\cup\Gamma_1)$.

Since $\mes(E)=0$ we write \begin{eqnarray*} M_Q(\Gamma) & =&
M_Q(\Gamma\setminus(\Gamma_0\cup\Gamma_1)) \leq
\int\limits_{\G}\rho^{Q}\,dx=\int\limits_{A}\rho^{\star}(f(x))^Q|D_Hf(x)|^Q\,dx
\\ & \leq &
K_O(f)\int\limits_{A}\rho^{\star}(f(x))^Q
J(x,f)\,dx=K_O(f)\int\limits_{\G}\rho^{\star}(y)^Q N(y,f,A)\,dy
\\ & \leq &
K_O(f)N(f,A)\int\limits_{\G}\rho^{\star}(y)^Q \,dy
\end{eqnarray*} by~(\ref{22}) and~(\ref{in3}). Taking infimum over
all $\rho^{\star}\in\mathcal F(f(\Gamma))$ we end the proof.
\end{proof}

\begin{remark}\label{r2} Within the proof of Theorem~\ref{K0ineq} we have derived the estimate
\begin{equation}\label{229} M_Q(\Gamma)\leq K_O(f)\int\limits_{\mathbb
G}\rho(y)^Q N(y,f,A)\,dy\end{equation} which will be used below.
Inequality~(\ref{229}) holds for any function $\rho\in\mathcal
F(f(\Gamma))$.
\end{remark}

We state here Poletski\u{\i} type lemma. Its complete proof can be
found in~\cite{Mar,Mar2}.

\begin{lemma}\label{Pollem}
Let $f:\Omega\to\mathbb G$ be a non-constant quasiregular mapping
and $U\subset\Omega$ be a domain, such that $\overline
U\subset\Omega$. Assume $\Gamma$ to be a family of curves in $U$
such that $\gamma^{\star}(s^{\star})=f(\gamma(s^{\star}))$ is
locally rectifiable and there exists a closed part
$\gamma^{\prime}(s^{\star})$ of $\gamma(s^{\star})$ that is not
absolutely continuous (the parameterization of $\Gamma$ and
$f(\Gamma)$ is correlated as in Remark~\ref{r3}). Then,
$M_Q(f(\Gamma))=0$
\end{lemma}

\begin{theorem}\label{Polineq}
Let $f:\Omega \to\overline{\mathbb G}$ be a nonconstant
quasimeromorphic mapping and $\Gamma$ be a family of curves in
$\Omega$. Then
\begin{equation}\label{233} M_Q(f(\Gamma))\leq
K_I(f)M_Q(\Gamma).\end{equation}
\end{theorem}

\begin{proof}
Let $\rho$ be an admissible function for a family $\Gamma$. We can
assume that $\int\limits_{\Omega}\rho^Q\,dx<\infty$. We take a
sequence $\Omega_1\Subset\Omega_2\Subset\ldots\Subset\Omega$ of
subdomains that exhausts $\Omega$. Then
$\int\limits_{\Omega\setminus\Omega_i}\rho^Q\,dx\to 0$ as $i\to
\infty$. Now we define an admissible function for a family
$f(\Gamma)$. If $x\in \Omega_i\setminus B_f$ then $f$ is $\mathcal
P$-differentiable almost everywhere with strictly positive
Jacobian $J(x,f)$. By $F$, we denote a Borel set of $\mes(F)=0$
containing all points $x\in \Omega_i\setminus B_f$ where $f$ is
not $\mathcal P$-differentiable. Notice that $\mes (F\cup
B_f)=\mes f(F\cup B_f)=0$. Let $E\subset\overline{\mathbb G}$ be a
Borel set of zero measure such that $f(F\cup
B_f)\cup\{\infty\}\subset E$. We set
$\lambda_f(x)=\frac{1}{\min\limits_{|\xi|_0=1,\xi\in
V_1}|D_Hf(x)\xi|_0}$ for $x\in \Omega_i\setminus(F\cup B_f)$ and
define
$$\rho^{\star}_i(y)=\left\{\array{cll}
\max\limits_{x\in f^{-1}(y)\cap
\Omega_i}\big(\rho(x)\lambda_f(x)\big)\quad &\text{for}\quad &
y\in f(\Omega_i)\setminus E,
\\
\infty\quad &\text{for}\quad & y\in E,
\\
0\quad &\text{for}\quad & y\notin f(\Omega_i).
\endarray\right.$$ The function $\rho^{\star}_i(y)$ is nonnegative
and Borel.

Since the $Q$-module of a family of non-rectifiable curves
vanishes we consider only rectifiable curves. We correlate the
parametrization of the curves in $\Gamma$ and
$\Gamma^{\star}=f(\Gamma)$ as in Remark~\ref{r3}. Denote by
$I_i^{\star}$ the maximal sub-interval of $[0,l(\gamma^{\star})]$
such that $\gamma\vert_{I_i^{\star}}\subset\Omega_i$. Let
$\Gamma_i$ be the family of curves
$\gamma_i=\gamma\vert_{I_i^{\star}}$ and
$\Gamma_i^{\star}=f(\Gamma_i)$. Lemma~\ref{Pollem} states that the
$Q$-module of a family $\Gamma_{0,i}^{\star}$ of curves
$\gamma^{\star}_i\in\Gamma_i^{\star}$, for which $\gamma_{i}$ is
not absolutely continuous, vanishes. As in Theorem~\ref{K0ineq} it
can be shown that the family of curves
$\Gamma_1^{\star}\subset\Gamma^{\star}$ where
$\int\limits_{\gamma^{\star}}\chi_E\,ds>0$, $E=f(F\cup
B_f)\cup\{\infty\}$, $\gamma^{\star}\in\Gamma_1^{\star}$, has the
$Q$-module zero. So, we restrict our attention to the family
$\Gamma_i^{\star}\setminus(\Gamma_{0,i}^{\star}\cup\Gamma_1^{\star})$.

Notice that
\begin{equation}\label{230}
1=\Big|\frac{df(\gamma_i(s^{\star}))}{ds^{\star}}\Big|_0\geq
\Big|\frac{1}{\lambda_f(\gamma_i)}\Big|\cdot\Big|\frac{d\gamma_i(s^{\star})}{ds^{\star}}\Big|_0
=\Big|\frac{1}{\lambda_f(\gamma_i)}\Big|\cdot\Big|\frac{ds}{ds^{\star}}\Big|
\end{equation} by~(\ref{234}).
The curves $\gamma_i$ are absolutely continuous. It implies
\begin{eqnarray}\label{231}
\int\limits_{f(\gamma)}\rho^{\star}_i\,ds^{\star} & \geq &
\int\limits_{f(\gamma_i)}\rho^{\star}_i\,ds^{\star}\geq
\int\limits_{I_i^{\star}}\rho(\gamma_i(s^{\star}))\lambda_f(\gamma_i(s^{\star}))\,ds^{\star}\nonumber
\\
& = &
\int\limits_{0}^{l(\gamma_i)}\rho(\gamma_i(s))\lambda_f(\gamma_i(s))\Big|\frac{ds^{\star}}{ds}\Big|\,ds
\geq \int\limits_{\gamma_i}\rho\,ds
\end{eqnarray} by~(\ref{230}). If we show that
$\lim\limits_{i\to\infty}\int\limits_{f(\gamma)}\rho^{\star}_i\,ds^{\star}$
exists, then~(\ref{231}) and
$\lim\limits_{i\to\infty}\int\limits_{\gamma_i}\rho\,ds\geq 1$
will imply that the limit function $\rho^{\star}$ belongs to
$\mathcal F(f(\Gamma))$. We argue as follows
\begin{eqnarray}\label{232}
\int\limits_{f(\Omega)}\big(\rho^{\star}_i\big)^Q\,dy & = &
\int\limits_{f(\Omega)}\max\limits_{x\in f^{-1}(y)\cap
\Omega_i}\big(\rho(x)\lambda_f(x)\big)^Q\,dy\nonumber
\\ & \leq &
\int\limits_{f(\Omega)}\sum\limits_{x\in f^{-1}(y)\cap
\Omega_i}\big(\rho(x)\lambda_f(x)\big)^Q\,dy\nonumber
\\ & \leq &
K_I(f)\int\limits_{f(\Omega)}\Big(\sum\limits_{x\in f^{-1}(y)\cap
\Omega_i}\rho(x)^QJ^{-1}(x,f)\Big)\,dy=K_I(f)\int\limits_{\Omega_i}\rho^{Q}\,dx,
\end{eqnarray} by~(\ref{id}) and~(\ref{in2}).
Since for $i\geq k>0$
$$\int\limits_{f(\Omega)}|\rho^{\star}_i-\rho^{\star}_k|^Q\,dy\leq
K_I(f)\int\limits_{\Omega_i\setminus\Omega_k}\rho^Q\,dx\to
0\quad\text{as}\quad k\to\infty,$$ we deduce that there exists a
function $\rho^{\star}\in L_Q(f(\Omega))$ such that
$$\lim\limits_{i\to\infty}\int\limits_{f(\Omega)}|\rho^{\star}_i-\rho^{\star}|^Q\,dy=0.$$ A
result of~\cite{Fug} implies that there is a subsequence of
$\rho^{\star}_i$ (for the simplicity we use the same symbol
$\rho^{\star}_i$) with
$\lim\limits_{i\to\infty}\int\limits_{\gamma^{\star}}|\rho^{\star}_i-\rho^{\star}|\,ds^{\star}=0$
for $\gamma^{\star}\in \Gamma^{\star}\setminus\Gamma^0$ where
$M_Q(\Gamma^0)=0$. From here it follows that
$\int\limits_{\gamma^{\star}}\rho^{\star}\,ds^{\star}\geq 1$ for
$\gamma^{\star}\in
\Gamma^{\star}\setminus\big(\Gamma^0\cup\Gamma_1^{\star}\cup(\cup_{i}\Gamma_{0,i}^{\star})\big)$.
Moreover, the inequality
$\int\limits_{f(\Omega)}(\rho^{\star})^Q\,dy\leq
K_{I}\int\limits_{\Omega}\rho^Q\,dx$ holds. Finally, we conclude
$$M_Q(f(\Gamma))=M_Q\big(f(\Gamma)\setminus\big(\Gamma^0\cup\Gamma_1^{\star}\cup(\cup_{i}\Gamma_{0,i}^{\star})\big)\big) \leq
K_I(f)\int\limits_{\Omega}\rho^Q\,dx.$$ We obtain~(\ref{233})
taking infimum over all admissible functions $\rho$ for $\Gamma$.
\end{proof}

\subsection{Lifting of curves}

Let $f:\Omega\to\overline{\mathbb G}$ be continuous discrete and
open mapping of a domain $\Omega\in \mathbb G$. Let
$\beta:[a,b[\in\overline{\mathbb G}$ be a curve and let $x\in
f^{-1}(\beta(a))$. A curve $\alpha:\ [a,c[\to\Omega$ is called an
{\it $f$-lifting of $\beta$ starting at point $x$} if
\begin{itemize}
\item[1)]{$\alpha(a)=x$,}
\item[2)]{$f\circ\alpha=\beta\vert_{[a,c[}$,}
\end{itemize}
We say that a curve $\alpha:\ [a,c[\to\Omega$ is a {\it maximal}
$f$-lifting of $\beta$ starting at point $x$ if both~1),~2) and
the following property hold:
\begin{itemize}
\item[3)]{if $c<c^{\prime}<b$ then there does not exist a curve
$\alpha^{\prime}:\ [a,c^{\prime}[\to\Omega$ such that
$\alpha=\alpha^{\prime}\vert_{[a,c[}$ and
$f\circ\alpha^{\prime}=\beta\vert_{[a,c^{\prime}[}$.}
\end{itemize}

Let $f^{-1}(\beta(a))=\{x_1,\ldots,x_k\}$ and
$m=\sum\limits_{j=1}^k i(x_j,f)$. We say that
$\alpha_1,\ldots,\alpha_m$ is a {\it maximal essentially separate}
sequence of $f$-liftings of $\beta$ starting at the points
$x_1,\ldots,x_k$ if
\begin{itemize}
\item[1)]{each $\alpha_j$ is a maximal lifting of $f$,}
\item[2)]{$\card\{j:\alpha_j(a)=x_l\}=i(x_l,f)$, $1\leq l\leq k$,}
\item[3)]{$\card\{j:\alpha_j(t)=x\}\leq i(x,f)$ for all
$x\in\Omega$ and all $t$.}
\end{itemize}

Similarly, we define a maximal sequence of $f$-liftings
terminating at $x_1,\ldots,x_k$ if $f:\ ]b,a]\to\overline{\mathbb
G}$.

\begin{theorem}\label{Rick}
Let $f:\Omega\to\overline{\mathbb G}$ be a quasimeromorphic
mapping, $\beta: [a,b[\to\overline{\mathbb G}$ $(\beta:
]b,a]\to\overline{\mathbb G})$ be a curve, and let
$x_1,\ldots,x_k$ be distinct points in $f^{-1}(\beta(a))$. Then
$\beta$ has a maximal sequence of $f$-liftings starting
$($terminating$)$ at $x_1,\ldots,x_k$.
\end{theorem}

\begin{proof} The theorem is formulated for quasimeromorphic mappings but
actually this is a topological assertion. For the proof we refer
to~\cite{Rick1} where it is shown the local existence of
$f$-liftings. Since topological properties of $S=f^{-1}(\infty)$
coincides with topological properties of the pre-image of a finite
point of $\mathbb G$, we can apply the arguments of~\cite{Rick1}
almost verbatim. Then, we can show that the local existence of
maximal $f$-liftings implies the global existence.
\end{proof}

In the next statement we present a generalization of the
inequality of J.~V\"{a}is\"{a}l\"{a}. The V\"{a}is\"{a}l\"{a}
inequality is an essential tool on the study of value distribution
of quasimeromorphic mappings.

\begin{theorem}\label{Vasineq}
Let $f:\Omega \to\overline{\mathbb G}$ be a nonconstant
quasimeromorphic mapping, $\Gamma$ be a family of curves in
$\Omega$, $\Gamma^{\star}$ be a family in $\overline{\mathbb G}$
and $m$ be a positive integer such that the following is true. For
every locally rectifiable curve $\beta: \langle a,b \rangle\to
\overline{\mathbb G}$ in $\Gamma^{\star}$ there exist curves
$\alpha_1,\ldots\alpha_m$ in $\Gamma$ such that
\begin{itemize}
\item[1)]{$(f\circ\alpha_j)\subset\beta$ for all $j=1,\ldots,m$,}
\item[2)]{$\card\{j:\ \alpha_j(t)=x\}\leq i(x,f)$ for all
$x\in\Omega$ and for all $t\in\langle a,b \rangle$.}
\end{itemize} Then
$$M_Q(\Gamma^{\star})\leq
\frac{K_I(f)}{m}M_Q(\Gamma).$$
\end{theorem}

\begin{proof}
Let $\rho$ be an admissible function for a family $\Gamma$. If
$x\in \Omega\setminus B_f$ then $f$ is $\mathcal P$-differentiable
almost everywhere with strictly positive Jacobian $J(x,f)$. By
$F$, we denote a Borel set of measure zero containing all points
$x\in \Omega\setminus B_f$ where $f$ is not $\mathcal
P$-differentiable. Since $\mes (F\cup B_f)=\mes f(F\cup B_f)=0$ we
find a Borel set $E$ such that $f(F\cup B_f)\cup\{\infty\}\subset
E\subset\overline{\mathbb G}$ and $\mes(E)=0$. We define a
function $\rho^{\star}(y)$ on $f(\Omega)$ by the following way:
$$\rho^{\star}(y)=\left\{\array{cll}
&\frac{1}{m}\sum\limits_{x\in f^{-1}(y)\cap\Omega
}\rho(x)\lambda_f(x)\quad & \text{if}\quad y\in f(\Omega)\setminus
E,
\\
& \infty \quad &\text{if}\quad y\in E,
\\
& 0\quad & \text{if}\quad y\notin f(\Omega),
\endarray\right.$$ where $\lambda_f(x)=\Big(\min\limits_{|\xi|_0=1,\xi\in
V_1}|D_Hf(x)\xi|_0\Big)^{-1}$. The function $\rho^{\star}(y)$ is a
Borel nonnegative function. We show that $\rho^{\star}\in\mathcal
F(\Gamma^{\star})$.

Let $\gamma^{\star}\in\Gamma^{\star}$ and
$\alpha_1,\ldots,\alpha_m$ be as in conditions of the theorem.
Lemma~\ref{Pollem} implies that the family $\Gamma_0^{\star}$ of
curves $\gamma^{\star}$ for which $\alpha_1,\ldots,\alpha_m$ are
not absolutely continuous, vanishes: $M_Q(\Gamma_0^{\star})=0$.
The family $\Gamma_1^{\star}\subset\Gamma^{\star}$ where
$\int\limits_{\gamma^{\star}}\chi_E\,ds>0$, $E=f(F\cup
B_f)\cup\{\infty\}$, $\gamma^{\star}\in\Gamma_1^{\star}$, has also
the $Q$-module zero (it can be shown as in Theorem~\ref{K0ineq}).
Throughout the proof we restrict our attention to the family
$\Gamma^{\star}\setminus \big(\Gamma_0^{\star}\cup
\Gamma_1^{\star}\big)$.

Let, for the moment, suppose that $\gamma^{\star}\in
\Gamma^{\star}\setminus \big(\Gamma_0^{\star}\cup
\Gamma_1^{\star}\big)$ be a closed curve: $\gamma^{\star}:[a,b]\to
\G$. We correlate the parameterization for $\gamma^{\star}$ and
$\alpha_1,\ldots,\alpha_m$ as in Remark~\ref{r3}. We follow the
notations: $s^{\star}$ is the length arc parameter on
$\gamma^{\star}$, $s^{\star}\in I^{\star}=[0,l(\gamma^{\star})]$;
for each $\alpha_k$ the interval $I^{\star}_{k}$ is such that
$f(\alpha_k(s^{\star}))\subset\gamma^{\star}$ when $s^{\star}\in
I^{\star}_{k}$; $s$ is the length arc parameter on $\alpha_k$,
$k=1,\ldots,m$, and $s\in I_k$. Then by~(\ref{230}) we see
\begin{eqnarray}\label{284}
1\leq\int\limits_{\alpha_k}\rho\,ds & = &
\int\limits_{I_k}\rho(\alpha_k(s))\,ds=
\int\limits_{I^{\star}_k}\rho(\alpha_k(s^{\star}))\Big|
\frac{ds}{ds^{\star}}\Big|\,ds^{\star}\nonumber \\ & \leq &
\int\limits_{I^{\star}_k}\rho(\alpha_k(s^{\star}))\lambda_f(\alpha_k(s^{\star}))\,ds^{\star}
\end{eqnarray} for each curve $\alpha_k$, $k=1,\ldots,m$.

Set $K_{s^{\star}}=\{k:\ s^{\star}\in I^{\star}_{k}\}$. Then for
almost all $s^{\star}\in I^{\star}$ the points
$\alpha_k(s^{\star})$, $k\in K_{s^{\star}}$, are distinct points
in $f^{-1}(\gamma^{\star}(s^{\star}))$, $\gamma^{\star}\in
\Gamma^{\star}\setminus \big(\Gamma_0^{\star}\cup
\Gamma_1^{\star}\big)$. Therefore,
\begin{equation}\label{285}
\rho^{\star}(\gamma^{\star}(s^{\star}))
\geq\frac{1}{m}\sum\limits_{k=1}^{m}\rho(\alpha_k(s^{\star}))\lambda_f(\alpha_k(s^{\star}))\chi_{I^{\star}_k},
\end{equation} where $\chi_{I^{\star}_k}$ is the characteristic
function of $I^{\star}_k$. We conclude
from~\eqref{284},~\eqref{285} and
\begin{eqnarray*}
1 & \leq &
\frac{1}{m}\sum\limits_{k=1}^{m}\int\limits_{\alpha_k}\rho\,ds
\leq  \frac{1}{m}
\sum\limits_{k=1}^{m}\int\limits_{I^{\star}_k}\rho(\alpha_k(s^{\star}))\lambda_f(\alpha_k(s^{\star}))\,ds^{\star}
\\ & \leq & \int\limits_{I^{\star}}\frac{1}{m}
\sum\limits_{k=1}^{m}\rho(\alpha_k(s^{\star}))\lambda_f(\alpha_k(s^{\star}))\chi_{I^{\star}_k}\,ds^{\star}
\leq\int\limits_{\gamma^{\star}}\rho^{\star}\,ds^{\star}.
\end{eqnarray*}
that $\rho^{\star}\in\mathcal F\big(\Gamma^{\star}\setminus
(\Gamma_0^{\star}\cup \Gamma_1^{\star})\big)$. If the curve
$\gamma^{\star}$ is not closed we obtain the same result taking
supremum over all closed parts of $\gamma^{\star}$.

Now we estimate $M_Q(\Gamma^{\star})$. We choose an exhaustion of
$\Omega$ by measurable sets
$\Omega_1\subset\Omega_2\subset\ldots\subset\Omega$,
$\cup_i\Omega_i=\Omega$. Since $\mes(E)=0$, we obtain
\begin{eqnarray*}\label{n41}
\int\limits_{\G}\rho^{\star}(y)^Q\chi_{\Omega_i}\,dy
 & \leq &
\int\limits_{\G}\Big(\frac{1}{m}\sum\limits_{x\in
f^{-1}(y)\cap\Omega}\rho(x)\chi_{\Omega_i} \lambda_f(x)\Big)^Q\,dy
\\ & \leq &
\frac{1}{m}\int\limits_{\G}\sum\limits_{x\in
f^{-1}(y)\cap\Omega}\rho(x)^Q\chi_{\Omega_i}\lambda_f^Q(x)\,dy
\\
& \leq & \frac{K_{I}(f)}{m}\int\limits_{\G}\sum\limits_{x\in
f^{-1}(y)\cap\Omega}\rho(x)^Q\chi_{\Omega_i}J^{-1}(x,f)\,dy
\\ & \leq &
\frac{K_I(f)}{m}\int\limits_{\Omega_i}\rho(x)^Q\,dx\leq
\frac{K_I(f)}{m}\int\limits_{\Omega}\rho(x)^Q\,dx ,\end{eqnarray*}
by~(\ref{id}) and~(\ref{in2}). Letting $i\to \infty$, we deduce
$$M_Q(\Gamma^{\star})=M_Q\big(\Gamma^{\star}\setminus
(\Gamma_0^{\star}\cup\Gamma_1^{\star})\big)\leq
\int\limits_{\G}\rho^{\star}(y)^Q\,dy\leq
\frac{K_I(f)}{m}\int\limits_{\Omega}\rho(x)^Q\,dx,$$ and the proof
is complete.
\end{proof}

\begin{corollary}\label{35} Let $f$ be as in
Theorem~$\ref{Vasineq}$, $U$ be a normal domain for $f$ with
$m=N(f,U)$, $\Gamma^{\star}$ be a family of curves in $f(U)$,
$\Gamma$ be a family of curves in $U$ such that
$f\circ\alpha\in\Gamma^{\star}$ for any $\alpha\in\Gamma$. Then
$$M_Q(\Gamma^{\star})\leq \frac{K_I(f)}{m}M_Q(\Gamma).$$
\end{corollary}
\begin{proof} Let $\beta:[a,b)\to f(U)$ be a curve from
$\Gamma^{\star}$ and $\{x_1,\ldots,x_k\}=U\cap f^{-1}(\beta(a))$.
Since $U$ is a normal domain, we have
$\sum\limits_{l=1}^{k}i(x_l,f)=m$. Theorem~\ref{Rick} implies that
there exists a maximal sequence $\alpha_1,\ldots,\alpha_m$ of
$f$-liftings starting at $x_1,\ldots,x_k$ defining on $[a,b)$. We
have $\sum\limits_{x\in f^{-1}(y)\cap U}i(x,f)=m$ for every $y\in
f(U)$. Therefore, the condition~2) of Theorem~\ref{Vasineq} also
holds.

If $\beta$ is defined on an arbitrary interval $\langle
a,b\rangle$, then a slight change in the proof give the same
conclusion.
\end{proof}

\section{Relations between module and counting functions}~\label{relations}

>From now on, we restrict our considerations on $\mathbb H$-type
Carnot groups, for which the Heisenberg group is the simplest
example. For the definition of $\mathbb H$-type Carnot groups see
the example~3 and references therein. As was mentioned in the
introduction, the presented value distribution theory can be
extended to "polarizable" Carnot group introduced in the
work~\cite{BT}. We give the necessary definitions.

Let $\Omega\subset\mathbb G$ be a domain. A function $u\in
W^1_2(\Omega)$ is said to be {\it harmonic} if it is a weak
solution to the equation
\begin{equation}\label{0-Lap}
\Delta_0u=\diverg(\nabla_0u)=\sum\limits_{j=1}^{n_1}X_{1j}^2u=0.
\end{equation}
The norm $|\cdot|=u^{1/(2-Q)}$, that we use below, is associated
to the fundamental solution of~\eqref{0-Lap}. The fundamental
solution exists and it is unique by a result of
Folland~\cite[Theorem 2.1]{Foll}. Denote a set of characteristic
points by $\mathcal Z=\{0\}\cup\{x\in\mathbb
G\setminus\{0\}:\nabla_0|x|=0\}$.

\begin{theorem}\cite{BT}\label{BalTys} Let $S=S(0,1)=\{x\in\mathbb G: |x|=1\}$.
There exists a unique Radon measure $\sigma^{\star}$ on
$S\setminus\mathcal Z$ such that for all $u\in L^1(\mathbb G)$
\begin{equation}\label{270}\int\limits_{\mathbb
G}u(x)\,dx=\int\limits_{S\setminus\mathcal
Z}\int\limits_{0}^{\infty} u(\varphi(s,y))s^{Q-1}\,
ds\,d\sigma^{\star}(y),\end{equation} where $dx$ denotes the Haar
measure on $\mathbb G$.
\end{theorem}
Here $\varphi:(0,\infty)\times \mathbb G\setminus\mathcal Z\to
\mathbb G\setminus\mathcal Z$ is a flow of radial rectifiable
curves. The flow $\varphi$ satisfies the following properties:
\begin{itemize}\item[(i)]{$|\varphi(s,x)|=s|x|$ for $s>0$ and $x\in \mathbb G\setminus\mathcal Z$;}
\item[(ii)]{$\big|\frac{\partial\varphi(s,x)}{\partial
s}\big|_0=\frac{|x|}{\big|\nabla_0|x|\big|_0}$ and, in
particulary, is independent of $s$;}
\item[(iii)]{$J(x,\varphi(s,x))=s^Q$ for $s>0$ and $x\in \mathbb
G\setminus\mathcal Z$.}\end{itemize} We present the value of the
constant $\kappa(\mathbb G,p)$ from~\eqref{272} (see~\cite{BT}).
If we use the notation
$\upsilon(x)=\frac{\big|\nabla_0|x|\big|_0}{|x|}$, $x\in \mathbb
G\setminus\mathcal Z$, then
\begin{equation}\label{273}\kappa(\mathbb
G,p)=\int\limits_{S\setminus\mathcal
Z}\upsilon(y)^{p}\,d\sigma^{\star}(y).\end{equation} The set
$\mathcal Z\cap S$ has Hausdorff dimension at most
$N-2$~\cite{FW}, where $N$ is the topological dimension of the
group.

The polar coordinates can be introduced on any Carnot group for
any homogeneous norm $|\cdot|^{\prime}$. The integration formula
for $u\in L^1(\mathbb G)$ is of the form
\begin{equation}\label{271}\int\limits_{\mathbb
G}u(x)\,dx=\int\limits_{S^{\prime}}\int\limits_{0}^{\infty}
u(\delta_sy)s^{Q-1}\,ds\,d\sigma^{\prime}(y),\end{equation} where
$d\sigma^{\prime}$ is a Radon measure on
$S^{\prime}=S^{\prime}(0,1)=\{x\in\mathbb G: |x|^{\prime}=1\}$ and
$\delta_s$ is the dilation that was introduced in Section~1 (more
details see in~\cite{FS}). The formula~\eqref{270} differs
from~\eqref{271} in one important respect: the curves
$\varphi(s,y)$ have the finite Carnot-Carath\'{e}odory length.
This is not the case for the curves $\gamma(s,y)=\delta_sy$ in the
most situations. The pointed out difference allows us to employ
the standard family curves arguments. On the Heisenberg group the
polar coordinates with a rectifiable radial flow were described
in~\cite{KR}.

Let $f:\Omega\to\overline{\mathbb G}$, $\Omega\subset\mathbb G$,
be a quasimeromorphic mapping. For a point $y\in \overline{\mathbb
G}$ and for a Borel set $E\subset\Omega$ such that $\overline E$
is a compact in $\Omega$ we set
$$n(E,y)=\sum\limits_{x\in f^{-1}(y)\cap E}i(x,f).$$ In the case
$E=B(0,r)$ we use the notation $n(r,y)$.
\begin{lemma}\label{15}
The function $y\mapsto n(E,y)$ is upper semicontinuous.
\end{lemma}

\begin{proof}
Since a quasimeromorphic mapping is discrete, any point $x\in E$
has a normal neighborhood. If $U(x)$ is a normal neighborhood of
$x\in E$, we have $$i(x,f)=\mu(f(x),f,U(x))=\mu(f(z),f,U(x))\geq
i(z,f)$$ for any $z\in U(x)$. It shows that the function $x\mapsto
i(x,f)$ is upper semicontinuous.
\end{proof}

If $S(z,s)$ is a sphere in $\mathbb G$ we denote by
$\nu(E,S(z,s))$ the average of $n(E,y)$ over the sphere $S(z,s)$
with respect to a measure $\sigma=\upsilon^Q\sigma^{\star}$, where
$\upsilon(x)$ is a function from~\eqref{273} and $\sigma^{\star}$
is the Radon measure on a sphere defined in Theorem~\ref{BalTys}.
The measure $\sigma$ is absolutely continuous with respect to the
measure $\sigma^{\star}$. In particular, we denote by $\nu(r,s)$
the average of $n(r,y)$ over $S(0,s)$. Hence
$$\nu(r,s)=\frac{1}{\sigma(S\setminus\mathcal Z)}\int\limits_{S\setminus\mathcal Z}n(r,\varphi(s,y))d\sigma(y),$$
where $\sigma(S\setminus\mathcal Z)=\sigma(S(0,1)\setminus\mathcal
Z)$ is the measure of the unit sphere coinciding with the constant
$\kappa(\mathbb G,Q)$. We observe that value of $\upsilon(x)$ is
invariant under the left translation. In fact, since
$\upsilon(x)=\big|\frac{\partial\varphi(s,x)}{\partial
s}\big|^{-1}_0=\big|\dot\varphi(s,x)\big|_0^{-1}$ we need to show
\begin{equation}\label{280}\big|\dot\varphi(s,x)\big|_0=\big|\dot{w\varphi}(s,x)\big|_0,\end{equation}
where $w\varphi(s,x)$ is the image of $\varphi(s,x)$ under the
left translation by the element $w$. We have
$(w\varphi)_{1j}(s,x)=\varphi_{1j}(s,x)+w_{1j}$, $j=1,\ldots,n_1$.
Thus $\dot{w\varphi}_{1j}(s,x)=\dot\varphi_{1j}(s,x)$. The curves
$\varphi(s,x)$, $w\varphi(s,x)$ are rectifiable. Then
\begin{eqnarray*}
\dot\varphi(s,x) & = &
\sum\limits_{j=1}^{n_1}\dot\varphi_{1j}(s,x)X_{1j}(\varphi(s,x)),
\\
\dot{w\varphi}(s,x) & = &
\sum\limits_{j=1}^{n_1}\dot{w\varphi}_{1j}(s,x)X_{1j}(w\varphi(s,x))
=\sum\limits_{j=1}^{n_1}\dot\varphi_{1j}(s,x)X_{1j}(w\varphi(s,x)).
\end{eqnarray*}
Since the left invariant basis is orthonormal at an arbitrary
point of $\mathbb G$ we deduce the necessary result. Therefore,
$\kappa(\mathbb G,Q)=\sigma(S\setminus\mathcal
Z)=\sigma(S(w,1)\setminus\mathcal Z)$, $w\in\mathbb G$.

We use below Definition~\ref{6} of an isolated essential
singularity of $f$.

\begin{lemma}\label{110}
Let $f:\mathbb G\to\overline{\mathbb G}$ be a quasimeromorphic
mapping with $\{\infty\}$ as an essential singularity and
$Y=S(w,s)$ is a sphere in $\overline{\mathbb G}$, then
$$\lim\limits_{r\to \infty}\nu(r,Y)=\infty. $$
\end{lemma}

\begin{proof}
By Lemma~\ref{19}, there is a $\mathcal F_{\sigma}$-set
$E\in\overline{\mathbb G}$ of the $Q$-capacity zero such that
$$N(y,f,\mathbb C\overline B(0,r)\setminus\{\infty\})=\infty\quad\text{for all}\quad
y\in\overline{\mathbb G}\setminus E\quad\text{and for all
}\quad{r>0}.
$$
Let $F_k(r)=\{y\in Y\setminus\mathcal Z:\ n(r,y)\geq k\}$ and
$E_0=(Y\setminus\mathcal Z)\cap E$. We claim that $\sigma(E_0)=0$
if $\capac_Q(E_0)=0$. Results of~\cite{ChV2,KV} imply that
$\mathcal H^{|\cdot|}_{\alpha}(E_0)=0$ for any $\alpha>0$. It
follows $\mathcal H^{|\cdot|}_{Q-1}(E_0)=0$. Here
$H^{|\cdot|}_{\alpha}$ is an $\alpha$-dimensional Hausdorff
measure with respect to the homogeneous norm $|\cdot|$. If we
denote the Riemannian area element on $Y\setminus\mathcal Z$ by
$dA$, then the connection between the Riemannian area element and
the Hausdorff measure is expressed by the formula
$dA=\frac{\big|\nabla|x|\big|}{\big|\nabla_0|x|\big|}d\mathcal
H^{|\cdot|}_{Q-1}$ (see, for instance,
~\cite[Proposition~4.9]{H1}). It follows that $A(E_0)=0$. The
Radon measure $\sigma^{\star}$ is absolutely continuous with
respect to $dA$, that gives $\sigma^{\star}(E_0)=\sigma(E_0)=0$.

We continue to estimate $\lim\limits_{r\to\infty}\nu(r,Y)$.
\begin{eqnarray*}\lim\limits_{r\to\infty}\nu(r,Y) & \geq &
\frac{1}{\sigma(Y\setminus\mathcal
Z)}\lim\limits_{r\to\infty}\int\limits_{F_k(r)}k\,d\sigma(y)=\frac{k}{\sigma(Y\setminus\mathcal
Z)} \lim\limits_{r\to\infty}\sigma(F_k(r))\\ & \geq &
\frac{k}{\sigma(Y\setminus\mathcal Z)}\sigma((Y\setminus\mathcal
Z)\setminus E_0))=k\end{eqnarray*} for every $k>0$. The lemma is
proved.
\end{proof}

>From now on, we use the symbol $S$ to denote the unit sphere
$S(0,1)$ centered at the identity of the group. Let us choose a
parametrization of radial curves $\varphi(s,x)$ in such a way that
$|\varphi(1,y)|=|y|=1$ for $y\in S\setminus \mathcal Z$. We fix
the notation $\varphi_s(y)$, $y\in S\setminus \mathcal Z$  for
such curve. Then any $x\in\mathbb G\setminus\mathcal Z$ can be
obtained as an image of $y\in S\setminus \mathcal Z$ under the map
$y\mapsto\varphi_s(y)$ for some $s>0$. Therefore
$|x|=|\varphi_s(y)|=s$.

\begin{lemma}\label{126} Let $E$ be a Borel set on $S\setminus\mathcal Z$ and let
$C$ be the cone $\{x\in\mathbb G: x/|x|\in E\}$. Set $\Gamma_E$ be
the family of all curves $\varphi_s(y):[a,b]\to\mathbb G$, $y\in
E$. Then
$$\sigma(E)\Big(\ln\frac{b}{a}\Big)^{1-Q}=M_Q(\Gamma_E).$$
\end{lemma}

\begin{proof} Let $\rho\in \mathcal F(\Gamma_E)$ and $\varphi_s(y):\ [a,b]\to\mathbb G$ be a radial
curve, introduced before Lemma~\ref{126}. Recall that
$\big|\frac{\partial\varphi_s(y)}{\partial
s}\big|_0=\frac{\big|\nabla_0 |y|\big|_0}{|y|}=\upsilon(y)^{-1}$.
The H\"{o}lder inequality implies
\begin{eqnarray*}
1 & \leq &
\Big(\int\limits_{\varphi_s(y)}\rho\,dt\Big)^Q=\Big(\int\limits_{a}^b\rho(\varphi_s(y))
\big|\frac{\partial\varphi_s(y)}{\partial s}\big|_0\,ds\Big)^Q \\
& \leq &
\int\limits_a^b\rho(\varphi_s(y))^Q\upsilon(y)^{-Q}s^{Q-1}\,ds
\Big(\int\limits_a^b\frac{ds}{s}\Big)^{Q-1}
\\ & = &
\Big(\ln\frac{b}{a}\Big)^{Q-1}\int\limits_a^b\rho(\varphi_s(y))^Q\upsilon(y)^{-Q}s^{Q-1}\,ds.
\end{eqnarray*}
Integrating over $y\in E$ with respect to the measure $\sigma$
yields
$$\sigma(E)\leq\Big(\ln\frac{b}{a}\Big)^{Q-1}\int\limits_E\int\limits_a^b\rho(\varphi_s(y))^Qs^{Q-1}\,ds\,d\sigma^{\star}(y)
=\Big(\ln\frac{b}{a}\Big)^{Q-1}\int\limits_C\rho(x)\,dx.$$ Taking
the infimum over all $\rho\in\mathcal F(\Gamma_E)$ we obtain the
inequality $\sigma(E)\Big(\ln\frac{b}{a}\Big)^{1-Q}\leq
M_Q(\Gamma_E)$.

To deduce the reverse inequality, we take the function
$$
\rho(x)=\left\{\array{cl} \frac{\upsilon(x)}{|x|\ln b/a}\quad &
\text{if $x\in C$},
\\
0\quad & \text{otherwise}.
\endarray\right.
$$
We claim that $\rho$ is admissible for $\Gamma_E$. Indeed, since
\begin{equation}\label{286}\upsilon(\varphi_s(y))^{-1}=\Big|\frac{\partial\varphi_s(\varphi_s(y))}{\partial
s}\Big|_0=\frac{\big|\nabla_0|\varphi_s(y)|\big|_0}{|\varphi_s(y)|}=\frac{\big|\nabla_0
s|y|\big|_0}{s|y|}=\upsilon(y)^{-1},\end{equation} we estimate
\begin{eqnarray*}\int\limits_{\varphi_s(y)}\rho(s)\,ds & = &
\int_a^b\rho((\varphi_s(y)))\Big|\frac{\partial\varphi_s(y)}{\partial
s}\Big|_0\,ds\\
& = &
\Big(\ln\frac{b}{a}\Big)^{-1}\int_a^b\frac{\upsilon(\varphi_s(y))}{|\varphi_s(y)|}\upsilon^{-1}(y)\,ds
=\Big(\ln\frac{b}{a}\Big)^{-1}\int_a^b\frac{ds}{s}=1.\end{eqnarray*}
Thus, by~\eqref{286} we see \begin{eqnarray*} M_Q(\Gamma_E) & \leq
& \int\limits_C\rho(x)^Q\,dx=
\Big(\ln\frac{b}{a}\Big)^{-Q}\int\limits_E\int\limits_a^b
\frac{\upsilon(\varphi_s(y))^Q}{|\varphi_s(y)|^Q}s^{Q-1}\,ds\,d\sigma^{\star}
\\ & = &
\Big(\ln\frac{b}{a}\Big)^{-Q}\int\limits_E\int\limits_a^b\frac{ds}{s}d\sigma
=\sigma(E)\Big(\ln\frac{b}{a}\Big)^{1-Q}.
\end{eqnarray*}
The proof is completed.
\end{proof}

We would like to compare averages of $n(r,y)$ over two distinct
concentric spheres on the Carnot group.

\begin{proposition}\label{51}
Let $m: S\setminus\mathcal Z\to\mathbb Z$, be a nonnegative
integer valued Borel function, $t,s>0$, and, for $y\in
S\setminus\mathcal Z$, $\beta_y(\tau)$ be a radial curve $\omega
\varphi_{\tau}(y)$ connecting the spheres $S(\omega,s)$ and
$S(\omega,t)$. Suppose $\Gamma$ is a family consisting of $m(y)$
essentially separate partial $f$-liftings of each $\beta_y$ when
$y$ runs over $S\setminus\mathcal Z$. Then
\begin{equation}\label{27}
\int\limits_{S\setminus\mathcal Z}m(y)\,d\sigma(y)\leq
K_I(f)\Big|\ln \frac{t}{s}\Big|^{Q-1}M_Q(\Gamma).\end{equation}
\end{proposition}

\begin{proof} We use the following notations: $E_k=\{y\in S\setminus\mathcal Z:\ m(y)=k\}$,
$k=0,1,2,\ldots$, $\Gamma^{\star}_k=\{\beta_y:\ y\in E_k\}$,
$\Gamma_k=\{f\text{-liftings of}\quad
\beta_y\in\Gamma_k^{\star}\}$. Theorem~\ref{Vasineq} and
Lemma~\ref{126} imply
\begin{equation}\label{28} k\sigma(E)\Big|\ln\frac{t}{s}\Big|^{1-Q}\leq kM_Q(\Gamma_k^{\star})\leq K_I(f)
M_Q(\Gamma_k).
\end{equation}
Since the curves in $\Gamma_k$ are separate, we deduce~(\ref{27})
summing~(\ref{28}) over $k$.
\end{proof}

\begin{theorem}\label{42}
Let $f: \Omega\to\overline{\mathbb G}$ be a quasimeromorphic
mapping, $\varrho>r>0$, and $t,s>0$. If $\overline
B(0,\varrho)\in\Omega$, then
\begin{equation}\label{29}
\nu(\varrho,S(\omega,t)) \geq \nu(r,S(\omega,s)) - K_I(f)\Big|\ln
\frac{t}{s}\Big|^{Q-1}\Big(\ln\frac{\varrho}{r}\Big)^{1-Q}
\end{equation}
for any $\omega\in\mathbb G$.
\end{theorem}

\begin{proof}
We may assume that $t>s>0$ and we set
$m(y)=\max\{0,n(r,\omega\varphi_s(y))-n(\varrho,\omega\varphi_t(y))\}$
for a point $y\in S\setminus\mathcal Z$. If $m(y)>0$, then there
exist at least $m(y)$ maximal $f$-liftings of
$\beta_y=\omega\varphi_{\tau}(y)$ starting in $\overline B(0,r)$
and terminating on $\partial B(0,\varrho)$ which are essentially
separate. We denote by $\Gamma$ the set of these $f$-liftings and,
making use of Proposition~\ref{51} and~\eqref{272}, deduce
\begin{eqnarray}\label{238}\int\limits_{S\setminus\mathcal Z}m(y)\,d\sigma(y) & \leq & K_I(f)
\Big(\ln\frac{t}{s}\Big)^{Q-1}M_Q(\Gamma)
\\ & \leq &
K_I(f)\kappa(\mathbb
G,Q)\Big(\ln\frac{t}{s}\Big)^{Q-1}\Big(\ln\frac{\varrho}{r}\Big)^{1-Q}.\nonumber
\end{eqnarray}
If $E=\{y\in S\setminus\mathcal Z:\ m(y)>0\}$, then
\begin{eqnarray*}
\int\limits_{S\setminus\mathcal
Z}n(\varrho,\omega\varphi_t(y))\,d\sigma(y) & = &
\int\limits_{(S\setminus\mathcal Z)\setminus
E}n(\varrho,\omega\varphi_t(y))\,d\sigma(y)+
\int\limits_{E}n(\varrho,\omega\varphi_t(y))\,d\sigma(y)
\\
& \geq & \int\limits_{(S\setminus\mathcal Z)\setminus
E}n(r,\omega\varphi_s(y))\,d\sigma(y)
\\ & +&
\int\limits_{E}\big(n(r,\omega\varphi_s(y))-m(y)\big)\,d\sigma(y)
\\
& \geq & \int\limits_{S\setminus\mathcal
Z}n(r,\omega\varphi_s(y))\,d\sigma(y)
-\int\limits_{S\setminus\mathcal Z}m(y)\,d\sigma(y).
\end{eqnarray*}
Dividing by $\kappa(\mathbb G,Q)$ and using~(\ref{238}), we finish
the proof.
\end{proof}

\begin{lemma}\label{111}
Let $\varrho>r>0$, $\theta>1$, and $f:\Omega\to\overline{\mathbb
G}$ be a quasimeromorphic mapping of a domain $\Omega$ with
$\overline B(0,\theta\varrho)\subset\Omega$. Let $a_1,\ldots,a_q$,
$q\geq 2$, be distinct finite points in $\mathbb G$. Set
$\sigma_m=\frac{1}{4}\min\limits_{i\neq j}\dist(a_i,a_j)$ and
$0<s<t\leq\sigma_m$. Assume that $F_1,\ldots, F_\lambda$,
$2\leq\lambda\leq q$, are disjoint compact sets in $\overline
B(0,\varrho)$ such that $f(F_j)\in \overline B(a_j,s)$ for each
$j\leq \lambda$, and $F_1,\ldots, F_\lambda$ intersect spheres
$S(0,\tau)$ for almost all $\tau\in[r,\rho]$. Then there are
positive constants $b_1$ and $b_2$ depending on $Q,\theta$ only
such that
\begin{equation}\label{210}
\big(M_Q(\Gamma_j)-b_1K_O(f)K_I(f)\big)\Big(\ln\frac{t}{s}\Big)^{Q-1}\leq
b_2K_O(f)\nu\big(\theta\varrho,S(a_j,t)\big),
\end{equation}
for all $j$. Here $\Gamma_j$ is a family of locally rectifiable
curves connecting $F_j$ with $\bigcup\limits_{i\neq j}F_i$ in
$B(0,\varrho)\setminus \overline B(0,r)$.
\end{lemma}

\begin{proof}
Fix $j\leq\lambda$ and assume that $a_j=0$. We put
\begin{equation}\label{213}\rho(z)=\left\{\array{cl}
\frac{\upsilon(z)}{|z|\ln t/s}\quad &\text{if}\quad z\in
B(0,t)\setminus B(0,s),
\\
0 & \text{elsewhere}.
\endarray\right.\end{equation}
The proof of Lemma~\ref{126} shows that the function $\rho(z)$ is
admissible for the family $f(\Gamma_j)$. Then
\begin{eqnarray}\label{211} M_Q(\Gamma_j) & \leq &
K_O(f)\int\limits_{\mathbb
G}\rho(z)^{Q}n(\varrho,z)\,dz=K_O(f)\kappa(\mathbb
G,Q)\Big(\ln\frac{t}{s}\Big)^{-Q}
\int\limits_{s}^{t}\nu(\varrho,\tau)\frac{d\tau}{\tau}\nonumber
\\
& \leq & K_O(f)\kappa(\mathbb G,Q)\nu(\theta
\varrho,t)\Big(\ln\frac{t}{s}\Big)^{1-Q}
+K_O(f)K_I(f)\kappa(\mathbb G,Q)(\ln\theta)^{1-Q}
\end{eqnarray} by Remark~\ref{r2} and
estimate~(\ref{29}) for $\tau<t$. We derive the lemma
from~(\ref{211}) with constants $b_1=\kappa(\mathbb
G,Q)(\ln\theta)^{1-Q}$ and $b_2=\kappa(\mathbb G,Q)$.
\end{proof}

\begin{corollary}\label{33}
Under the assumptions of Lemma~$\ref{111}$, there exist positive
constants $b_1$ and $b_2$ depending  on $Q,\theta$ only such that
$$\big(\ln\frac{\varrho}{r}-b_1K_O(f)K_I(f)\big)\Big(\ln\frac{t}{s}\Big)^{Q-1}\leq
b_2K_O(f)\nu\big(\theta\varrho,S(a_j,t)\big)$$ for arbitrary
$j=1,\ldots,\lambda$.
\end{corollary}

\begin{proof} Let $\Gamma_j$ be as in the condition of Lemma~\ref{111}.
Then the following estimate $M_Q(\Gamma_j)\geq
c(Q)\ln\frac{\rho}{r}$ is true for arbitrary
$j=1,\ldots,\lambda$~\cite{H}. We deduce the desired result
applying this inequality to the left hand side of~\eqref{210}.
Here $b_1=c^{-1}(Q)\kappa(\mathbb G,Q)(\ln\theta)^{1-Q}$ and
$b_2=c^{-1}(Q)\kappa(\mathbb G,Q)$.
\end{proof}

We would like to obtain an analogue of Lemma~\ref{111} when one of
the points $a_j$ is~$\{\infty\}$.

\begin{lemma}\label{112}
Let $\varrho>r>0$, $\theta>1$, and $f:\Omega\to\overline{\mathbb
G}$ be a quasimeromorphic mapping of a domain $\Omega$ with
$\overline B(0,\theta\varrho)\subset\Omega$. Let $a_0=\infty$ and
$a_1,\ldots,a_q$ be distinct finite points in $\mathbb G$. Set
$\sigma_M=4\max\limits_{1\leq j\leq q}\{|a_j|,1\}$ and
$\sigma_m=\frac{1}{4}\min\limits_{1\leq i\neq j\leq q
}\dist(a_i,a_j)$. Assume that $F_0,\ldots, F_\lambda$,
$1\leq\lambda\leq q$, are disjoint compact sets in $\overline
B(0,\varrho)$ such that $f(F_0)\in\mathbb C \overline B(0,t)$ for
$t>s>\sigma_M$, $f(F_j)\in \overline B(a_j,s^{\prime})$ for each
$1\leq j\leq \lambda$ and $0<s^{\prime}<t^{\prime}<\sigma_m$. We
suppose that $F_0,F_1,\ldots, F_\lambda$ intersect spheres
$S(0,\tau)$ for almost all $\tau\in[r,\rho]$. Then there are
positive constants $b_1$ and $b_2$ depending only on $Q,\theta$
such that
\begin{equation}\label{212}
\big(M_Q(\Gamma_0)-b_1K_O(f)K_I(f)\big)\Big(\ln\frac{t}{s}\Big)^{Q-1}\leq
b_2K_O(f)\nu\big(\theta\varrho,S(0,s)\big).
\end{equation}
Here $\Gamma_0$ is the family of locally rectifiable curves in
$B(0,\varrho)\setminus \overline B(0,r)$ that connect the
compact~$F_0$ to~$\bigcup\limits_{1\leq i\leq\lambda}F_i$.
\end{lemma}

\begin{proof}
The function~(\ref{213}) for $s,t$ such that $\sigma_M<s<t<\infty$
is admissible for the family $f(\Gamma_0)$. We have
\begin{eqnarray*}
\nu(\varrho,\tau) & \leq & \nu(\theta
\varrho,s)+K_I(f)\Big(\frac{|\ln\frac{\tau}{s}|}{\ln\theta}\Big)^{Q-1}
\leq  \nu(\theta\varrho,s) +
\frac{K_I(f)}{(\ln\theta)^{Q-1}}\Big|\ln\frac{t}{s}\Big|^{Q-1}
\end{eqnarray*}
for all $\tau\in(s,t)$ by Theorem~\ref{42}. Then, using the same
estimates as in~\eqref{211}, we finish the proof as in
Lemma~\ref{111} with $b_1=\kappa(\mathbb G,Q)(\ln\theta)^{1-Q}$
and $b_2=\kappa(\mathbb G,Q)$.
\end{proof}

\begin{corollary}\label{34}
Under the assumptions of Lemma~$\ref{112}$, there exist positive
constants $b_1$ and $b_2$ depending on $Q,\theta$ only such that
$$\big(\ln\frac{\varrho}{r}-b_1K_O(f)K_I(f)\big)\Big(\ln\frac{t}{s}\Big)^{Q-1}\leq
b_2K_O(f)\nu\big(\theta\varrho,S(0,t)\big),$$
\end{corollary}

\begin{proof}
We argue as in Corollary~\ref{33}.
\end{proof}

The average $\nu(r,s)$ over a sphere of radius $s$ is a
discontinuous function in the variable $r$. We need an auxiliary
continuous function $A(r)$ related to a quasimeromorphic mapping
$f:\Omega\to\overline{\mathbb G}$. In the classical value
distribution theory for analytic functions, it is considered the
integral $A(r)=\int_En(r,y)\,d\mu(y)$ with respect to a
nonnegative measure distributed over a closed set $E$ in a target
plane of an analytic function. Then $A(r)$ has a
geometric-physical significance: it equals the total mass
distributed over the Riemann surface of the analytic function onto
which this function maps the disk $B(0,r)$. In the Euclidean space
$\mathbb R^n$, S.~Rickman employed the $n$-dimensional normalized
spherical measure $\frac{c(n)\,dy}{(1+|y|^2)^n}$ as a measure
$\mu(y)$ which is invariant under conformal mappings. In this
case, $A(r)$ has a geometric meaning also. In the arguments below,
our version of $A(r)$ is used as an auxiliary tool only.

We define $A(r)$ as
\begin{equation}\label{239} A(r)=\frac{2Q}{\pi\kappa(\mathbb G,Q)}\int\limits_{\mathbb
G}\frac{n(r,y)\upsilon(y)^Q}{1+|y|^{2Q}}\,dy
=\frac{2Q}{\pi\kappa(\mathbb
G,Q)}\int\limits_{B(0,r)}\frac{J(x,f)\upsilon(f(x))^Q}{1+|f(x)|^{2Q}}\,dx.\end{equation}
In the definition of $A(r)$ we have taken into account that
$\mes(B_f)=\mes(f(B_f))=0$.

\begin{lemma}\label{113}
If $\theta>1$, $r>0$, $\overline B(0,\theta  r)\subset\Omega$, and
$Y=S(\omega,t)$ is a sphere with radius $t>0$, then
\begin{eqnarray}\label{214}
\nu\Big(\frac{r}{\theta},Y\Big) & - &
\frac{K_I(f)c_1}{(\ln\theta)^{Q-1}}\big(|\ln t|^{Q-1}+c_0\big)
\leq \ A(r)\ \nonumber
\\ & \leq & \nu(\theta
r,Y)  +  \frac{K_I(f)c_1}{(\ln\theta)^{Q-1}}\big(|\ln
t|^{Q-1}+c_0\big),
\end{eqnarray}
where $c_0$ and $c_1$ some positive constants depending on $Q$
only.
\end{lemma}

\begin{proof}
Writing~(\ref{29}) in the form
\begin{equation}\label{287}\nu(\theta r,t) \geq \nu(r,s)
-\frac{K_I(f)|\ln\frac{t}{s}|^{Q-1}}{(\ln\theta)^{Q-1}},\end{equation}
multiplying both sides of~\eqref{287} by
$\frac{s^{Q-1}}{1+s^{2Q}}$, and integrating from $0$ to $\infty$,
we obtain
\begin{equation}\label{215} \nu(\theta r,t)
\int\limits_{0}^{\infty}\frac{s^{Q-1}}{1+s^{2Q}}\,ds\geq
\int\limits_{0}^{\infty}\frac{\nu(r,s)s^{Q-1}}{1+s^{2Q}}\,ds
-\frac{K_I(f)
}{(\ln\theta)^{Q-1}}\int\limits_{0}^{\infty}\frac{|\ln\frac{t}{s}|^{Q-1}s^{Q-1}}{1+s^{2Q}}\,ds.
\end{equation}
The integral in the left hand side of~(\ref{215}) equals
$\frac{\pi}{2Q}$. The first integral in the right hand side gives
\begin{eqnarray*}\int\limits_{0}^{\infty}\frac{\nu(r,s)s^{Q-1}}{1+s^{2Q}}\,ds
& = &  \frac{1}{\kappa(\mathbb
G,Q)}\int\limits_{S\setminus\mathcal
Z}\int\limits_0^{\infty}\frac{n(r,\varphi_s(y))\upsilon(\varphi_s(y))^Qs^{Q-1}}{1-|\varphi_s(y)|^{2Q}}ds\,d\sigma^{\star}(y)
\\ & = & \frac{1}{\kappa(\mathbb G,Q)}\int\limits_{\mathbb
G}\frac{n(r,z)\upsilon(z)^Q}{1+|z|^{2Q}}\,dz=\frac{\pi}{2Q}A(r).\end{eqnarray*}
To estimate the second one we use the inequality
$$\Big|\ln\frac{t}{s}\Big|^{Q-1}\leq 2^{Q-2}(|\ln t|^{Q-1}+|\ln
s|^{Q-1})$$ and the finiteness of the integral
$\int\limits_{0}^{\infty}\frac{|\ln
s|^{Q-1}s^{Q-1}}{1+s^{2Q}}\,ds$. To simplify notation, we write
$\frac{\pi}{2Q}c_0$ for the last integral. Thus
\begin{equation}\label{216}\frac{K_I(f)
}{(\ln\theta)^{Q-1}}\int\limits_{0}^{\infty}\frac{|\ln\frac{t}{s}|^{Q-1}s^{Q-1}}{1+s^{2Q}}\,ds
\leq
\frac{K_I(f)2^{Q-2}}{(\ln\theta)^{Q-1}}\Big(\frac{\pi}{2Q}|\ln
t|^{Q-1}+\frac{\pi}{2Q}c_0\Big).\end{equation} Joining the
estimates of all integrals, we get the right hand side
of~(\ref{214}) with $c_1=2^{Q-2}$.

To obtain the left hand side of~(\ref{214}) we multiply
$$\nu(r,t) \geq \nu(r/ \theta ,s) -\frac{K_I(f)
|\ln\frac{t}{s}|^{Q-1}}{(\ln\theta)^{Q-1}}$$ by
$\frac{s^{Q-1}}{1+s^{2Q}}$, and integrate the result from $0$ to
$\infty$. Arguing as above, we prove the lemma.
\end{proof}

\begin{corollary}\label{36} Let $f:\mathbb G\to\overline{\mathbb G}$ be a quasimeromorphic
mapping with $\{\infty\}$ as an essential singularity and let
$A(r)$ be as in~$(\ref{239})$. Then $\lim\limits_{r\to
\infty}A(r)=\infty$.
\end{corollary}

\begin{proof}
The corollary follows from Lemmas~\ref{110} and~\ref{113}.
\end{proof}

We need an inequality that will permit us to compare the average
over spheres with different centers and radii.

\begin{lemma}\label{117}
If $\vartheta>1$, $r>0$, $\overline B(0,\vartheta r)\subset
\Omega$, and if $Z=S(\omega,u)$ and $Y=S(w,v)$ are spheres in
$\overline{\mathbb G}$ with radii $u$ and $v$ then
$$\nu(\vartheta r,Z)\geq\nu(r,Y)-\frac{c_1K_I(f)}{(\log\vartheta)^{Q-1}}
\big(|\ln u|^{Q-1}+|\ln v|^{Q-1}+2c_0\big)$$ where $c_1=2^{2Q-3}$.
\end{lemma}

\begin{proof}
To show Lemma~\ref{117} we replace $\theta$ by $\vartheta^{1/2}$
and $r$ by $r\vartheta^{1/2}$ in the inequality~\eqref{214}.
\end{proof}

\section{Auxiliary statements for the proof of the main theorem}\label{mt}

Recall the statement of Theorem~\ref{45}. In the formulation we
use the notation $u_{+}=\max\{0,u\}$.

{\it Let $f:\ \mathbb G\to\overline{\mathbb G}$ be a nonconstant
$K$-quasimeromorphic mapping. Then there exists a set
$E\subset[1,\infty)$ and a constant $C(Q,K)<\infty$ such that
\begin{equation*}
\lim\limits_{r\to \infty}\sup\limits_{r\notin E}\sum\limits_{j=0
}^{q}\Big(1-\frac{n(r,a_j)}{\nu(r,1)}\Big)_+\leq
C(Q,K)\quad\text{with}\quad \int\limits_{E}\frac{dr}{r}<\infty,
\end{equation*}
whenever $a_0,a_1,\ldots,a_q$ are distinct points in
$\overline{\mathbb G}$. }

To prove our main result we need to realize the following steps.
\begin{itemize}
\item[1.]{To set up the exceptional set $E\subset[1,\infty)$
mentioned in Theorem~\ref{45}.} \item[2.]{To construct a special
decomposition of the ball $B(0,s)$ to some smaller sets $U_i$,
$i=1,2,\ldots,p$. The sets $U_i$ constitute a finite-to-one
covering: there is $M>0$ such that $\sum_i\chi_{U_i}\leq M$ for
all $x\in B(0,s)$. Here $\chi_{U_i}$ is the characteristic
function of $U_i$.} \item[3.]{To apply Lemmas~\ref{111}
and~\ref{112}. For realizing this, we consider a chain of
increasing balls $U_i\subset V_i\subset W_i\subset X_i\subset
Y_i\subset Z_i\subset\overline B(0,s^{\prime})$ having finite-to
one covering property. Then we calculate the number of
$f$-liftings of radial curves connecting points $z\in
S(0,\sigma_M)$ with $a_j\in B(0,\sigma_M/2)$, $j=1,\ldots,q$, or
$a_0=\infty$.} \item[4.]{To estimate $Q$-module of families
related to these liftings in the spherical ring type domains
$V_i\setminus\overline U_i$ and $W_i\setminus\overline V_i$}
\item[5.]{To sum obtained results over all $U_i$ obtaining by this
way a bound from above for the sum $\sum_{j\in
J\cup\{0\}}\Delta_j$ (the definition of $\Delta_j$ see at the end
of Subsection~\ref{Constr}).}
\end{itemize}

\subsection {Construction of the exceptional set}\label{Constr}

In the present section we construct the exceptional set $E\in
[1,\infty)$ mentioned in Theorem~\ref{45}.

\begin{theorem}\label{44}
There exists a set $E\in[1,\infty)$ such that
$\int\limits_{E}\frac{dr}{r}<\infty$ and the following is true: if
$0<\varepsilon_0<1/5$ and we write
$$
s^{\prime}=s+\frac{s}{\varepsilon_0(A(s))^{\frac{1}{Q-1}}}\qquad\text{for}\quad
s>0,$$ then this is an increasing function
$\omega:[0,\infty)\to[D_{\varepsilon_0},\infty]$ such that
\begin{equation}\label{218}
\Big|\frac{\nu(s,Y)}{A(s^{\prime})}-1\Big|<\varepsilon_0
\end{equation}
and
\begin{equation}\label{219}
\frac{\nu(s,Y)}{\nu(s^{\prime},Y)}>1-\varepsilon_0
\end{equation} hold whenever $Y=S(w,t)$ is a sphere of
radius $t$, and $s^{\prime}\in[\omega(|\ln t|),\infty]\setminus
E$. Here $D_{\varepsilon_0}>0$,
$A(D_{\varepsilon_0})>\varepsilon_0^{1-Q}$.
\end{theorem}

\begin{proof}
First, we construct a set $E$. Set
$\phi(r)=m^{-2}A(r)^{\frac{1}{Q-1}}$ for each integer $m\geq 2$.
We can choose $r_0^{\prime\prime}=r_0^{\prime\prime}(m)\geq 1$
such that $\phi(r_0^{\prime\prime})\geq 1$, by Corollary~\ref{36}.
Let $F_m=\bigl\{r>r_0^{\prime\prime}:\
A\bigr(r+\frac{2r}{\phi(r)}\bigl)>\frac{m}{m-1}A(r)\bigr\}$ and assume that
$F_m\neq\emptyset$. We define inductively a sequence
$$r_0^{\prime\prime}\leq r_1<r_1^{\prime\prime}\leq
r_2<r_2^{\prime\prime}\leq\ldots$$ by
$$r_k=\inf\{r>r_{k-1}^{\prime\prime}:\ r\in F_m\}\quad\text{and}\quad
r_k^{\prime\prime}=r_k+\frac{2r_k}{\phi(r_k)}.$$ Then
$F_m\in\bigcup\limits_{k\geq 1}[r_k,r_k^{\prime\prime}]$. Also we
define
$\rho_k=r_k^{\prime\prime}+\frac{2r_k^{\prime\prime}}{\phi(r_k)}$
and put $H_m=\bigcup\limits_{k\geq 1}[r_k,\rho_k]$. We estimate
the logarithmic measure of $H_m$ as follows:
\begin{eqnarray*}\int\limits_{H_m}\frac{dr}{r}
& \leq & \sum\limits_{k\geq 1}\frac{\rho_k-r_k}{r_k}
=\sum\limits_{k\geq
1}\frac{1}{r_k}\Big(r_k^{\prime\prime}+\frac{2r_k^{\prime\prime}}{\phi(r_k)}-r_k\Big)
\\
& = & \sum\limits_{k\geq
1}\frac{1}{r_k}\Big(r_k+\frac{2r_k}{\phi(r_k)}
+\frac{2r_k+\frac{4r_k}{\phi(r_k)}}{\phi(r_k)}-r_k\Big)
\\
& \leq & \sum\limits_{k\geq
1}\frac{8}{\phi(r_k)}=\sum\limits_{k\geq
1}\frac{8m^2}{(A(r_k))^{\frac{1}{Q-1}}}.
\end{eqnarray*} Since the function $A(r)$
increases, we have $A(r_{k+1})\geq
A(r_k^{\prime\prime})>\frac{m}{m-1}A(r_k)$ and
$$\int\limits_{H_m}\frac{dr}{r}\leq\frac{8m^2}{(A(r_1))^{\frac{1}{Q-1}}}\sum\limits_{k\geq
1}\Big(\frac{m-1}{m}\Big)^{\frac{k}{Q-1}}<\infty.$$ If
$F_m=\emptyset$, then we set $H_m=\emptyset$.

Further, we choose a sequence $d_2<d_3<\ldots$ of numbers such
that
$$d_m\geq 3r_0^{\prime\prime}(m)\quad\text{and}\quad\int\limits_{E_m}\frac{dr}{r}<\frac{1}{2^m},$$
where $E_m=H_m\cap[d_m,\infty)$. If we denote  the union $\bigcup
E_m$ by $E$ then we get $\int\limits_{E}\frac{dr}{r}<\infty$. Let
$\varepsilon_0\in(0,1/5)$ and let $Y=S(w,t)$ be a sphere with
radius $t$. We choose the least integer $m$ satisfying
\begin{itemize}
\item[(i)]{$m\geq 4$,}
\item[(ii)]{$\frac{m^2}{(m-1)^2}<1+\frac{\varepsilon_0}{2}$, and}
\item[(iii)]{$c_1K_I(f)\big(|\ln t|^{Q-1}+c_0\big)\leq
\frac{m}{2^{Q-1}}$.}
\end{itemize}
In this case the inequality $\frac{2}{m^2}<\varepsilon_0$ also
holds. We take $s^{\prime}\geq d_m$ and $s^{\prime}\notin E$. Then
there is $r\geq r_0^{\prime\prime}$ such that
$$s^{\prime}=s+\frac{s}{\varepsilon_0(A(s))^{\frac{1}{Q-1}}}\quad\text{with}\quad
s=r+\frac{r}{\phi(r)}.$$ We claim that $r\notin F_m$. Indeed, if
we suppose that $r\in [r_k,r_k^{\prime\prime}]$ for some $k$, then
we obtain
\begin{eqnarray}\label{282}
r_k<r<s^{\prime} & = & r+\frac{r}{\phi(r)}
+\frac{r+\frac{r}{\phi(r)}}{\varepsilon_0\Big(A(r+\frac{r}{\phi(r)})\Big)^{\frac{1}{Q-1}}}\leq
r+\frac{r}{\phi(r)}+\frac{r+\frac{r}{\phi(r)}}{2\phi(r)}\nonumber
\\
& \leq & r+\frac{2r}{\phi(r)}\leq
r_k^{\prime\prime}+\frac{2r_k^{\prime\prime}}{\phi(r_k)}=\rho_k
\end{eqnarray} from the inequalities $A(r+\frac{r}{\phi(r)})>A(r)$, $\phi(r)>1$, and
$\varepsilon_0>2/m^2$. The estimates $r_k<s^{\prime}<\rho_k$ imply
that $s^{\prime}\in H_m$, but it contradicts to choice of
$s^{\prime}\notin E$.

Now, we apply Lemma~\ref{113} with $\theta =1+\frac{1}{\phi(r)}$.
Using the condition~(iii), the estimates $\ln \theta \geq
\frac{1}{2\phi(r)}$, $s^{\prime}<r+\frac{2r}{\phi(r)}$, and the
definition of $F_m$, we get
\begin{eqnarray*}
\nu(s,Y) & = & \nu(r\theta,Y)\geq A(r)-\frac{c_1K_I(f)}{(\ln
\theta)^{Q-1}}\big(|\ln t|^{Q-1}+c_0\big)\geq
A(r)-m\phi(r)^{Q-1} \\
& = & A(r)-\frac{m}{m^{2Q-2}}A(r)\geq
\frac{m-1}{m}A(r)\geq\frac{(m-1)^2}{m^2}A\Big(r+\frac{2r}{\phi(r)}\Big)
\\
& \geq &
\frac{(m-1)^2}{m^2}A(s^{\prime}).
\end{eqnarray*}
This and the condition~(ii) imply
\begin{equation}\label{220}
\frac{\nu(s,Y)}{A(s^{\prime})}\geq\Big(\frac{m-1}{m}\Big)^{2}>1-\frac{\varepsilon_0}{2}.
\end{equation}

Notice, that since $s^{\prime}\notin E$ and $s^{\prime}\geq d_m$,
we have $s^{\prime}\notin H_m$. Therefore, $s^{\prime}\notin F_m$.
Now, we apply Lemma~\ref{113} with
$\theta=1+\frac{1}{\phi(s^{\prime})}$. Arguing as above we deduce
$$\nu(s^{\prime},Y) \leq
A\Big(s^{\prime}+\frac{s^{\prime}}{\phi(s^{\prime})}\Big) +
m\phi(s^{\prime})^{Q-1}\leq
\Big(\frac{m}{m-1}+\frac{1}{m}\Big)A(s^{\prime})\leq
\Big(\frac{m}{m-1}\Big)^2A(s^{\prime}).$$ Finally,
\begin{equation}\label{221}
\frac{\nu(s^{\prime},Y)}{A(s^{\prime})}\leq\Big(\frac{m}{m-1}\Big)^{2}<1+\frac{\varepsilon_0}{2}.
\end{equation}

The inequalities~(\ref{220}) and~(\ref{221}) imply~(\ref{218}).
Moreover, we have
$$\nu(s,Y)>\frac{2-\varepsilon_0}{2}A(s^{\prime})>\frac{2-\varepsilon_0}{2+\varepsilon_0}\nu(s^{\prime},Y)
\geq (1-\varepsilon_0) \nu(s^{\prime},Y),$$ that
prove~(\ref{219}).

By $m_0\geq 4$, we denote the least integer satisfying the
condition~(iii) and $\frac{2}{m^2}<\varepsilon_0$. Then since
$s^{\prime}\geq d_{m_0}$ we can put $D_{\varepsilon_0}=d_{m_0}$.
In this case
$$1\leq\frac{\big(A(d_{m_0})\big)^{\frac{1}{Q-1}}}{m_0^2}<\varepsilon_0\big(A(D_{\varepsilon_0})\big)^{\frac{1}{Q-1}}$$
that gives the last statement of the theorem.
\end{proof}

Let $E\subset[1,\infty)$ be the exceptional set constructed in the
previous theorem. The points $a_1,\ldots,a_q$ belong to the ball
$B(0,\sigma_M/2)$, $\sigma_M=4\max_{1,\ldots,q}\{1,|a_j|\}$, and
$a_0=\infty$. To apply Theorem~\ref{44} we fix a positive
$\varepsilon_0\leq \min(\frac{1}{5},\frac{1}{8q+9})$.
Then~(\ref{219}) implies that
\begin{equation}\label{224}
\frac{\nu(s^{\prime},\sigma_M)}{\nu(s,\sigma_M)}<1+\frac{\varepsilon_0}{1-\varepsilon_0}<\frac{3}{2}
\end{equation}
and
$\frac{\nu(s^{\prime},1)}{\nu(s,1)}<\frac{1}{1-\varepsilon_0}$.
Moreover, we have $\frac{\nu(s,1)}{A(s^{\prime})}<1+\varepsilon_0$
and
$\frac{A(s^{\prime})}{\nu(s,\sigma_M)}<\frac{1}{1-\varepsilon_0}$
from~\eqref{218}. Three last inequalities yield
\begin{equation}\label{223}
\frac{\nu(s^{\prime},1)}{\nu(s,\sigma_M)}<\frac{1+\varepsilon_0}{(1-\varepsilon_0)^2}
\leq\Big(1-\frac{2\varepsilon_0}{1-\varepsilon_0}\Big)^2 \leq 1+
\frac{8\varepsilon_0}{1-\varepsilon_0}< 1+\frac{1}{q+1}.
\end{equation}
Estimates~\eqref{224} and~\eqref{223} hold whenever $s>0$ is such
that $s^{\prime}\in[\hat\kappa,\infty)\setminus E$ with
$\hat\kappa>\max\{\omega(|\ln 2\sigma_M|),\omega(0)\}$. The
numbers $\varepsilon_0$ will be specified later. Set
$J=\{1,\ldots,q\}$. To prove Theorem~\ref{45} it suffices to show
$$
\sum\limits_{j\in
J\cup\{0\}}\Big(1-\frac{n(s^{\prime},a_j)}{\nu(s^{\prime},1)}\Big)_+\leq
C(Q,K)<\infty.$$ We may assume that
$\frac{n(s^{\prime},a_j)}{\nu(s^{\prime},1)}<1$ for all $j\in
J\cup\{0\}$. Introduce the auxiliary value
$$\Delta_j=1-\frac{n(s^{\prime},a_j)}{\nu(s,\sigma_M)}.$$ Then
\begin{eqnarray}\label{225}
\sum\limits_{j\in
J\cup\{0\}}\Big(1-\frac{n(s^{\prime},a_j)}{\nu(s^{\prime},1)}\Big)
& = & \sum\limits_{j\in J\cup\{0\}}\Delta_j+\sum\limits_{j\in
J\cup\{0\}}\frac{n(s^{\prime},a_j)}{\nu(s^{\prime},1)}
\Big(\frac{\nu(s^{\prime},1)}{\nu(s,\sigma_M)}-1\Big)\nonumber
\\
& \leq & \sum\limits_{j\in
J\cup\{0\}}\Delta_j+\frac{1}{q+1}\sum\limits_{j\in
J\cup\{0\}}\frac{n(s^{\prime},a_j)}{\nu(s^{\prime},1)}
\\ & \leq &
\sum\limits_{j\in J\cup\{0\}}\Delta_j+1.\nonumber
\end{eqnarray}

\subsection{Decomposition of the ball $B(0,s)$}

Here we construct the precise decomposition of the ball $B(0,s)$
into finitely overlapping sets. Let $d=s^{\prime}-s$. We start
from the ball $B(0,s)$ and, inside of it, we construct rings of
increasing diameters. The diameters increase when we move from the
sphere $\partial B(0,s)$ to the center of $B(0,s)$.
$$R_0=\overline B(0,s)\setminus B(0,s-d),
R_1=\overline B(0,s-d)\setminus B(0,s-2d),\ldots $$
$$R_n=\overline B(0,s-2^{n-1}d)\setminus B(0,s-2^{n}d),\ldots.$$
We continue this process up to the step $L$ when
$B(0,s-2^Ld)\subset B(0,2^Ld)$.

Every ring $R_n$ has a diameter $2^{n-1}d$, $n=1,\ldots,L$. Then
we use Wiener lemma (see for instance~\cite[p.~53]{FS}) and cover
the ring $R_n$ by balls
$B(x_i,\frac{2^{n-2}d}{100\varpi\varkappa})$, $x_i\in R_n$, such
that every point of $R_n$ belongs to at most $M$ balls. Here
$\varpi$ is the constant from the generalized triangle inequality
and $\varkappa>6$ is a number specifying later. Also we can
suppose that $\varpi\geq 1$. We cover the ring $R_0$ by balls
$B(x_i,\frac{d}{100\varpi\varkappa})$ and the ball $B(0,2^Ld)$ by
$B(x_i,\frac{2^{L-1}d}{100\varpi\varkappa})$. Notice, that the
balls $B(x_i,\frac{2^{n-2}d}{100\varpi\varkappa})$ intersect only
three rings $R_{n-1}$, $R_n$, and $R_{n+1}$.

Estimate the quantity of balls covering $B(0,s)$. The volume
$V(R_n)$ of one ring $R_n$ is estimated by
$$V(R_n)=(s-2^{n-1}d)^Q-(s-2^{n}d)^Q\leq C_1(Q)s^{Q-1}2^{n-1}d.$$
Since the volume of a ball
$B(x_i,\frac{2^{n-2}d}{100\varpi\varkappa})\cap R_n\ne\emptyset$
is $\Big(\frac{2^{n-2}d}{100\varpi\varkappa}\Big)^Q$, we deduce
the upper bound for the number $p_n$ of balls that have nonempty
intersection with $R_{n+1}$, $R_n$ or $R_{n-1}$:
$$p_n\leq\frac{V(R_{n+1})+V(R_n)+V(R_{n-1})}{\Big(\frac{2^{n-2}d}{100\varpi\varkappa}\Big)^Q}
\leq
\frac{C_2(Q,M,\varpi)\varkappa^Q}{(2^{(Q-1)(n-1)})}\Big(\frac{s}{d}\Big)^{Q-1}.$$
Summing over $n$, we get the estimate for number $p$ of balls
covering $B(0,s)$: $$p  \leq \sum\limits_{n=0}^{L}p_n\leq
C_3(Q,M,\varkappa,\varpi)\Big(\frac{s}{d}\Big)^{Q-1}=
C_3\varepsilon_0^{Q-1}A(s)\leq 2C_3
\varepsilon_0^{Q-1}\nu(s,\sigma_M).$$

Now, we denote by $B(x_i,r_i)$ the balls constructed in the
decomposition of $B(0,s)$. Then we write
$$U_i=B(x_i,r_i)\cap B(0,s),\quad V_i=B(x_i,2r_i),\quad W_i=B(x_i,4r_i),$$
$$X_i=B(x_i,6r_i),\quad Y_i=B(x_i,2r_i\varkappa)\quad Z_i=B(x_i,4r_i\varkappa).
$$ We obtain that $Z_i\in \overline B(0,s^{\prime})$. Let us
estimate the multiplicity of overlapping of balls $Z_i$. It is
sufficient to estimate the multiplicity in one ring, for instance,
in $R_n$. Firstly, we observe that
$Z_k=B(x_k,\frac{2^{n-2}d}{25\varpi})$, $x_k\in R_n$, can
intersect only the balls with the same radius and center from
$R_n$, either the balls $Z_i=B(x_i,\frac{2^{n-1}d}{25\varpi})$,
$x_i\in R_{n+1}$, or $Z_j=B(x_j,\frac{2^{n-3}d}{25\varpi})$,
$x_j\in R_{n-1}$. Indeed, if $y\in
B(x_m,\frac{2^{n}d}{25\varpi})\cap
B(x_k,\frac{2^{n-2}d}{25\varpi})$, $x_m\in R_{n+2}$, $x_k\in
R_{n}$, then
$$
2^{n}d\leq d(x_m,x_k)\leq
\varpi\Big(\frac{2^{n}d}{25\varpi}+\frac{2^{n-2}d}{25\varpi}\Big)\leq
\frac{2^{n-2}d}{5}.
$$ We get a contradiction. If $y\in
B(x_l,\frac{2^{n-4}d}{25\varpi})\cap
B(x_k,\frac{2^{n-2}d}{25\varpi})$, $x_l\in R_{n-2}$, $x_k\in
R_{n}$, then
$$2^{n-2}d\leq d(x_l,x_k)\leq
\varpi\Big(\frac{2^{n-4}d}{25\varpi}+\frac{2^{n-2}d}{25\varpi}\Big)\leq
\frac{2^{n-4}d}{5},$$ which is a contradiction again. We note
also, that if $y\notin B\bigl(x_i,\frac{2^{n-2}d}{5}\bigr)$,
$x_i\in R_{n}$, then the ball
$B\bigl(x_i,\frac{2^{n-2}d}{25\varpi}\bigr)$ does not meet the
balls $B\bigl(y,\frac{2^{n-1}d}{25\varpi}\bigr)$, $y\in R_{n+1}$,
$B\bigl(y,\frac{2^{n-2}d}{25\varpi}\bigr)$, $y\in R_{n}$, and
$B\bigl(y,\frac{2^{n-3}d}{25\varpi}\bigr)$, $y\in R_{n-1}$. It
follows that the multiplicity $\widetilde M$ of overlapping of
$Z_i=B\bigl(x_i\frac{2^{n-2}d}{25\varpi}\bigr)$ cannot exceed
\begin{equation}\label{283}\widetilde M\leq
M\Big(\frac{2^{n-2}d}{5}\Big)^Q\Big/\Big(\frac{2^{n-3}d}{100\varkappa\varpi}\Big)^Q=(40\varkappa\varpi)^Q
M.\end{equation} Now we can choose $\varkappa$. Setting
$$\theta =2\exp\left(\Big(\frac{\kappa(\mathbb
G,Q)K_O(f)K_I(f)}{c(Q)\ln 6/5 }\Big)^{\frac{1}{Q-1}}\right)$$ we
put $\varkappa=3\theta$. Here the constant $c(Q)$ is from the
proof of Corollary~\ref{33}.

\begin{remark}\label{r4} Notice, that by the choice of $\theta$ we
have
$$\ln\theta >\Big(\frac{\kappa(\mathbb G,Q)K_O(f)K_I(f)}{c(Q)\ln 6/5
}\Big)^{\frac{1}{Q-1}}$$ and
$$\ln\frac{3}{2}-\frac{\kappa(\mathbb G,Q)K_O(f)K_I(f)}{c(Q)(\ln \theta)^{Q-1}}>\ln\frac{5}{4}.$$
\end{remark}

\subsection{Estimates for partial liftings}\label{71}

We remember that the points $a_1,\ldots,a_q$ belong to
$B(0,\sigma_M/2)$ and $a_0=\infty$. Now we describe a rectifiable
curve $\gamma_y^{j}(t)$ that connect the points $a_j$,
$j=0,\ldots,q$ with $z\in S(0,\sigma_M)\setminus\mathcal Z$.

Let $G_y(t):[0,\sigma_M]\to[0,t_y]$ be an affine mapping with
$G_y(0)=t_y$, $G_y(\sigma_M)=0$, and such that
$a_j\varphi_{G_y(\sigma_M)}(y)=a_j\varphi_{0}(y)=a_j$,
$a_j\varphi_{G_y(0)}(y)=a_j\varphi_{t_y}(y)=z$ where $z\in
S(0,\sigma_M)$, $j=1,\ldots,q$, and $\varphi_s(y)$ is a radial
curve for $y\in S\setminus\mathcal Z$. We put
$\gamma_y^{j}(t)=a_j\varphi_{G_y(t)}(y)$, $j=1,\ldots,q$. Then the
rectifiable curve $\gamma_y^{j}(t):[0,\sigma_M]\to \overline
B(0,\sigma_M)$ connect $a_j$ with points $z\in
S(0,\sigma_M)\setminus\mathcal Z$ such that
$a_j=\gamma_y^{j}(\sigma_M)$, $z=\gamma_y^{j}(0)$. In view of
$t_y\in[\sigma_M/2,2\sigma_M]$ and~\eqref{280}, we deduce
\begin{equation}\label{281}|\dot\gamma_y^{j}(t)|_0=|\dot\varphi_{G_y(t)}(y)|_0=\upsilon^{-1}(y)|\dot
G_y(t)|\leq 2\upsilon^{-1}(y).\end{equation} Also, we consider
locally rectifiable curves $\gamma_y^{0}(t):[\sigma_M,\infty)\to
\mathbb C B(0,\sigma_M)$, $\gamma_y^{0}(t)=\varphi_t(y)$, $y\in
S\setminus\mathcal Z$. These curves joint the point
$a_0=\{\infty\}$ with $z\in S(0,\sigma_M)\setminus\mathcal Z$,
$a_0=\lim\limits_{t\to\infty}\gamma_y^{0}(t)$,
$z=\gamma_y^{0}(\sigma_M)$.

>From now on, up to the end of the paper, we suppose that the
points $y\in S\setminus\mathcal Z$ correspond to the points $z\in
S(0,\sigma_M)\setminus\mathcal Z$ as it was described in the
preceding paragraph.

Let $f_0=f\vert_{B(0,2s^{\prime})}$. We fix $i\in
I=\{1,\ldots,p\}$ and $j\in J\cup\{0\}=\{0,1,\ldots,q\}$. Put
$\{x_1,\ldots,x_k\}=f^{-1}(z)\cap U_i$, $z\in
S(0,\sigma_M)\setminus\mathcal Z$, $z=\gamma_y^{j}(0)$ for $j\in
J$ or $z=\gamma_y^{0}(\sigma_M)$. Let us consider the maximal
sequence of essentially separate $f_0$-liftings
$\beta_1,\ldots,\beta_m$, $m=n(U_i,z)$, of a curve $\gamma_y^j$
starting at $\{x_1,\ldots,x_k\}$. In this section we will work
only with those curves of $\beta_{\mu}$ whose locus is not
contained in $\overline B(0, s^{\prime})$. We denote them by
$\alpha_1,\ldots,\alpha_{\mu_y}$ for each $y\in S\setminus\mathcal
Z$. We correlate the parametrization of $\gamma_y^j$, $j\in J$,
and $\alpha_{\mu}$ as in Remark~\ref{r3}. Then $\alpha_{\mu}:\
[0,\sigma_M]\to\mathbb G$. The curve
$\gamma^0_y(t)\vert_{[\sigma_M,R]}$, $\sigma_M<R<\infty$, is
rectifiable. Making use of correlation of Remark~\ref{r3}, we may
assume that $f_0$-liftings of $\gamma^0_y(t)\vert_{[\sigma_M,R]}$
are curves $\alpha_{\mu}:\ [\sigma_M,R]\to\mathbb G$,
$\sigma_M<R<\infty$. Taking the supremum over all closed parts of
$\gamma_y^0$ we get that $f_0$-liftings of $\gamma_y^0$ are curves
$\alpha_{\mu}:[\sigma_M,\infty)\to\G$. We fix values of parameters
$0\leq u_{\mu,y}<v_{\mu,y}<w_{\mu,y}<\sigma_M$ if $j\in J$ and
$\sigma_M\leq u_{\mu,y}<v_{\mu,y}<w_{\mu,y}<\infty$ if $j=0$ for
each $\alpha_{\mu}$ such that
\begin{equation}\label{258}\alpha_{\mu}(u_{\mu,y})\in\partial U_i,\quad
\alpha_{\mu}(v_{\mu,y})\in\partial V_i,\quad
\alpha_{\mu}(w_{\mu,y})\in\partial W_i.\end{equation} We use the
notation
$$\alpha_{\mu}^{(1)}=\alpha_{\mu}\vert_{[u_{\mu,y},v_{\mu,y}]},\quad
\alpha_{\mu}^{(2)}=\alpha_{\mu}\vert_{[v_{\mu,y},w_{\mu,y}]}.
$$
We want to estimate the number of $f_0$-liftings of different
parts of $\gamma_y^j$, moving from $z=\gamma_y^j(0)$ to $a_j$ in
the case when $j\in J$ or advancing from $z=\gamma_y^0(\sigma_M)$
to $\infty$ for $j=0$.

Set $\varsigma_0=\frac{1}{16}\min\limits_{1\leq j\neq k\leq
q}|a_j^{-1}a_k|$. We introduce
$$L^j_i(y)=\{\mu\in\{1,\ldots,\mu_y\}:\
\frac{\sigma_M-u_{\mu,y}}{\sigma_M-v_{\mu,y}}\leq \frac{\sigma_M}{\varsigma_0}\}
\qquad\text{for}\qquad j=1,\ldots,q$$ and
$$L^0_i(y)=\Big\{\mu\in\{1,\ldots,\mu_y\}:\
\frac{v_{\mu,y}}{u_{\mu,y}}\leq
\frac{\sigma_M}{\varsigma_0}\Big\}.$$ Since we can not say
anything about measurability of $\card L^j_i(y)$ we need to
estimate $\card L_i^j(y)$ by some measurable function.
\begin{lemma}\label{118}
There exists a nonnegative measurable function $l_i^j(y):\
S\setminus\mathcal Z\to\mathbb R^1$ such that
\begin{equation}\label{241}
\card L^j_i(y)\leq l_i^j(y)\qquad \text{for any}\qquad y\in
S\setminus\mathcal Z,\quad  j\in J\cup\{0\},
\end{equation}
\begin{equation}\label{242}
\int\limits_{S\setminus\mathcal Z}l_i^j(y)\,d\sigma(y) \leq
c_5K_I(f)\Big(\ln \frac{\sigma_M}{\varsigma_0}\Big)^{Q-1}.
\end{equation}
\end{lemma}

\begin{proof}
The properties of a quasimeromorphic mapping imply that there is a
Borel set $F\subset\mathbb G$ with $\mes(F)=0$ containing the set
of points, where $f$ is not $\mathcal P$-differentiable. Then
$\mes(f(F))=0$. In this case the one dimensional Hausdorff measure
of intersection $\gamma^j_y(t)\cap f(F)$ vanishes for almost all
points $z\in S(0,\sigma_M)\setminus\mathcal Z$ with respect to the
$\sigma$-measure. Indeed, if we assume that the $\sigma$-measure
of supposed exceptional set $E_1\in S(0,\sigma_M)\setminus\mathcal
Z$ is positive, then, making use of the Fubini theorem and that
$\sigma$ is absolutely continuous with respect to
$\sigma^{\star}$, we get a contradiction with $\mes(f(F))=0$. Let
$\Gamma_0$ be a family of locally rectifiable curves
$\gamma^j_y(t)\setminus a_j$, $y\in S\setminus\mathcal Z$, that
have closed parts with the $f_0$-liftings $\alpha_{\mu,y}$ of
which are not absolutely continuous. Lemma~\ref{Pollem} states
that the $Q$-module of $\Gamma_0$ vanishes. We conclude that for
almost all points $z\in S(0,\sigma_M)\setminus\mathcal Z$ curves
$\alpha_{\mu,y}$, $\mu=1,\ldots,\mu_y$, are absolutely continuous
for all $j=0,1,\ldots,q$ by Lemma~\ref{126}. We denote by $E_2$
this exceptional set of $S(0,\sigma_M)$. Put $E\subset
S(0,\sigma_M)$ be a Borel set such that $E_1\cup E_2\subset E$ and
$\sigma(E)=0$.

We start from the case when $j=1,\ldots,q$. Let $\rho$ be an
admissible function for the family of curves connecting $U_i$ with
$\partial V_i$ in $V_i$ such that $\rho\vert_{\mathbb
G\setminus(V_i\setminus \overline{U_i})}=0$ and
$\int\limits_{V_i}\rho^Q\,dx\leq 2M_Q(\overline{U_i},V_i)$. We fix
$z\in S(0,\sigma_M)\setminus (E\cup\mathcal Z)$. For the
corresponding $y\in S$ and $\mu=1,\ldots,\mu_y$ we define
$\rho^{\star}_{\mu}$ on $f(\alpha_{\mu})$ by
\begin{equation}\label{243}\rho^{\star}_{\mu}(z)=\rho^{\star}_{\mu}(f(x))=\left\{\array{ll}
\rho(x)\lambda_f(x),\quad  & x\in\alpha_{\mu}\cap\mathbb C F,
\\
0 & \text{otherwise},
\endarray\right.\end{equation} where $1/\lambda_f(x)=\min\limits_{|\xi|_0=1,\xi\in
V_1}|D_Hf(x)\xi|_0$. Then by~\eqref{281}
\begin{eqnarray}\label{274}
2\int\limits_{u_{\mu,y}}^{v_{\mu,y}}\rho^{\star}_{\mu}(\gamma_y^j(t))\upsilon(y)^{-1}\,dt
& \geq &
\int\limits_{u_{\mu,y}}^{v_{\mu,y}}\rho^{\star}_{\mu}(\gamma_y^j(t))|\dot\gamma_y^j(t)|_0\,dt\nonumber
\\ & = &
\int\limits_{u_{\mu,y}}^{v_{\mu,y}}\rho(\alpha_{\mu}^{(1)}(t))\lambda_f(\alpha_{\mu}^{(1)})
|\dot f(\alpha_{\mu}^{(1)}(t))|_0\,dt\nonumber
\\ & \geq & \int\limits_{u_{\mu,y}}^{v_{\mu,y}} \rho(\alpha_{\mu}^{(1)}(t))|\dot\alpha^{(1)}(t)|_0\,dt
=\int\limits_{\alpha^{(1)}_{\mu}} \rho\,dt\geq 1\end{eqnarray} for
each point $y\in S$ such that the corresponding point $z$ belongs
to $S(0,\sigma_M)\setminus(E\cup\mathcal Z)$. Applying the
H\"{o}lder  inequality we obtain
$$1\leq\Big(2^Q\int\limits_{u_{\mu,y}}^{v_{\mu,y}}\big(\rho^{\star}_{\mu}(\gamma_y^j(t))\big)^Q
\upsilon(y)^{-Q}(\sigma_M-t)^{Q-1}\,dt\Big)
\Big(\int\limits_{u_{\mu,y}}^{v_{\mu,y}}\frac{dt}{\sigma_M-t}\Big)^{Q-1}.$$
We get the next estimate for the last integral
$$\Big(\ln\frac{\sigma_M-u_{\mu,y}}{\sigma_M-v_{\mu,y}}\Big)^{1-Q}
\leq
2^Q\int\limits_{u_{\mu,y}}^{v_{\mu,y}}\big(\rho^{\star}_{\mu}(\gamma_y^j(t))\big)^Q\upsilon(y)^{-Q}(\sigma_M-t)^{Q-1}\,dt.$$
Summing over $L_i^j(y)$ yields
\begin{eqnarray*}
\card L_i^j(y)\Big(\ln \frac{\sigma_M}{\varsigma_0}\Big)^{1-Q} &
\leq &
\sum\limits_{\mu=1}^{\mu_y}\Big(\ln\frac{\sigma_M-u_{\mu,y}}{\sigma_M-v_{\mu,y}}\Big)^{1-Q}
\\
& \leq &
2^Q\sum\limits_{\mu=1}^{\mu_y}\int\limits_{u_{\mu,y}}^{v_{\mu,y}}
\big(\rho^{\star}_{\mu}(\gamma_y^j(t))\big)^Q\upsilon(y)^{-Q}(\sigma_M-t)^{Q-1}\,dt
\\
& \leq & 2^Q\int\limits_{0}^{\sigma_M}\sum\limits_{x\in
f^{-1}(\gamma_y^j(t))}(\rho(x)\lambda_f(x))^Q\upsilon^{-Q}(\sigma_M-t)^{Q-1}\,dt.
\end{eqnarray*}
Now, we define the measurable function $l_i^j(y):\
S\setminus\mathcal Z\to\mathbb R$ satisfying~(\ref{241}), by
$$l_i^j(y)=2^Q
\Big(\ln
\frac{\sigma_M}{\varsigma_0}\Big)^{Q-1}\int\limits_{0}^{\sigma_M}\sum\limits_{x\in
f^{-1}(\gamma_y^j(t))}(\rho(x)\lambda_f(x))^Q\upsilon(y)^{-Q}(\sigma_M-t)^{Q-1}\,dt
$$ if $z\in S(0,\sigma_M)\setminus (E\cup\mathcal Z)$, and $l_i^j(y)=n(U_i,z)$ if $z\in E$.
Let us show that $l_i^j(y)$ satisfies~(\ref{242}). We denote by
$\phi(w)$ an auxiliary map from $\overline B(0,\sigma_M)$ to
itself such that
$\phi(w)=\phi(\varphi_{(\sigma_M-t)}(y))=\gamma^j_y(t)$ for
$t\in[0,\sigma_M]$, $y\in S\setminus\mathcal Z$. Then integrating
over $S\setminus\mathcal Z$ with respect to the $\sigma$-measure
we deduce
\begin{eqnarray*}
\int\limits_{S\setminus\mathcal Z}l_i^j(y)\,d\sigma(y) & = &
2^Q\Big(\ln
\frac{\sigma_M}{\varsigma_0}\Big)^{Q-1}\int\limits_{S\setminus\mathcal
Z}d\sigma^{\star}(y)\int\limits_{0}^{\sigma_M}\sum\limits_{x\in
f^{-1}(\phi(w))}(\rho(x)\lambda_f(x))^Qt^{Q-1}\,dt
\\
& \leq & c_4\Big(\ln
\frac{\sigma_M}{\varsigma_0}\Big)^{Q-1}\int\limits_{B(0,\sigma_M)}\sum\limits_{x\in
f^{-1}(z)}(\rho(x)\lambda_f(x))^Q\,dz.
\end{eqnarray*}
For any point $z\in B(0,\sigma_M)\setminus f(B_f)$, there exists
some neighborhood $W$ where the inverse mapping $f^{-1}$ exists
and homeomorphic. We write $V_1,\ldots, V_k$ for the components of
$f^{-1}(W)\cap \overline B(s^{\prime})$. The mappings
$f_j=f\vert_{V_j}:\ V_j\to W$ are quasiconformal. Then
\begin{eqnarray*}\int\limits_{W}\sum\limits_{x\in
f^{-1}(z)}(\rho(x)\lambda_f(x))^Q\,dz & = &
\sum\limits_{j=1}^{k}\int\limits_{V_j}(\rho(x)\lambda_f(x))^QJ(x,f)\,dx
\\
& \leq &
K_I(f)\int\limits_{f^{-1}(W)}\rho(x)^Q\,dx.\end{eqnarray*} The set
$B(0,\sigma_M)\setminus f(B_f)$ can be covered up to the a set of
measure zero by disjoint neighborhoods of this kind.  It follows
$$\int\limits_{S\setminus\mathcal Z}l_i^j(y)\,d\sigma(y)\leq c_4 K_I(f)\Big(\ln
\frac{\sigma_M}{\varsigma_0}\Big)^{Q-1}\int\limits_{\mathbb
G}\rho(x)^Q\,dx.$$ We have~(\ref{242}) with $c_5^{\prime}=c_4
2M_Q(\overline{U_i},V_i)$.

Now we consider a locally rectifiable curve
$\gamma_y^{0}(t):[\sigma_M,\infty)\to \mathbb C B(0,\sigma_M)$
joining $z\in S(0,\sigma_M)\setminus\mathcal Z$ with $\infty$,
such that $\gamma_y^{0}(\sigma_M)=z$ and
$\lim\limits_{t\to\infty}\gamma_y^{0}(t)=\infty$. We will define a
measurable function $l^0_i(y):\ S\setminus\mathcal Z\to\mathbb R$
satisfying~(\ref{241}) and~(\ref{242}). As at the beginning of the
proof, for the functions $\rho$ and $\rho^{\star}$, we
obtain~\eqref{274} for any $y\in S\setminus\mathcal Z$, $j=0$. The
H\"{o}lder inequality implies $$1\leq
\Big(\int\limits_{u_{\mu,y}}^{v_{\mu,y}}\big(\rho^{\star}_{\mu}(\gamma_y^0(t))\big)^Q\upsilon(y)^{-Q}t^{Q-1}\,dt\Big)
\Big(\int\limits_{u_{\mu,y}}^{v_{\mu,y}}\frac{dt}{t}\Big)^{Q-1}.$$
Making use of the estimate
$\Big(\ln\frac{v_{\mu,y}}{u_{\mu,y}}\Big)^{1-Q} \leq
\int\limits_{u_{\mu,y}}^{v_{\mu,y}}\big(\rho^{\star}_{\mu}(\gamma_y^0(t))\big)^Q\upsilon(y)^{-Q}t^{Q-1}\,dt$
and summing over $\mu\in L^0_i(y)$, we get
\begin{eqnarray*}
\frac{\card L_i^0(y)}{(\ln (\sigma_M/\varsigma_0))^{Q-1}} & \leq &
\sum\limits_{\mu=1}^{\mu_y}\Big(\ln\frac{v_{\mu,y}}{u_{\mu,y}}\Big)^{1-Q}
\leq
\sum\limits_{\mu=1}^{\mu_y}\int\limits_{u_{\mu,y}}^{v_{\mu,y}}
\big(\rho^{\star}_{\mu}(\gamma_y^0(t))\big)^Q\upsilon(y)^{-Q}t^{Q-1}\,dt
\\
& \leq & \int\limits_{\sigma_M}^{\infty}\sum\limits_{x\in
f^{-1}(\gamma_y^0(t))}(\rho(x)\lambda_f(x))^Q\upsilon(y)^{-Q}t^{Q-1}\,dt.
\end{eqnarray*}
We define the function $l^0_i(y)$ by
$$l_i^0(y)=\Big(\ln
\frac{\sigma_M}{\varsigma_0}\Big)^{Q-1}\int\limits_{\sigma_M}^{\infty}\sum\limits_{x\in
f^{-1}(\gamma_y^0(t))}(\rho(x)\lambda_f(x))^Q\upsilon(y)^{-Q}t^{Q-1}\,dt$$
if $z\in S(0,\sigma_M)\setminus(\mathcal Z\cup E)$, and
$l_i^0(y)=n(U_i,z)$ if $z\in E$. Obviously, that $l^0_i(y)$
satisfies~(\ref{241}). To show~(\ref{242}) we argue as in the
previous case, integrate over $S\setminus\mathcal Z$, and  obtain
$$\int\limits_{S\setminus\mathcal Z}l_i^0(y)\,d\sigma(y) = \Big(\ln
\frac{\sigma_M}{\varsigma_0}\Big)^{Q-1}\int\limits_{\mathbb C
B(0,\sigma_M)}\sum\limits_{x\in
f^{-1}(z)}(\rho(x)\lambda_f(x))^Q\,dz.$$ Covering $\mathbb
CB(0,\sigma_M)$ by disjoint neighborhoods, where the mapping
$f^{-1}$ is homeomorphic, we deduce
$$\int\limits_{S\setminus\mathcal Z}l_i^0(y)\,d\sigma(y)\leq K_I(f)\Big(\ln
\frac{\sigma_M}{\varsigma_0}\Big)^{Q-1}\int\limits_{\mathbb
G}\rho^Q\,dx\leq c_5^{\prime\prime}K_I(f)\Big(\ln
\frac{\sigma_M}{\varsigma_0}\Big)^{Q-1}.$$ We have~(\ref{242})
with $c_5^{\prime\prime}=2M_Q(\overline{U_i},V_i)$. Setting
$c_5=\max\{c^{\prime}_5,c^{\prime\prime}_5\}$ we end the proof.
\end{proof}

Now we estimate the cardinality of $\alpha_{\mu}^{(2)}$. For each
$i\in I$, we let $\sigma_i\in]0,\varsigma_0]$ be a number defined
by
\begin{equation}\label{250}\Big(\ln\frac{\varsigma_0}{\sigma_i}\Big)^{Q-1}
=A_i\nu(U_i,\sigma_M),\end{equation} where $A_i$ will be given
later. Set
$$M^j_i(y)=\Big\{\mu\in\{1,\ldots,\mu_y\}:\frac{\sigma_M-v_{\mu,y}}{\sigma_M-w_{\mu,y}}
\leq \frac{\varsigma_0}{\sigma_i}\Big\} \qquad\text{for}\qquad
j\in J$$ and
$$M^0_i(y)=\Big\{\mu\in\{1,\ldots,\mu_y\}:\ \frac{w_{\mu,y}}{v_{\mu,y}}\leq
\frac{\varsigma_0}{\sigma_i}\Big\},$$ where index $0$ corresponds
to the curve $\gamma^0_i$ joining $\infty$ with $z\in
S(0,\sigma_M)\setminus\mathcal Z$.
\begin{lemma}\label{119}
There exists a nonnegative measurable function $m_i^j(y):\
S\setminus\mathcal Z\to\mathbb R$ such that
\begin{equation}\label{244}
\card \Big(M^j_i(y)\setminus L^j_i(y)\Big)  \leq  m_i^j(y)\quad
\text{for any $y\in S\setminus\mathcal Z$, $j\in
J\cup\{0\}$,}\quad \text{and}
\end{equation}
\begin{equation}\label{245}
\sum\limits_{j\in J\cup\{0\}}\int\limits_{S\setminus\mathcal
Z}m_i^j(y)\,d\sigma(y) \leq c_6 K_I\Big(\ln
\frac{\varsigma_0}{\sigma_i}\Big)^{Q-1}.
\end{equation}
\end{lemma}

\begin{proof}
Let $E\subset S(0,\sigma_M)$ be the exceptional set that we
defined in the previous lemma. Let $\rho$ be an admissible
function for $\Gamma(\overline V_i,W_i)$ such that
$\rho\vert_{\mathbb G\setminus(W_i\setminus \overline{V_i})}=0$.

Fix $z\in S(0,\sigma_M)\setminus(\mathcal Z\cup E)$ and put the
function $\rho^{\star}_{\mu}$ as in~(\ref{243}). If $\mu\in
M^j_i(y)\setminus L^j_i(y)$ for $j\in J$ then
$\sigma_M-v_{\mu,y}<\varsigma_0$. We estimate
\begin{equation*}
\card \big(M^j_i(y)\setminus L^j_i(y)\big)  \leq 2^Q
\Big(\ln\frac{\varsigma_0}{\sigma_i}\Big)^{Q-1}
\sum\limits_{\mu\in M^j_i\setminus
L^j_i}\int\limits_{v_{\mu,y}}^{w_{\mu,y}}
\big(\rho^{\star}_{\mu}(\gamma_y^j(t))\big)^Q\upsilon^{-Q}(\sigma_M-t)^{Q-1}\,dt
\end{equation*}
\begin{equation*}\leq 2^Q
\Big(\ln\frac{\varsigma_0}{\sigma_i}\Big)^{Q-1}
\int\limits_{\sigma_M-\varsigma_0}^{\sigma_M}\sum\limits_{x\in
f^{-1}(\gamma_y^j(t))}
\big(\rho(x)\lambda_f(x))\big)^Q\upsilon(y)^{-Q}(\sigma_M-t)^{Q-1}\,dt.
\end{equation*}
Define the measurable function $m^j_i(y):\ S\setminus\mathcal
Z\to\mathbb R$ by
$$m_i^j(y)=2^Q\Big(\ln
\frac{\varsigma_0}{\sigma_i}\Big)^{Q-1}\int\limits_{\sigma_M-\varsigma_0}^{\sigma_M}\sum\limits_{x\in
f^{-1}(\gamma_y^j(t))}(\rho(x)\lambda_f(x))^Q\upsilon(y)^{-Q}(\sigma_M-t)^{Q-1}\,dt$$
if $z\in S(0,\sigma_M)\setminus(\mathcal Z\cup E)$ and
$m_i^j(y)=n(U_i,z)$ if $z\in E$. The property~(\ref{244}) is
evident. To show~(\ref{245}) we integrate over $S\setminus\mathcal
Z$ and obtain
$$\int\limits_{S\setminus\mathcal Z}m_i^j(y)\,d\sigma(y)\leq \hat c_4\Big(\ln
\frac{\varsigma_0}{\sigma_i}\Big)^{Q-1}\int\limits_{B(a_j,3\varsigma_0/2)}\sum\limits_{x\in
f^{-1}(w)}(\rho(x)\lambda_f(x))^Q\,dw.$$ Here we used an auxiliary
map $\phi(w)=\phi(\varphi_{(\sigma_M-t)}(y))=\gamma_y^j(t)$ as in
the previous lemma. Since the balls $B(a_j,3\varsigma_0/2)$ are
disjoint, summing over $j\in J$, we deduce
\begin{equation}\label{246} \sum\limits_{j\in
J}\int\limits_{S\setminus\mathcal Z}m_i^j(y)\,d\sigma(y) \leq \hat
c_4\Big(\ln
\frac{\varsigma_0}{\sigma_i}\Big)^{Q-1}\int\limits_{B(0,\sigma_M)}\sum\limits_{x\in
f^{-1}(w)}(\rho(x)\lambda_f(x))^Q\,dw. \end{equation} As in the
proof of Lemma~\ref{118} the estimate~\eqref{246}
gives~(\ref{245}) with the constant $c_6^{\prime}=\hat
c_42M_Q(\overline{V_i},W_i)$.

Now, we show~(\ref{244}) and~(\ref{245}) for $j=0$. Arguing as in
Lemma~\ref{118} we obtain
\begin{eqnarray*}
& \card\big(M^0_i(y)\setminus L^0_i(y)\big)  \leq
\Big(\ln\frac{\varsigma_0}{\sigma_i}\Big)^{Q-1}\sum\limits_{\mu\in
M^0_i\setminus L^0_i}\int\limits_{v_{\mu,y}}^{w_{\mu,y}}
\big(\rho^{\star}_{\mu}(\gamma_y^0(t))\big)^Q\upsilon(y)^{-Q}t^{Q-1}\,dt
\\
& \leq  \Big(\ln\frac{\varsigma_0}{\sigma_i}\Big)^{Q-1}
\int\limits_{\sigma_M}^{\infty}\sum\limits_{x\in
f^{-1}(\gamma_y^0(t))}
\big(\rho(x)\lambda_f(x))\big)^Q\upsilon(y)^{-Q}t^{Q-1}\,dt.
\end{eqnarray*}

We define $m^0_i(y)$ by
$$m_i^0(y)=\Big(\ln
\frac{\varsigma_0}{\sigma_i}\Big)^{Q-1}\int\limits_{\sigma_M}^{\infty}\sum\limits_{x\in
f^{-1}(\gamma_y^0(t))}(\rho(x)\lambda_f(x))^Q\upsilon(y)^{-Q}t^{Q-1}\,dt$$
if  $z\in S(0,\sigma_M)\setminus(\mathcal Z\cup E)$ and
$m_i^0(y)=n(U_i,z)$ if $z\in E$. The function $m^0_i(y):\
S\setminus\mathcal Z\to\mathbb R$ is measurable and
satisfies~(\ref{244}). Integrating over $S\setminus\mathcal Z$ and
arguing as in Lemma~\ref{118}, we have
\begin{eqnarray}\label{247}
\int\limits_{S\setminus\mathcal Z}m_i^0(y)\,d\sigma(y) & = &
\Big(\ln
\frac{\varsigma_0}{\sigma_i}\Big)^{Q-1}\int\limits_{\mathbb C
B(0,\sigma_M)}\sum\limits_{x\in
f^{-1}(w)}(\rho(x)\lambda_f(x))^Q\,dw\nonumber
\\
&\leq & K_I(f)\Big(\ln
\frac{\varsigma_0}{\sigma_i}\Big)^{Q-1}\int\limits_{\mathbb
G}\rho(x)^Q\,dx
\\ & \leq & K_I(f)2M_Q(\overline V_i,W_i)\Big(\ln \frac{\varsigma_0}{\sigma_i}\Big)^{Q-1}\nonumber
\\ & = & c_6^{\prime\prime}K_I(f)\Big(\ln \frac{\varsigma_0}{\sigma_i}\Big)^{Q-1}.\nonumber
\end{eqnarray}
We obtain~(\ref{245}) from~(\ref{246}) and~(\ref{247}) with the
constant $c_6=\max\{c_6^{\prime},c_6^{\prime\prime}\}$.
\end{proof}

\subsection{Estimates for extremal maximal sequences of liftings}

In the previous section we estimated the number of different parts
of $f_0$-liftings of the curve $\gamma_y^j$. But we do not know
precisely how to chose these parts. In the present section we give
a rule how to choose the $f_0$-liftings.

Fix $j\in J\cup \{0\}$ for a moment and consider the curves
$\gamma_y^j$. For each point $y\in S\setminus\mathcal Z$ we
consider maximal essentially separate sequences
$\Lambda=(\lambda_1,\ldots,\lambda_g)$ of $f_0$-liftings of
$\gamma_y^j$ starting at the points
$(x_1,\ldots,x_l)=f^{-1}(z)\cap\overline B(0,s)$. Here
$z=\gamma^{j}_y(0)$ for $j\in J$ or $z=\gamma^{0}_y(\sigma_M)$. In
this case $g=n(s,z)=\sum\limits_{k=1}^{l}i(x_k,f)$. Let $\Omega_y$
denote the set of all these sequences. We introduce measurable
functions $\psi_i^j(y)=\psi_i(y,a_j)$, $\psi_i^j(y):\
S\setminus\mathcal Z\to \mathbb R$, $i=1,\ldots,p$, that help us
to calculate the number of $f_0$-liftings of $\gamma_y^j$,
starting from $U_i$, locus of which belongs to $\overline
B(0,s^{\prime})$. We recall that $p$ is the number of sets $U_i$
in the decomposition of $B(0,s)$.

For $\Lambda\in\Omega_y$ we define $$N(\Lambda)=\card\{\nu:\
\lambda_{\nu}\subset \overline
B(0,s^{\prime})\},\qquad\psi^j(y)=\sup\{N(\Lambda):\
\Lambda\in\Omega_y\}.$$ By $\Omega(y)$, we denote the set of
sequences $\Lambda\in \Omega_y$ for which $N(\Lambda)=\psi^j(y)$.
Set
$$U=B(0,s)\setminus \Big(\bigcup\limits_{i=1}^{p}\partial U_i\Big).$$

\begin{lemma}\label{121}
Under the above-mentioned notations, the function $\psi^j(y)$ is
upper semicontinuous at points $y\in S\setminus(\mathcal Z\cup
f(\partial U))$.
\end{lemma}

\begin{proof}
Let $y_0\in S\setminus(\mathcal Z\cup f(\partial U))$ and
$\{y_h\}_{h=1}^{\infty}$ be some sequence of points in $S\setminus
(\mathcal Z\cup f(\partial U))$ that converges to $y_0$. We want
to show
\begin{equation}\label{248}
\limsup\limits_{h\to\infty}\psi^j(y_h)\leq
\psi^j(y_0).\end{equation}

Firstly, we prove that if $y_h\to y_0$, then the maximal sequences
of $f_0$-liftings of $\gamma^j_{y_h}$ converge to some sequence of
$f_0$-liftings of $\gamma^j_{y_0}$. Then we show that the obtained
sequence is essentially separate and deduce the
inequality~(\ref{248}).

Passing to a subsequence of $\{y_h\}$(that we denote by the same
symbol $\{y_h\}$) we can assume that, for some integer $m$ such
that $\psi^j(y_h)=m$ and for all $h\geq 1$, the following
properties hold.
\begin{itemize}
\item[(i)]{There exist normal neighborhoods $V_1,\ldots,
V_l\subset U$ of points $(x_1,\ldots,x_l)=f^{-1}(z_0)\cap U$ such
that $z_h\in f(V_1)\cap\ldots\cap f(V_l)$ for all $h\geq 1$. Here
$z_0=\gamma^j_{y_0}(0)$, $z_h=\gamma^j_{y_h}(0)$ for $j\in J$ or
$z_0=\gamma^0_{y_0}(\sigma_M)$, $z_h=\gamma^0_{y_h}(\sigma_M)$.}
\item[(ii)]{There is a maximal essentially separate sequence
$\Lambda_h=(\lambda_{h,1},\ldots,\lambda_{h,g})\in\Omega(y_h)$,
$g=n(s,z_0)$, such that each $\lambda_{h,\nu}$ starts at a point
in $f^{-1}(z_h)\cap V_{\mu}$ for some $\mu$ that depends on $\nu$
and is independent of $h$. Moreover, the locus of each curve
$\lambda_{h,\nu}$ belongs to $\overline B(0,s^{\prime})$ for
$1\leq\nu\leq m$.}\item[(iii)]{The bound $\upsilon(y_h)^{-1}\leq
K$ holds for all $h=0,1,\ldots$.}
\end{itemize} This choice is possible because of properties of the
local index and the topological degree. Next, we consider two
cases: {\bf 1.} $1\leq \nu\leq m$ and {\bf 2.} $m+1\leq\nu\leq g$.

{\it Case 1.} We claim that the family
$\{\lambda_{h,\nu}\}_{h=1}^{\infty}$ is equicontinuous for
any~$\nu$. Here we need to separate the consideration for $j\in J$
and $j=0$. We start from the indexes $j\in J$. Since
$\lambda_{h,\nu}\subset\overline B(0,s^{\prime})$ for
$\nu=1,\ldots,m$, each curve $\lambda_{h,\nu}$ is defined on
$[0,\sigma_M]$. Fix $\varepsilon>0$. For each $t\in[0,\sigma_M]$
there is $\delta(t)>0$ such that for all $\rho\in(0,\delta(t)]$,
$U(\xi,\rho)=f^{-1}(B(f(\xi),\rho))$ is a normal neighborhood of
$\xi\in f^{-1}(\gamma^j_{y_0}(t))\cap\overline B(0,s^{\prime})$
with $\diam(U(\xi,\rho))<\varepsilon$. We cover the curve
$\gamma^j_{y_0}$ by finite number of balls
$B(\gamma^j_{y_0}(t),\delta(t)/2\varpi)$, that we denote by
$B(\eta_u,\rho_u)$, $u=1,\ldots,v$. Here $\varpi$ is the constant
from the generalized triangle inequality. We may suppose that the
curves $\gamma^j_{y_h}$, $h\geq 1$, belong to the tube
$\bigcup_{u=1}^{v}B(\eta_u,\rho_u)$. For any $t\in[0,\sigma_M]$,
there is $u$ such that $\gamma^j_{y_h}(t)\in B(\eta_u,\rho_u)$. We
recall that the distances $d(\cdot,\cdot)$ and $d_c(\cdot,\cdot)$
are equivalent with a constant $\tilde c>0$: $\tilde
c^{-1}d(\cdot,\cdot)\leq d_c(\cdot,\cdot)\leq \tilde c
d(\cdot,\cdot)$. We have, by~\eqref{281},
\begin{eqnarray*}d(\gamma^j_{y_h}(t),\gamma^j_{y_h}(t^{\prime})) & \leq & \tilde c
d_c(\gamma^j_{y_h}(t),\gamma^j_{y_h}(t^{\prime}))\leq \tilde c
\int\limits_{t}^{t^{\prime}}|\dot\gamma^j_{y_h}(s)|_0\,ds=2\tilde
c\upsilon^{-1}(y_h)|t-t^{\prime}|\\ & \leq &2 \tilde c
K|t-t^{\prime}|.\end{eqnarray*} Thus, if
$|t-t^{\prime}|\leq\frac{\rho_u}{2\tilde c K}$, then
$$d(\gamma^j_{y_0}(t),\gamma^j_{y_h}(t^{\prime}))\leq
\varpi\big(d(\gamma^j_{y_0}(t),\gamma^j_{y_h}(t))+d(\gamma^j_{y_h}(t),\gamma^j_{y_h}(t^{\prime}))\big)\leq
2\varpi\rho_u$$ and $\gamma^j_{y_h}(t^{\prime})\in
B(\eta_u,2\varpi\rho_u)$ for all $h\geq 1$. Then there exists
$\xi\in f^{-1}(\eta_u)\cap \overline B(0,s^{\prime})$ such that if
$|t-t^{\prime}|<\frac{\rho_u}{2\tilde c K}$, then
$\lambda_{h,\nu}(t^{\prime})\in U(\xi,2\varpi\rho_u)$ for $h\geq
1$. This means that the considered family
$\{\lambda_{h,\nu}\}_{h=1}^{\infty}$ is equicontinuous. Since the
families $\{\lambda_{h,\nu}\}_{h=1}^{\infty}$, $\nu=1,\ldots,m$,
are also uniformly bounded, by the Ascoli theorem, we can find a
subsequence of $\{\lambda_{h,\nu}\}_{h=1}^{\infty}$ that converges
uniformly to the curve $\lambda_{\nu}:\ [0,\sigma_M]\to \overline
B(0,s^{\prime})$. The curve $\lambda_{\nu}$ is a maximal
$f_0$-lifting of $\gamma^j_{y_0}$ starting in $\overline B(0,s)$.

We continue now, studying the $f_0$-liftings of the curve
$\gamma^0_{y_0}(t):\ [\sigma_M,\infty]\to\mathbb C B(0,\sigma_M)$
joining $z_0$ with $\infty$, provided
$\gamma^0_{y_0}(\infty)=\infty$. We show that the family
$\{\lambda_{h,\nu}\}_{h=1}^{\infty}$ is equicontinuous for each
$\nu=1,\ldots,m$. Here $\lambda_{h,\nu}$ are the $f_0$-liftings of
$\gamma^0_{y_h}$ as was described above. Each of the curves
$\lambda_{h,\nu}$ is defined in $[\sigma_M,\infty)$. Fix
$\varepsilon>0$. For every $t\in[\sigma_M,\infty)$, there exists
$\delta(t)>0$ such that, for all $\rho\in(0,\delta(t)]$,
$U(\xi,\rho)=f^{-1}(B(f(\xi),\rho))$ is a normal neighborhood of
$\xi\in f^{-1}(\gamma^0_{y_0}(t))\cap\overline B(0,s^{\prime})$
with $\diam(U(\xi,\rho))<\varepsilon$. We find a ball
$B(0,R)\in\mathbb G$, $R>\sigma_M$, such that $f^{-1}(\mathbb C
B(0,R))=W_1\cup W_2\cup\ldots\cup W_s$, $W_j$ are disjoint normal
domains in $\overline B(0,s^{\prime})$ with
$\diam(W_j)<\varepsilon$, $j=1,\ldots,s$. We cover the
intersection $\gamma^0_{y_0}\cap \overline B(0,R)$ by finite
number of balls $B(\gamma^0_{y_0}(t),\delta(t)/2)$, that we still
denote by $B(\eta_u,\rho_u)$, $u=1,\ldots,v$. We may suppose that
the curves $\gamma^0_{y_h}$, $h\geq 1$, belong to
$\bigcup_{u=1}^{v}B(\eta_u,\rho_u)\cup \mathbb C\overline B(0,R)$.
For every $t\in[\sigma_M,R]$, there exists $u$ such that
$\gamma^0_{y_h}(t)\in B(\eta_u,\rho_u)$. We argue as in the
previous case and deduce that if
$|t-t^{\prime}|<\frac{\rho_u}{2\tilde c K}$,
$t,t^{\prime}\in[\sigma_M,R]$, then
$\lambda_{h,\nu}(t^{\prime})\in U(\xi,2\varpi\rho_u)$ for $h\geq
1$. This means that the considered family is equicontinuous in
$[\sigma_M,R]$. For $t>R$, the curves $\gamma^0_{y_h}(t)$ belong
to $\mathbb C\overline B(0,R)$. Then for each $h>1$ one can find
$W_j$, $\diam(W_j)<\varepsilon$, $j=1,\ldots,s$, such that
$\lambda_{h,\nu}(t)\in W_j$. This proves that the family under
consideration is equicontinuous for $t\in[\sigma_M,\infty]$.
Applying the Ascoli theorem we conclude that the limit curve
$\lambda_{\nu}$ is a maximal $f_0$-lifting of $\gamma^0_{y_0}$
starting in $\overline B(0,s)$.

{\it Case 2.} $m+1\leq\nu\leq g$. Fix $j\in J$. If
$\lambda_{h,\nu}$ is half opened, it extends to the closed curve
in $\overline B(0,s^{\prime} +1)$. Let $\overline\lambda_{h,\nu}:\
[0,t_h]\to\overline B(0,s^{\prime} +1)$ be extended curves. We may
assume (if it necessary passing to the subsequence) that $t_h\to
t_0\in (0,1]$. Let $G_h:\ [0,t_0]\to[0,t_h]$ be an affine mapping
with $G(0)=0$. Arguing as above we deduce that the sequence
$\{\overline\lambda_{h,\nu}\}_{h=1}^{\infty}$ converges uniformly
to $\overline\lambda_{\nu}:\ [0,t_0]\to\overline B(0,s^{\prime}
+1)$ which is lifting of $\gamma_{y_0}^j\vert_{[0,t_0]}$. If
$\Delta\subset[0,t_0]$ is the largest interval such that
$0\in\Delta$ and $\overline\lambda_{\nu}\vert_{\Delta}\subset
B(0,s^{\prime} +1)$, then
$\lambda_{\nu}=\overline\lambda_{\nu}\vert_{\Delta}$ is a maximal
$f_0$-lifting of $\gamma_{y_0}^j$ starting in $\overline B(0,s)$.

It remains to consider the $f_0$-liftings $\lambda_{h,\nu}$ of
$\gamma^0_{y_h}$ for $\nu=m+1,\ldots,g$. Fix~$\nu$. We may assume
that the locus of $\lambda_{h,\nu}$ is not contained in $\overline
B(0,s^{\prime}+2)$. (If it is contained in $\overline
B(0,s^{\prime}+2)$, then we argue as above for $\lambda_{h,\nu}$,
$\nu\leq m$, and find a maximal $f_0$-lifting of $\gamma^0_{y_0}$
starting in $\overline B(0,s)$.) By $t_h$, for every $h$, we
denote  the first value of parameter when $\lambda_{h,\nu}(t)$
intersects $\partial B(0,s^{\prime}+1)$. Then
$[\sigma_M,t_h]\subset[\sigma_M,L]$ for some $\sigma_M<L<\infty$
and for all $h>1$. We extend each $\lambda_{h,\nu}(t)$ to the
closed curve $\overline\lambda_{h,\nu}$ in $\overline
B(0,s^{\prime} +1)$. Assuming that $t_h\to t_0\in (\sigma_M,L]$,
we argue as in the previous paragraf and find the limit sequence
$\lambda_{\nu}=\overline\lambda_{\nu}\vert_{\Delta}$ that is a
maximal $f_0$-lifting of $\gamma^0_{y_0}$ starting in $\overline
B(0,s)$. Here $\Delta\subset[\sigma_M,t_0]$ is the largest
interval such that $\sigma_M\in\Delta$ and
$\overline\lambda_{\nu}\vert_{\Delta}\subset B(0,s^{\prime} +1)$.

If we show that the constructed limit sequence
$\Lambda_0=(\lambda_1,\ldots,\lambda_g)$ is essentially separate,
then we conclude that $\Lambda_0\in\Omega_{y_0}$. Fix a point
$x\in B(0,s^{\prime})$. Let $A=\{\nu:\
\lambda_{\nu}(t)=x\}\neq\emptyset$ and $U(x,r)$ be a normal
neighborhood of $x$. Fix $h_0$ such that $\lambda_{h,\nu}\cap
U(x,r)\neq\emptyset$ for all $h\geq h_0$ and $\nu\in A$. We choose
$h\geq h_0$ and, then, find a point
$\eta=\gamma^j_{y_h}(t^{\prime})$ in $\bigcap\limits_{\nu\in
A}f(\lambda_{h,\nu}\cap U(x,r))$. Let $\xi_1,\ldots,\xi_w$ be the
points in $\{\lambda_{h,\nu}(t^{\prime}):\ \nu\in A\}\subset
f^{-1}(\eta)\cap U(x,r)$. Since the curves
$\lambda_{h,1},\ldots,\lambda_{h,g}$ are essentially separate, we
have $$\theta_u=\card\{\nu:\
\lambda_{h,\nu}(t^{\prime})=\xi_u\}\leq i(\xi_u,f),\qquad
u=1,\ldots,w.$$ Hence $$\card
A=\sum\limits_{u=1}^{w}\theta_u\leq\sum\limits_{u=1}^{w}i(\xi_u,f)\leq
i(x,f).$$ The claim is proved.

In the limit sequence $\Lambda_0=(\lambda_1,\ldots,\lambda_g)$ we
have $\lambda_{\nu}\subset\overline B(0,s^{\prime})$ for
$1\leq\nu\leq m$ with guarantee. By the limit process, it can
happened that $\lambda_{\nu}\subset\overline B(0,s^{\prime})$ for
some other $\nu>m$. Thus, we have $\psi^j(y_0)\geq
N(\Lambda_0)\geq m$. The lemma is proved.
\end{proof}

Since $\mes(\partial U)=0$ and a quasimeromorphic mapping
possesses the Luzin property, we have $\mes(f(\partial U))=0$.
Then $\sigma\big(S(0,t)\cap f(\partial U)\big)=0$ for almost all
$t>0$. If it holds for $t=\sigma_M$, then $\psi^j(y)$ is a
measurable function. If it is not so, then we choose $t$
sufficiently close to $\sigma_M$ such that Theorem~\ref{44} still
holds. Thus, from now on we can think, that
$\sigma\big(S(0,\sigma_M)\cap f(\partial U)\big)=0$ and
$\psi^j(y)$ is a measurable function.

We are interested in estimating of the quantity of $f_0$-liftings
starting from different sets $U_i$. For this we define the
functions $\psi^j_i(y)$ which will calculate the number of
$f_0$-liftings of $\gamma_y^j$ starting in $B(0,s)\setminus U_i$.
For sequences $\Lambda\in\Omega(y)$, we let
$$N(i,\Lambda)=\card\Big\{\nu:\ \lambda_{\nu}\subset\overline B(0,s^{\prime}),
\ \ \lambda_{\nu}
\quad\text{starts in}\quad B(0,s)\setminus U_i\Big\}$$ and
$$\psi^j_{i}(y)=\sup\{N(i,\Lambda):\ \Lambda\in \Omega(y)\}.$$
By $\Omega(i,y)$, we denote the set of sequences $\Lambda\in\Omega(y)$ for which
$\psi^j_{i}=N(i,\Lambda)$.

\begin{lemma}\label{122}
The functions $\psi^j_{i}:\ S\setminus\mathcal Z\to\mathbb N$,
$i=1,\ldots,p$, are measurable.
\end{lemma}

\begin{proof}
We prove that $\psi^j_{i}$ is measurable. Fix $j\in J\cup\{0\}$
and $i\in I$. Since we know that $\psi^j$ is measurable, it is
enough to show that the restriction of $\psi^j_{i}$ to each set
$$A_m=\{y\in S\setminus\mathcal Z: \psi^j(y)=m\},\qquad m=0,...,m_{max},$$ is measurable.
Here $m_{max}=\max\{n(s,z),\ z\in S(0,\sigma_M)\setminus\mathcal
Z\}$. Fix $m$ and verify that $\psi^j_{i}$ is upper semicontinuous
in $B_m=A_m\setminus f(\partial U)$. Let $y_0\in B_m$ and
$\{y_h\}\in B_m$ be a sequence that converges to $y_0$. We may
assume that for some integer $m_1\leq m$ we have $\psi^j_{i}=m_1$
for all $h\geq 1$ and the following properties hold.
\begin{itemize}
\item[(i)]{There exist normal neighborhoods $V_1,\ldots,
V_l\subset U$ of points $(x_1,\ldots,x_l)=f^{-1}(z_0)\cap U$ such
that $z_h\in f(V_1)\cap\ldots\cap f(V_l)$ for all $h\geq 1$. Here
$z_0=\gamma^j_{y_0}(0)$, $z_h=\gamma^j_{y_h}(0)$ for $j\in J$ or
$z_0=\gamma^0_{y_0}(\sigma_M)$, $z_h=\gamma^0_{y_h}(\sigma_M)$.}
\item[(ii)]{There is a maximal essentially separate sequence
$\Lambda_h=(\lambda_{h,1},\ldots,\lambda_{h,g})\in\Omega(i,z_h)$,
$g=n(s,z_0)$, such that every $\lambda_{h,\nu}$ starts at a point
in $f^{-1}(z_h)\cap V_{\mu}$ for some $\mu$ dependent on $\nu$ and
independent of $h$. Moreover, the locus of every curve
$\lambda_{h,\nu}$ belongs to $\overline B(0,s^{\prime})$ for
$1\leq\nu\leq m_1$.}\item[(iii)]{The bound $\upsilon(y_h)^{-1}\leq
K$ holds for all $h=0,1,\ldots$.}
\end{itemize}
As in the proof of Lemma~\ref{121} we find a maximal essentially
separate sequence $\Lambda_0\in\Omega_{y_0}$ such that
$N(\Lambda_0)\geq m$. Since $y_0\in A_m$, we have
$N(\Lambda_0)\leq \psi^j(y_0)=m$. Thus $N(\Lambda_0)=m$ and
therefore $\Lambda_0\in\Omega(y_0)$. The construction of
$\Lambda_0$ implies that $\psi^j_{i}(y_0)\geq N(i,\Lambda_0)\geq
m_1$. It means the upper semicontinuity of $\psi^j_{i}$.
\end{proof}

\subsection{Estimates for $\Delta_j$}

If the sum $\sum\limits_{j\in
J\cup\{0\}}\Delta_j=\sum\limits_{j\in
J\cup\{0\}}1-\frac{n(s,a_j)}{\nu(s,\sigma_M)}$ is bounded, then it
is nothing to prove. In this section we show that if the sum
$\sum\limits_{j}\Delta_j$ is large we come to the situation when
Lemma~\ref{111} can be applied. We may assume that
$\sum\limits_{j}\Delta_j> 20M$.

Let $i\in I=\{1,\ldots,p\}$ and $j\in
J\cup\{0\}=\{0,1,\ldots,q\}$. For each $y\in S\setminus\mathcal Z$
and $U_i\in B(0,s)$ we choose a maximal essentially separate
$f_0$-liftings $\Lambda_y=(\lambda_{y,1},\ldots,\lambda_{y,g})$,
$g=n(s,z)$, $\Lambda_y\in\Omega(i,y)$. Here $z=\gamma^j_y(0)$ for
$j\in J$ or $z=\gamma^0_y(\sigma_M)$. Those curves
$\lambda_{y,\nu}$ in $\Lambda_y$ that start in $U_i$ form a
maximal sequence $(\beta_1,\ldots,\beta_m)$, $m=n(U_i,z)$, of
essentially separate $f_0$-liftings of $\gamma^j_y$. For this
sequence we consider the sequence $\alpha_1,\ldots,\alpha_{\mu_y}$
of those liftings of $\beta_{\mu}$ whose locus $|\beta_{\mu}|$ is
not contained in $\overline B(0,s^{\prime})$. From now on, we
denote the quantity of such $\alpha_{\mu}$ by $n_{i}^{j}(y)$. The
number of curves $\beta_{\mu}$, starting on $U_i$ and the locus of
which is contained in $\overline B(0,s^{\prime})$, equals the
difference $\psi^{j}(y)-\psi_{i}^j(y)$. We have
$$n(U_i,z)=n^{j}_{i}(y)+(\psi^{j}(y)-\psi_{i}^j(y)).$$ Since the
functions $\psi^j$, $\psi^j_i$ are measurable by
Lemmas~\ref{121},~\ref{122}, and $n(U_i,z)$ are upper
semicontinuous, we get that $n^{j}_{i}(y)$ are measurable for
arbitrary $i$ and $j$.

For $i\in I$, set \begin{equation}\label{249}J_i=\Big\{j:\
\frac{1}{\kappa(\mathbb G,Q)}\int\limits_{S\setminus\mathcal
Z}n_i^j(y)\,d\sigma>\frac{1}{2M}\nu(U_i,\sigma_M)\Delta_j\Big\},\end{equation}
where $M$ is the constant of multiplicity in the decomposition of
$B(0,s)$. Observe, that the overlapping of the decomposition with
the multiplicity $M$ imply the following inequality
\begin{equation}\label{265}\frac{1}{M}\sum\limits_{i=1}^{p}\nu(U_i,\sigma_M)\leq
\nu(s,\sigma_M)\leq\sum\limits_{i=1}^{p}\nu(U_i,\sigma_M).\end{equation}

We start to estimate the sum $\sum\Delta_j$ with the next lemma.

\begin{lemma}\label{120}
The function $n_i^j(y)$ satisfy
$$\sum\limits_{i\in I}\sum\limits_{j\in J_i}\int\limits_{S\setminus\mathcal
Z}n_i^j(y)\,d\sigma\geq \frac{\kappa(\mathbb
G,Q)}{2}\nu(s,\sigma_M)\sum\limits_{j\in J\cup\{0\}}\Delta_j.$$
\end{lemma}

\begin{proof}
Summing the functions $n^{j}_{i}$ over $i$, we obtain
$$\sum\limits_{i\in I}n_i^j(y)\geq\sum\limits_{i\in I}n(U_i,z)-\psi^j(y)
\geq n(s,z)-n(s^{\prime},a_j).$$ Integrating $n_i^j(y)$ over
$S\setminus\mathcal Z$ we have
$$\frac{1}{\kappa(\mathbb
G,Q)}\int\limits_{S\setminus\mathcal Z}\sum\limits_{i\in
I}n_i^j(y)\,d\sigma\geq\nu(s,\sigma_M)-n(s^{\prime},a_j)
=\nu(s,\sigma_M)\Delta_j.$$ From the last estimate,~(\ref{249}),
and~\eqref{265} we deduce
\begin{align*}
 \frac{\kappa(\mathbb G,Q)}{2} &\nu(s,\sigma_M) \sum\limits_{j\in
J\cup\{0\}}\Delta_j \geq \frac{\kappa(\mathbb
G,Q)}{2M}\sum\limits_{i\in I}\sum\limits_{j\in
J\cup\{0\}}\nu(U_i,\sigma_M)\Delta_j
\\ & \geq
\frac{\kappa(\mathbb G,Q)}{2M}\sum\limits_{i\in
I}\sum\limits_{j\in (J\cup\{0\})\setminus
J_i}\nu(U_i,\sigma_M)\Delta_j
 \geq
\sum\limits_{i\in I}\sum\limits_{j\in (J\cup\{0\})\setminus
J_i}\int\limits_{S\setminus\mathcal Z}n_i^j(y)\,d\sigma
\\ & =
\sum\limits_{i\in I}\sum\limits_{j\in
J\cup\{0\}}\int\limits_{S\setminus\mathcal Z}n_i^j(y)\,d\sigma-
\sum\limits_{i\in I}\sum\limits_{j\in
J_i}\int\limits_{S\setminus\mathcal Z}n_i^j(y)\,d\sigma
\\ & \geq
\kappa(\mathbb G,Q)\nu(s,\sigma_M)\sum\limits_{j\in
J\cup\{0\}}\Delta_j - \sum\limits_{i\in I}\sum\limits_{j\in
J_i}\int\limits_{S\setminus\mathcal Z}n_i^j(y)\,d\sigma.
\end{align*}
\end{proof}

Now, we try to estimate the average number over $S(0,\sigma_M)$ of
$f_0$-liftings of parts of $\gamma_y^j$ such as
$\alpha^{(1)}_{\nu}$ and $\alpha^{(2)}_{\nu}$. This helps us to
refine the upper bound for $\sum\limits_{j\in J\cup
\{0\}}\Delta_j$. We set
$$J^i=\Bigl\{j\in J_i: \int\limits_{S\setminus\mathcal Z}n_i^j(y)\,d\sigma\leq
3\int\limits_{S\setminus\mathcal Z}l_i^j(y)\,d\sigma\quad
\text{or}\quad \int\limits_{S\setminus\mathcal
Z}n_i^j(y)\,d\sigma\leq 3 \int\limits_{S\setminus\mathcal
Z}m_i^j(y)\,d\sigma\Bigr\}.$$ Derive the bound on average number
over $S(0,\sigma_M)$ for indexes $J^i$ from estimates for curves
$\alpha^{(1)}_{\nu}$ and $\alpha^{(2)}_{\nu}$. By Lemmas~\ref{118}
and~\ref{119}, we deduce
\begin{eqnarray}\label{269}\sum_{j\in J^i}\int\limits_{S\setminus\mathcal Z}n_i^j(y)\,d\sigma &
\leq & 3\sum_{j\in J^i}\int\limits_{S\setminus\mathcal
Z}(l_i^j(y)+m^j_i(y))\,d\sigma\nonumber
\\ & \leq &
3c_{5}K_{I}(f)q\Big(\ln\frac{\sigma_M}{\varsigma_0}\Big)^{Q-1}+3c_{6}K_I(f)A_i\nu(U_i,\sigma_M),
\end{eqnarray}
where $A_i$ are constants from~(\ref{250}).

We continue to estimate, making use of Lemma~\ref{120}, the
inequality~\ref{269} and upper bound for $p$: $p\leq
2C_3\varepsilon_0^{Q-1}\nu(s,\sigma_M)$,
\begin{eqnarray}\label{253}\sum_{i\in I}\ \sum\limits_{j\in J_i\setminus J^i}
\int\limits_{S\setminus\mathcal Z}n_i^j(y)\,d\sigma & = &
\sum\limits_{i\in I}\sum\limits_{j\in
J_i}\int\limits_{S\setminus\mathcal
Z}n_i^j(y)\,d\sigma-\sum\limits_{i\in I}\sum\limits_{j\in
J^i}\int\limits_{S\setminus\mathcal Z}n_i^j(y)\,d\sigma\nonumber
\\ & \geq &
\frac{\kappa(\mathbb G,Q)}{2}\nu(s,\sigma_M)\sum\limits_{j\in
J\cup\{0\}}\Delta_j\nonumber
\\ & - &
6C_3c_5K_I(f)q\Big(\ln\frac{\sigma_M}{\varsigma_0}\Big)^{Q-1}\varepsilon_0^{Q-1}\nu(s,\sigma_M)
\\
& - & 3c_{6}K_I(f)\sum_{i\in I}A_i\nu(U_i,\sigma_M).\nonumber
\end{eqnarray}

We need to choose the constants $A_i$ to obtain an effective lower
bound for the sum $\sum\limits_{i\in I}\sum\limits_{j\in
J_i\setminus J^i} \int\limits_{S\setminus\mathcal
Z}n_i^j(y)\,d\sigma$. We introduce the value
$\lambda_i=\card(J_i\setminus J^i)$. The set $J_i\setminus J^i$
contains the indexes $j\in J$ for which the following inequalities
hold:
\begin{eqnarray}\label{251}\int\limits_{S\setminus\mathcal Z}n_i^j(y)\,d\sigma & > &
\frac{\kappa(\mathbb G,Q)}{2M}\nu(U_i,\sigma_M)\Delta_j,\nonumber
\\
\int\limits_{S\setminus\mathcal Z}n_i^j(y)\,d\sigma & > &
3\int\limits_{S\setminus\mathcal Z}l_i^j(y)\,d\sigma,
\\
\int\limits_{S\setminus\mathcal Z}n_i^j(y)\,d\sigma & > &
3\int\limits_{S\setminus\mathcal Z}m_i^j(y)\,d\sigma.\nonumber
\end{eqnarray}

We recall the definition of $A_i$:
$A_i\nu(U_i,\sigma_M)=\Big(\ln\frac{\varsigma_0}{\sigma_i}\Big)^{Q-1}$,
where $\sigma_i\in (0,\varsigma_0]$. When $A_i$ increases in the
range $[0,\infty)$ the number $\lambda_i$ decreases from
$\lambda_i^0$ to some value $\lambda_i^{\infty}$. From~(\ref{251})
we have
$$
\lambda_i^0=\card\Bigl\{j\in J_i:\ \int\limits_{S\setminus\mathcal
Z}n_i^j(y)\,d\sigma > 3\int\limits_{S\setminus\mathcal
Z}l_i^j(y)\,d\sigma\Bigr\}.
$$
We may assume that the jumps at the discontinuities of the
function $A_i\mapsto\lambda_i$ equal~$1$. If it is not so, we can
make a small variation of function $m_i^j$ changing $\sigma_i$ for
different $j$'s. Finally, we conclude that we can choose $A_i\geq
0$ with
\begin{equation}\label{252} \lambda_i-1\leq
9\kappa^{-1}(\mathbb G,Q)c_6K_I(f)A_i\leq\lambda_i.
\end{equation}

We use the value of the constants $A_i$ to terminate the
estimation of~(\ref{253}). From~(\ref{253}) and~(\ref{252}) we
obtain
\begin{eqnarray}\label{254}\sum_{i\in I}\lambda_i\nu(U_i,\sigma_M)
& = &\sum_{i\in I}\sum_{j\in J_i\setminus J^i}\nu(U_i,\sigma_M)
\geq \sum_{i\in I}\sum_{j\in J_i\setminus
J^i}\frac{1}{\kappa(\mathbb G,Q)}\int\limits_{S\setminus\mathcal
Z}n_i^j(y)\,d\sigma\nonumber
\\ & \geq &
\frac{\nu(s,\sigma_M)}{2}\sum\limits_{j\in
J\cup\{0\}}\Delta_j\nonumber
\\ & - &\frac{6C_3c_5K_I(f)q\varepsilon_0^{Q-1}}{\kappa(\mathbb G,Q)}
\Big(\ln\frac{\sigma_M}{\varsigma_0}\Big)^{Q-1}\nu(s,\sigma_M) \\
& - & \frac{1}{3}\sum_{i\in I}\lambda_i\nu(U_i,\sigma_M).\nonumber
\end{eqnarray}
We are free in the choice of $\varepsilon_0$. Take $\varepsilon_0$
such that
\begin{equation}\label{278}
\frac{6C_3c_5K_I(f)q\varepsilon_0^{Q-1}}{\kappa(\mathbb
G,Q)}\Big(\ln\frac{\sigma_M}{\varsigma_0}\Big)^{Q-1}<M.\end{equation}
Then~(\ref{254}) implies
\begin{equation}\label{257}\sum_{i\in
I}\lambda_i\nu(U_i,\sigma_M)\geq\frac{\nu(s,\sigma_M)}{4}
\sum\limits_{j\in J\cup\{0\}}\Delta_j.\end{equation}

Now we exclude some set of indexes from $I$ where
$\sum\lambda_i\nu(U_i,\sigma_M)$ can be bound from above by
$\sum_{j\in J\cup\{0\}}\Delta_j$. Namely, let
$$I_1=\Bigl\{i\in I;\
\lambda_i\leq \frac{\sum_{j\in
J\cup\{0\}}\Delta_j}{10M}\quad\text{or}\quad \nu(U_i,\sigma_M)\leq
P\Bigr\}.$$ Here $P=\frac{c_1K_I(f)}{(\ln
2)^{Q-1}}(|\ln\frac{\sigma_M^2}{\sigma_m}|^{Q-1}+|\ln
\sigma_M|^{Q-1}+2c_0)$ with $\sigma_m$ defined in Lemma~\ref{111}.
We need the first choice for using Lemmas~\ref{111} and~\ref{112}.
The second one serves for applying of Lemma~\ref{117}.

\begin{lemma}\label{123} Under the previous notations we have
\begin{equation}\label{256}\sum_{i\in I\setminus
I_1}\lambda_i\nu(U_i,\sigma_M)\geq\frac{\nu(s,\sigma_M)}{8}
\sum\limits_{j\in J\cup\{0\}}\Delta_j.\end{equation}
\end{lemma}

\begin{proof}
Summing over $I_1$ we get \begin{eqnarray*}\sum_{i\in
I_1}\lambda_i\nu(U_i,\sigma_M) & \leq &
\frac{1}{10}\nu(s,\sigma_M)\sum\limits_{j\in
J\cup\{0\}}\Delta_j+Pqp
\\
& \leq & \frac{1}{10}\nu(s,\sigma_M)\sum\limits_{j\in
J\cup\{0\}}\Delta_j+Pq2C_3\varepsilon_0^{Q-1}\nu(s,\sigma_M)\end{eqnarray*}
from\eqref{265}. Let us add one more restriction on
$\varepsilon_0$. We choose $\varepsilon_0$ such that
\begin{equation}\label{279} 2PqC_3\varepsilon_0^{Q-1}<M/2.\end{equation} Then
\begin{equation}\label{255}\sum_{i\in
I_1}\lambda_i\nu(U_i,\sigma_M)\leq \frac{\nu(s,\sigma_M)}{8}
\sum\limits_{j\in J\cup\{0\}}\Delta_j.\end{equation} Joining the
estimates~(\ref{257}) and~(\ref{255}) we deduce~(\ref{256}).
\end{proof}

\begin{lemma}\label{125}
If $\sum\limits_{j\in J\cup\{0\}}\Delta_j>20M$, then, for $i\in
I\setminus I_1$, the following inequality
\begin{equation}\label{264}\lambda_i\nu(U_i,\sigma_M)\leq
c_{7}K_O(f)K_I(f)\nu(Z_i,\sigma_M).\end{equation} holds.
\end{lemma}

\begin{proof}
Let $i\in I\setminus I_1$. The inequalities~(\ref{251}) imply for
$j\in J_i\setminus J^i$ that
$$\int\limits_{S\setminus\mathcal Z}(n_i^j-l_i^j-m_i^j)(y)\,d\sigma
\geq\frac{1}{3}\int\limits_{S\setminus\mathcal Z}n_i^j(y)\,d\sigma
>\frac{\kappa(\mathbb G,Q)}{6M}\nu(U_i,\sigma_M)\Delta_j>0.$$ We conclude that
the function $n_i^j(y)-l_i^j(y)-m_i^j(y)$ is positive for some
$y\in S\setminus\mathcal Z$. Therefore
$$0<n_i^j(y)-l_i^j(y)-m_i^j(y)<n_i^j(y)-\card(M^j_i(y)\cup
L^j_i(y))$$ by Lemmas~\ref{118} and~\ref{119}. It means that for
these indexes $i$ and $j\in (J_i\setminus J^i)\setminus\{0\}$
there are indexes $\nu\in\{1,\ldots,n_i^j(y)\}$ such that we have
$\frac{\sigma_M-u_{y,\nu}}{\sigma_M-v_{y,\nu}}>\frac{\sigma_M}{\varsigma_0}$
and
$\frac{\sigma_M-v_{y,\nu}}{\sigma_M-w_{y,\nu}}>\frac{\varsigma_0}{\sigma_i}$.
We also can say that the corresponding curves has left $W_i$ and
reached $\partial B(0,s^{\prime})$. Here the values $u_{y,\nu}$,
$v_{y,\nu}$, $w_{y,\nu}$ were defined in~(\ref{258}). The two last
inequalities imply that $\sigma_M-w_{y,\nu}\leq\sigma_i$ and for
the corresponding $f_0$-lifting $\alpha_\nu$ of $\gamma^j_y$ the
following properties hold
\begin{itemize}
\item[(i)]{the restriction
$\alpha_i^j=\alpha_{\nu}\vert_{[w_{y,\nu},t]}$ is a curve in $X_i$
connecting $\partial W_i$ and $\partial X_i$ for some $t\leq
\sigma_M$;} \item[(ii)]{$f(\alpha^j_{i})\subset
B(a_j,3\sigma_i/2)$.}
\end{itemize}

If $j\in (J_i\setminus J^i)\setminus\{0\}$, then we apply
Lemma~\ref{111} to the mapping $f_0$ and the set
$F_j=\alpha_i^{j}$. We consider the balls $W_i\subset X_i\subset
Y_i$ with radii $r=4r_i$, $\rho=6r_i$, and $\theta\rho$, where
$\theta$ was defined in the decomposition of $B(0,s)$. We obtain
\begin{eqnarray}\label{259} \Big(M_Q(\Gamma_j)&-&\frac{\kappa(\mathbb G,Q)K_O(f)K_I(f)}{(\ln
\theta)^{Q-1}}\Big)
\Big(\ln\frac{2\sigma_m}{3\sigma_i}\Big)^{Q-1}\nonumber
\\ & \leq &\kappa(\mathbb G,Q)K_O(f)\nu\big(6\theta r_i,
S(a_j,\sigma_m)\big).\end{eqnarray} Note, that we use the
relations $\sigma_i\leq\varsigma_0<\frac{2}{3}\sigma_m$. Thanks to
our choice of $r$ and $\rho$ we have $M_Q(\Gamma_j)\geq c(Q)\ln
\frac{\rho}{r}=c(Q)\ln\frac{3}{2}$. By Remark~\ref{r4} we get
\begin{equation}\label{260}
c(Q)\ln\frac{5}{4}\Big(\ln\frac{2\sigma_m}{3\sigma_i}\Big)^{Q-1}\leq
\Big(\ln\frac{3}{2}-\frac{\kappa(\mathbb
G,Q)K_O(f)K_I(f)}{c(Q)(\ln \theta)^{Q-1}}\Big)
\Big(\ln\frac{2\sigma_m}{3\sigma_i}\Big)^{Q-1}.\end{equation}
Making use of Lemma~\ref{117} with $Z=S(0,\sigma_M)$,
$Y=S(a_j,\sigma_m)$, $r=6\theta r_i$ and $\vartheta=2$ we deduce
\begin{eqnarray}\label{261}\nu\big(6\theta r_i, S(a_j,\sigma_m)\big) &
\leq & \nu\big(12\theta r_i, S(0,\sigma_M)\big)\nonumber
\\ & + & \frac{c_1K_I(f)}{(\ln 2)^{Q-1}}\big(|\ln
\sigma_m|^{Q-1}+|\ln \sigma_M|^{Q-1}+2c_0\big)
\\ & \leq &
\nu(Z_i,\sigma_M)+P=2\nu(Z_i,\sigma_M).\nonumber\end{eqnarray} The
last inequality was possible because of the choice of
$P<\nu(U_i,\sigma_M)\leq\nu(Z_i,\sigma_M)$ for $i\in I\setminus
I_1$ and $\theta=\varkappa/3$. We conclude
\begin{equation}\label{262}
\Big(\ln\frac{2\sigma_m}{3\sigma_i}\Big)^{Q-1} \leq
\frac{2\kappa(\mathbb G,Q)}{c(Q)\ln 5/4}K_O(f)\nu(Z_i,\sigma_M)
\end{equation}
from~(\ref{259}),~(\ref{260}), and~(\ref{261}). Now, we use the
definition~(\ref{252}) of $A_i$. We can  assume also that
$\lambda_i\geq 2$. We have
\begin{equation}\label{263}\frac{\kappa(\mathbb
G,Q)\lambda_i\nu(U_i,\sigma_M)}{18c_6K_{I}(f)}\leq
A_i\nu(U_i,\sigma_M)
=\Big(\ln\frac{\varsigma_0}{\sigma_i}\Big)^{Q-1}
<\Big(\ln\frac{2\sigma_m}{3\sigma_i}\Big)^{Q-1}.\end{equation}
Finally,~(\ref{262}) and ~(\ref{263}) imply~(\ref{264}).

It remains to consider the case $j=0$. We should slightly change
the arguments. If $0<n_i^0(y)-\card(M^0_i(y)\cup L^0_i(y))$ then
there are indexes $\nu\in(1,\ldots,n_i^0(y))$ such that
$\frac{v_{y,\nu}}{u_{y,\nu}}>\frac{\sigma_M}{\varsigma_0}$ and
$\frac{w_{y,\nu}}{v_{y,\nu}}>\frac{\varsigma_0}{\sigma_i}$, where
$\sigma_M\leq u_{y,\nu}<v_{y,\nu}<w_{y,\nu}<\infty$. We deduce
$$w_{y,\nu}>\frac{v_{y,\nu}\varsigma_0}{\sigma_i}>\frac{u_{y,\nu}\sigma_M}{\sigma_i}
>\frac{\sigma_M^2}{\sigma_i}>\frac{2\sigma_M^2}{3\sigma_i}>\frac{\sigma_M^2}{\sigma_m}>\sigma_M.$$
For the corresponding $f_0$-lifting $\alpha_{\nu}$ the following
properties hold for some $t>w_{y,\nu}$
\begin{itemize}
\item[(i)]{the restriction
$\alpha_i^0=\alpha_{\nu}\vert_{[w_{y,\nu},t]}$ is a curve in $X_i$
connecting $\partial W_i$ and $\partial X_i$;}
\item[(ii)]{$f(\alpha^0_{i})\subset \mathbb
CB(0,\frac{2\sigma_M^2}{3\sigma_i})$.}
\end{itemize}
Instead of Lemma~\ref{111} we apply Lemma~\ref{112} to the sets
$F_0=\alpha_i^0$, $F_j=\alpha_j^i$ in the balls $W_i\subset
X_i\subset Y_i$, $t=\frac{2\sigma_M^2}{3\sigma_i}$,
$s=\frac{\sigma_M^2}{\sigma_m}$. We get
\begin{eqnarray}\label{266}\Big(M_Q(\Gamma_0) & - & \frac{\kappa(\mathbb G,Q)K_O(f)K_I(f)}{(\ln
\theta)^{Q-1}}\Big)
\Big(\ln\frac{2\sigma_m}{3\sigma_i}\Big)^{Q-1}\nonumber
\\
& \leq & \kappa(\mathbb G,Q)K_O(f)\nu\Big(6\theta r_i,
S\big(0,\frac{\sigma_M^2}{\sigma_m}\big)\Big).\end{eqnarray} Then
we estimate $M_Q(\Gamma_0)$ from below $M_Q(\Gamma_0)\geq c(Q)\ln
3/2$ and deduce~(\ref{260}). Applying Lemma~\ref{117} with
$Z=S(0,\sigma_M)$, $Y=S(0,\frac{\sigma_M^2}{\sigma_m}))$,
$r=6\theta r_i$ and $\vartheta=2$ we have
\begin{eqnarray}\label{267}\nu\big(6\theta r_i, S(0,\frac{\sigma_M^2}{\sigma_m}))\big) &
\leq & \nu\big(12\theta r_i, S(0,\sigma_M)\big)\nonumber
\\ & + & \frac{c_1K_I(f)}{(\ln 2)^{Q-1}}\big(|\ln
\sigma_m|^{Q-1}+\big|\ln
\frac{\sigma_M^2}{\sigma_m}\big|^{Q-1}+2c_0\big)
\\ & \leq &
2\nu(Z_i,\sigma_M).\nonumber\end{eqnarray}
Combining~(\ref{266}),~(\ref{260}), and~(\ref{267}) we
obtain~(\ref{262}) and continue as in the previous case.
\end{proof}

\section{Proof of Theorem~\ref{45}}
We recall that at the end of Subsection~\ref{Constr} we observed
that for the proof of the main theorem it is sufficient to show
the finiteness of the sum $\sum\limits_{j\in J\cup\{0\}}\Delta_j$.
If we assume that $\sum\limits_{j\in J\cup\{0\}}\Delta_j>20M$,
then we can apply Lemma~\ref{125}. We get
\begin{eqnarray*}\frac{\nu(s,\sigma_M)}{8} \sum\limits_{j\in J\cup\{0\}}\Delta_j & \leq & \sum_{i\in
I\setminus I_1}\lambda_i\nu(U_i,\sigma_M)\leq
c_{7}K_O(f)K_I(f)\sum_{i\in I}\nu(Z_i,\sigma_M) \\
& = & c_7\widetilde MK_O(f)K_I(f)\nu(s^{\prime},\sigma_M) < 3/2c_8
K_O(f)K_I(f)\nu(s,\sigma_M)\end{eqnarray*} by Lemma~\ref{123}, the
inequality~\eqref{224}, the estimate~\eqref{283}, and the
decomposition of $B(0,s^{\prime})$.

The final conclusion is that $$\sum\limits_{j\in
J\cup\{0\}}\Delta_j\leq C(Q,K_O(f),K_I(f)).$$ We proved
Theorem~\ref{45}.

\section{Proof of Theorem~\ref{46}}

Let $\mathbb G$ be a $\mathbb H$-type Carnot group and let
$B(0,1)$ be the unit ball in the group~$\mathbb G$. We recall the
statement of Theorem~\ref{46}: {\it Let
$f:B(0,1)\to\overline{\mathbb G}$ be a nonconstant
$K$-quasimeromorphic mapping such that
\begin{equation*}
\limsup\limits_{r\to 1}(1-r)A(r)^{\frac{1}{Q-1}}=\infty.
\end{equation*}
Then there exists a set $E\subset(0,1)$ satisfying
\begin{equation*}
\liminf\limits_{r\to 1}\frac{\mes_1(E\cap [r,1))}{(1-r)}=0.
\end{equation*}
and a constant $C(Q,K)<\infty$ such that
\begin{equation*}
\lim\limits_{r\to 1}\sup\limits_{r\notin E}\sum\limits_{j=0
}^{q}\Big(1-\frac{n(r,a_j)}{\nu(r,1)}\Big)_+\leq C(Q,K),
\end{equation*}
whenever $a_0,a_1,\ldots,a_q$ are distinct points in
$\overline{\mathbb G}$.}

To prove Theorem~\ref{46} we need an analogue of Theorem~\ref{44}
and a construction of the exceptional set $E$. In spite of the
different definition of the function $A(r)$, the proof of the next
lemma repeats the proof from~\cite{Rick1} almost verbatim, because
we have used the continuity of $A(r)$ only. We present the proof
for the completeness.

\begin{lemma}\label{127}
Let $B(0,1)\in \mathbb G$ and $f:B(0,1)\to\overline{\mathbb G}$ be
a quasimeromorphic mapping with the property that
\begin{equation}\label{275}
\limsup\limits_{r\to 1}(1-r)A(r)^{\frac{1}{Q-1}}=\infty.
\end{equation}
Then there exists a set $E\subset[0,1)$ such that
\begin{equation}\label{276}
\liminf\limits_{s\to 1}\frac{\mes_1(E\cap[s,1))}{(1-s)}=0
\end{equation}
and such that the following is true: If $\varepsilon_0\in(0,1/5)$
and if for $s\in(0,1)$ we write
$$s^{\prime}=s+\frac{s}{\varepsilon_0A(s)^{1/(Q-1)}},$$
then there exists an increasing function
$\omega:[0,\infty[\to[D_{\varepsilon_0},1)$ such that for any
sphere $Y=S(w,t)$ in $\mathbb G$ and any
$s^{\prime}\in[\omega(|\log t|),1)\setminus E$ there is an
$s\in(0,1)$, for which inequality~\eqref{218}$:$
\begin{equation*}
\Big|\frac{\nu(s,Y)}{A(s^{\prime})}-1\Big|<\varepsilon_0
\end{equation*}
and inequality~\eqref{219}$:$
\begin{equation*}
\frac{\nu(s,Y)}{\nu(s^{\prime},Y)}>1-\varepsilon_0.
\end{equation*} hold.
Moreover $A(D_{\varepsilon_0})>\varepsilon^{2-2Q}$.
\end{lemma}

\begin{proof}
We inductively define an increasing sequence $t_3,t_4,t_5,\ldots$
in the interval $(0,1)$ tending to $1$ such that
$1-t_m<\frac{(1-t_{m-1})}{m}$. At the same time we define a
sequence $E_3,E_4,\ldots$ of subsets of $(0,1)$. Set $t_3=3/4$ and
$E_3=[3/4,1)$. Suppose that $m\geq 4$ and $t_{m-1}$ is defined.
Set $\phi(r)=\frac{A(r)^{1/(Q-1)}}{m^2}$. There exists $t_m$ such
that
\begin{equation}\label{277}
1-t_m<\frac{1-t_{m-1}}{m}\quad\text{and}\quad
(1-t_m)A(t_m)^{1/(Q-1)}>\frac{8m^3}{1-(\frac{m-1}{m})^{1/(Q-1)}}
\end{equation} by~\eqref{275}. Set
$t^{\star}_m=t_m+\frac{1-t_m}{m}$, $\overline
t_m=1-\frac{1-t_m}{m}$, and
$$F_m=\Big\{r\in(t_m,1):\ A\Big(r+\frac{2r}{\phi(r)}\Big)>\frac{m}{m-1}A(r)
\quad\text{or}\quad r+\frac{2r}{\phi(r)}>1\Big\}.$$ Let us assume
that $F_m\cap(t_m,t_m^{\star}]\neq\emptyset$ and put
$t_m=r_0^{\prime\prime}$. Now we define inductively the sequence
$r_0^{\prime\prime}\leq r_1<r_1^{\prime\prime}\leq
r_2<r_2^{\prime\prime}\ldots\leq r_h<r_h^{\prime\prime}$ of
numbers in $[t_m,1)$ by
$$r_k=\inf\{r\in(r_{k-1}^{\prime\prime},1):r\in F_m\},\qquad r_k^{\prime\prime}=r_k+\frac{2r_k}{\phi(r_k)}.$$
Here the number $h$ is the last index $k$ for which
$F_m\cap(r_{k-1}^{\prime\prime},1)\neq\emptyset$ and $r_k\leq
t^{\star}_m$. If we denote
$\rho_k=r_k^{\prime\prime}+\frac{2r_k^{\prime\prime}}{\phi(r_k)}$
then we obtain $\rho_k<1$ for $k=1,\ldots,h$. In this case, we put
$E_m^1=\bigcup\limits_{k=1}^{h}(r_k,\rho_k)$. If
$F_m\cap(t_m,t_m^{\star}]=\emptyset$, then we set
$E_m^1=\emptyset$. To estimate the $1$-measure of $E_m^1$ we
use~\eqref{277} and the definition of $F_m$. Thus
\begin{eqnarray*}\mes_1(E_m^1) & < & \sum\limits_{k=1}^{h}(\rho_k-r_k)<\sum\limits_{k=1}^{h}
\frac{4}{\phi(r_k)}
=\sum\limits_{k=1}^{h}\frac{4m^2}{A(r_k)^{1/(Q-1)}} \\ & \leq &
\frac{4m^2}{A(t_m)^{1/(Q-1)}\big(1-\big(\frac{m-1}{m}\big)^{1/(Q-1)}\big)}<\frac{1}{2m}(1-t_m).\end{eqnarray*}
Set $E_m=[t_m,t^{\star}_m]\cup[\overline t_m,1)\cup E_m^1$. Then
$\mes_1(E_m)<\frac{3(1-t_m)}{m}$. Making use of definitions of
sequences $\{t_m\}$ and $\{E_m\}$ we set
$$E=\bigcup\limits_{m\geq 3}\big(E_m\cap[t_m,t_{m+1}]\big).$$
Then, clearly
$$\lim\limits_{m\to\infty}\frac{\mes_1(E\cap[t_m,1))}{1-t_m}=0.$$

Let $\varepsilon\in(0,1/5)$, let $S(w,t)$ be a sphere in $\mathbb
G$, and let $m_1$ be the integer $m$ satisfying the conditions
$(i)-(iii)$ in the proof of Theorem~\ref{44}. Suppose
$s^{\prime}\in[t_{m_1},1)\setminus E$. Then $s^{\prime}$ belongs
to some interval $[t_m,t_{m+1}]$ with $m>m_1$. By the definition
of $E$, we have $s^{\prime}\in(t^{\star}_{m},\overline
t_{m})\setminus E_m^1$. We claim that there exists
$r\in(t_m,\overline t_m)$ such that
$$s^{\prime}=s+\frac{s}{\varepsilon_0(A(s))^{1/(Q-1)}}\quad\text{with}\quad
s=r+\frac{r}{\phi(r)}$$ and $r\notin F_m$. In fact, if we suppose
that $r\leq t_m$ from~\eqref{277} we get
\begin{eqnarray*}
s^{\prime} & = &
r+\frac{r}{\phi(r)}+\frac{r+\frac{r}{\phi(r)}}{\varepsilon_0A(r+\frac{r}{\phi(r)})^{1/(Q-1)}}
\leq t_m+\frac{t_m}{\phi(t_m)}+
\frac{t_m+\frac{t_m}{\phi(t_m)}}{\varepsilon_0A(t_m)^{1/(Q-1)}}
\\ & \leq & t_m+\frac{1-t_m}{8m}+\frac{1+(1-t_m)/8m}{16m}(1-t_m)
 \leq  t_m+\frac{1-t_m}{4m}<t^{\star}_m.
\end{eqnarray*}
It follows that $r>t_m$. If we assume that $r>\overline t_m$, then
$\overline t_m<r<s<s^{\prime}$ and we get a contradiction with
$s^{\prime}\in(t^{\star}_m,\overline t_m)$. To show that $r\notin
F_m$ we argue as in~\eqref{282}. Since, in addition, $F_m\subset
E_m^1\subset E_m$ we obtain inequalities~\eqref{218}
and~\eqref{219} as in the proof of Theorem~\ref{44}. If $m_0$ is
the least positive integer with $\frac{2}{m_0^2}<\varepsilon$, we
can put $D(\varepsilon_0)=t_{m_0}$. From~\eqref{277} it follows
that
$$A(D_{\varepsilon_0})>\Big(\frac{8m_0^3}{1-D_{\varepsilon_0}}\Big)^{Q-1}>(8m_0^4)^{Q-1}>\varepsilon_0^{2-2Q}.$$
The lemma is proved.
\end{proof}

We are at the point to show Theorem~\ref{46}.

{\it Proof of Theorem~\ref{46}.} Let $E\subset(0,1)$ be the set
constructed in Lemma~\ref{127}. We fix
$\varepsilon_0\in\bigl(0,\min\{\frac{1}{5},\frac{1}{8q+9}\}\bigr)$
such that, in addition, $\varepsilon_0$ satisfies~\eqref{278}
and~\eqref{279}. By Lemma~\ref{127}, we find $\kappa\in(0,1)$ such
that, for every $s^{\prime}\in(\kappa,1)\setminus E$, there exists
$s$ with $s^{\prime}=s+\frac{s}{\phi(s)}$ and such that
 the estimates~\eqref{224}, \eqref{223} hold. Fix  such  an $s$. Then we
denote by $f_0$ the restriction
$f\vert_{B(0,s^{\prime}+(1-s^{\prime})/2)}$. Now, for proving of
Theorem~\ref{46},  we repeat the arguments used in the proof of
Theorem~\ref{45}.

\end{document}